\documentclass[11pt]{article}

\usepackage{inputenc}
\usepackage{enumerate}
\usepackage{amssymb}
\usepackage{amsmath}
\usepackage{mathrsfs}
\usepackage{amsthm}
\usepackage{tikz-cd}
\usepackage{authblk}
\usepackage{mathpazo}
\usepackage{abstract}
\usepackage{extarrows}
\usepackage{color}
\usepackage{setspace}
\usepackage[colorlinks,linkcolor = blue,anchorcolor = red,urlcolor = blue,citecolor = blue]{hyperref}
\usepackage{cleveref}
\usepackage{graphicx, subfig}
\usepackage[numbers,sort]{natbib}
\usepackage{setspace}
\usepackage{lipsum}
\usepackage{xpatch}
\xpatchcmd{\author}{\relax#1\relax}{\relax\detokenize{#1}\relax}{}{}

\usepackage{titlesec}
\usepackage{titletoc}
\usepackage[titletoc]{appendix}

\theoremstyle{definition}

\newtheorem{Def}{Definition}[section]

\newtheorem{Rem}[Def]{Remark}
\newtheorem{Cons}[Def]{Construction}
\theoremstyle{plain}
\newtheorem{Thm}[Def]{Theorem}
\newtheorem{Lem}[Def]{Lemma}
\newtheorem{Pro}[Def]{Proposition}
\newtheorem{Cor}[Def]{Corollary}

\newtheorem*{ECBCcon}{The equivariant coarse Baum-Connes Conjecture}
\newtheorem*{ECNcon}{The equivariant coarse Novikov Conjecture}
\newtheorem*{con}{Conjecture}

\usepackage[left=2.6cm, right=2.6cm, top=3cm, bottom=3cm]{geometry}
\usepackage[all]{xy}
\usepackage{bm}

\newcommand\blfootnote[1]{%
\begingroup \renewcommand\thefootnote{}\footnote{#1}%
\addtocounter{footnote}{-1}%
\endgroup }

\def\IN{\mathbb N}\def\IR{\mathbb R}\def\IC{\mathbb C}\def\IQ{\mathbb Q}\def\IZ{\mathbb Z}

\def\A{\mathcal A}\def\C{\mathcal C}\def\B{\mathcal B}\def\E{\mathcal E}\def\K{\mathcal K}\def\F{\mathcal F}\def\P{\mathcal P}\def\L{\mathcal L}\def\R{\mathcal R}\def\H{\mathcal H}\def\S{\mathcal S}\def\N{\mathcal N}\def\U{\mathcal U}\def\V{\mathcal V}
\def\sC{\mathscr C}

\def\MRdnZ{\widetilde{P_{d,n}}(Z)}
\def\MRdnX{\widetilde{P_{d,n}}(X)}

\def\supp{\textup{supp}}
\def\ind{\textup{Ind}}
\def\prop{\textup{Prop}}
\def\diam{\textup{diam}}
\def\Cl{\textup{Cliff}_{\IC}}

\def\ox{\otimes}
\def\wh{\widehat}
\def\wt{\widetilde}
\def\ox{\otimes}
\def\EG{\underline{\E}\Gamma}
\def\Ga{\Gamma}
\def\car{\curvearrowright}

\begin{document}

\title{Hilbert-Hadamard spaces and the equivariant coarse Novikov conjecture
\thanks{
The first two authors are supported by NSFC (No.~12171156), Key Laboratory of MEA, Ministry of Education, and Shanghai Key Laboratory of PMMP, the Science and Technology Commission of Shanghai (No.~22DZ2229014).
The third author is partially supported by
National Key R\&D Program of China 2022YFA100700 and 
NSFC Key Program No.~12231005. 
The fourth author is supported by NSF~2247313, and the Simons Fellows Program. 
}}

\author[a]{Liang Guo}
\author[b]{Qin Wang}
\author[c]{Jianchao Wu}
\author[d]{Guoliang Yu}
\affil[a]{\small\it{Shanghai Institute for Mathematics and Interdisciplinary Sciences, Shanghai, P.~R.~China}}
\affil[b]{\small\it{Research Center for Operator Algebras, School of Mathematical Sciences, East China Normal University, Shanghai, P.~R.~China.}}
\affil[c]{\small\it{Shanghai Center for Mathematical Sciences, Fudan University, Shanghai, P.~R.~China}}
\affil[d]{\small\it{Department of Mathematics, Texas A\&M University, College Station, Texas, United States}}

\maketitle

\begin{abstract}
	The equivariant coarse Novikov conjectures stand among a handful profound $K$-theoretic conjectures in noncommutative geometry. 
	Roughly speaking, these Novikov-type conjectures provide an algorithm to determine when an elliptic operator on a (possibly non-compact) manifold has non-vanishing higher indices, the latter being refinements of the classical Fredholm index that reveal deep insights about the topology and geometry of manifolds. 
Motivated by the quest to verify Novikov-type conjectures for groups of diffeomorphisms, we study in this paper the equivariant coarse Novikov conjectures for spaces that equivariantly and coarsely embed into admissible Hilbert-Hadamard spaces, which are a type of infinite-dimensional nonpositively curved spaces. The paper is split into two parts.

We prove in the first part that for any metric space $X$ with bounded geometry and with a proper isometric action $\alpha$ by a countable discrete group $\Gamma$, if $X$ admits an equivariant coarse embedding into an admissible Hilbert-Hadamard space and $\Gamma$ is torsion-free, then the equivariant coarse strong Novikov conjecture holds rationally for $(X, \Gamma, \alpha)$.

In the second part, we extend the result in the first part by dropping the torsion-free assumption on $\Gamma$. To this end, we introduce, 
for a proper $\Gamma$-space $X$ with equivariant bounded geometry, 
a new Novikov-type conjecture that we call the rational analytic equivariant coarse Novikov conjecture, which generalizes the rational analytic Novikov conjecture and asserts the rational injectivity of a certain assembly map associated with a coarse analog of the classifying space $E\Gamma$. We show that for a proper $\Gamma$-space $X$ with equivariant bounded geometry, if $X$ admits an equivariant coarse embedding into an admissible Hilbert-Hadamard space, then the rational analytic equivariant coarse Novikov conjecture holds for $(X,\Gamma,\alpha)$, i.e., the assembly map is a rational injection.
\end{abstract}

\blfootnote{\emph{Emails: {liangguo@simis.cn (L.~Guo)}, {qwang@math.ecnu.edu.cn (Q.~Wang)}, {jianchao\_wu@fu\-dan.edu.cn (J.~Wu)}, {guoliangyu@tamu.edu (G.~Yu)}}}

\begin{spacing}{0.1}
{\tableofcontents}
\end{spacing}

\section{Introduction}

A major theme in noncommutative geometry is centered around a circle of conjectures that are inspired by the Novikov conjecture and include the (rational) strong Novikov conjecture, the Baum-Connes conjecture, the coarse Baum-Connes conjecture, and several other variants. 
The Novikov conjecture is a major open problem in differential topology that has stimulated the development of a wide range of mathematical areas including geometric topology, higher index theory, and coarse geometry. Among the different approaches to (what is known about) the Novikov conjecture, noncommutative geometry provides some of the strongest results to date, largely through the verification of the (rational) strong Novikov conjecture for many classes of countable groups \cite{Misc1974, Kas88, CM90, CGM93, KS1991, KS2003, HK2001, GHW05, Yu1998, Yu2000, Higson2000, STY2002, mathai03, hankeschick08, MattheyOyono-OyonoPitsch2008Homotopy, FukayaOguni2012, FukayaOguni2015, FukayaOguni2020, GWY2018, GWXY2023}. 
Moreover, 
this approach through the (rational) strong Novikov conjecture also
attains some of the best results on the Gromov-Lawson conjecture in differential geometry, which claims aspherical manifolds do not support Riemannian metrics with everywhere positive scalar curvatures \cite{Rosenberg1983}. 
For further references, see \cite{RosenbergWeinberger1990,  Weinberger1990, FerryWeinberger1991, Valette2002, ChangWeinberger2006, DranishnikovFerryWeinberger2008, MattheyOyono-OyonoPitsch2008Homotopy, KaminkerTang2009, OY2009, RYI2012, RYII2012, CWY2013,WeinbergerXieYu2021, Weinberger2023, GWXY2023, WangXieYuZhu2024}.

A key development in this line of research came with the entry of coarse geometric ideas as a toolbox for attacking the Novikov conjecture. Most notably, the concept of coarse embeddings, introduced by M.~Gromov \cite{Gromov1993} as a natural notion of inclusion in coarse geometry, has played a significant role in the theory for the last quarter century, starting with a result of the last author that countable groups that coarsely embed into a Hilbert space satisfy the coarse Baum-Connes conjecture and the strong Novikov conjecture \cite{Yu2000,STY2002}. 

Vast classes of groups were shown to satisfy these conjectures by verifying their coarse embeddability into Hilbert spaces. 
Attempts to handle wider classes of groups have since led researchers to consider coarse embeddability into other types of ``model'' metric spaces, such as certain Banach spaces, nonpositively curved Riemannian manifolds, etc, leading to deep results
in coarse geometry and the geometry of Banach spaces \cite{KY2006,KY2012,CWY2015}. 

Despite these progresses, for a long time, it was not clear how to deal with groups coarsely embeddable into spaces which are ``nonpositively curved'' and ``infinite-dimensional'' at the same time.  
We point out that the need to consider ``nonpositively curved metric spaces'' (e.g., CAT(0) spaces) as opposed to merely linear spaces is linked to the fact that while there are obstructions to coarse embeddability of groups into Hilbert spaces and various classes of Banach spaces, as of now there is no known examples of metric spaces with bounded geometry that are not coarsely embeddable into CAT(0) spaces \cite[Question~9]{EskenazisNaor2019}.

Motivated by this, 
a recent paper of the last two authors together with S.~Gong (\cite{GWY2018}) introduced the notion of admissible Hilbert-Hadamard spaces, a kind of CAT(0) metric spaces that are ``manifold-like'' but typically ``infinite-dimensional'' (see \Cref{HHspace}), and includes as basic examples separable Hilbert spaces and Hadamard manifolds. It is proved there that any countable group admitting an isometric and metrically proper action on such a space satisfies the rational strong Novikov conjecture. 
This includes geometrically discrete groups of volume-form-preserving\footnote{The volume-form-preserving condition was removed in \cite{GWXY2023}.} diffeomorphisms on smooth closed manifolds. Here the condition of geometric discreteness can be understood as a polar opposite to the condition of being isometric (with regard to some Riemannian metric; see \cite[the discussion below Theorem~1.3]{GWY2018}), the latter being a much better understood scenario in terms of the Novikov conjecture, for the group of isometries on a closed Riemannian manifold is a compact Lie group and thus its countable subgroups embed into Hilbert spaces \cite{GHW05}. 

In order to bridge between these polar opposite cases for the study of the Novikov conjecture on groups of diffeomorphisms, a natural unifying approach is to consider coarse embeddings into not only Hilbert spaces, but admissible Hilbert-Hadamard spaces. 
Indeed, from a coarse-geometric point of view, the condition in \cite{GWY2018} that there is an isometric and metrically proper action of a countable group $\Gamma$ on an admissible Hilbert-Hadamard space $M$ can be reformulated as requiring the existence of a coarse embedding $\Gamma \to M$ that is equivariant with regard to the translation $\Gamma \curvearrowright \Gamma$ and an isometric action $\Gamma \curvearrowright M$.  
This immediately leads to the following question: if we only require coarse embeddability into an admissible Hilbert-Hadamard space without requiring equivariance, what conclusions can we still draw with regard to Novikov-type conjectures? 

In this paper, we provide one answer to this question: under said assumption, the group $\Gamma$ (or more generally, a metric space $X$ with bounded geometry) satisfies the rational analytic coarse Novikov conjecture. 
Actually, we can reach a more general result that unifies both the new, purely coarse geometric, setting and the old setting in \cite{GWY2018}: namely, as we gradually weaken the equivariance condition on coarse embeddings into a Hilbert-Hadamard space, we obtain a range of conclusions that can be collectively called the \emph{rational analytic equivariant coarse Novikov conjecture}. 

To be more precise, let us first pin down the statement of this conjecture and discuss its relationship with other prominent Novikov-type conjectures. For the sake of exposition, we shall start with its better known close relative, the \emph{rational equivariant coarse strong Novikov conjecture}. 
Let $X$ be a metric space with bounded geometry and let $\Ga$ be a countable discrete group. We say $X$ is a proper $\Ga$-space if there is a proper and isometric action $\alpha$ of $\Ga$ on $X$. 
This gives rise to an equivariant coarse index map 
$$\mathrm{Ind}_X^{\Ga}: K_n^{\Ga}(X)\to K_n(C^*(X)^{\Ga}) , \quad \text{ for } n \in \{0,1\} ,$$
that goes from the (topological) $K$-homology of $X$ to the operator $K$-theory of the \emph{equivariant Roe algebra}, the latter encoding information given by the large-scale geometry of $X$ and analytic properties of $\Ga$. 
Using the naturality of this map and the equivariant coarse invariance of the right-hand side, one ``takes an inductive limit'' on the left-hand side across all proper $\Ga$-spaces equivariantly coarsely equivalent to $X$ and arrives at the \emph{equivariant coarse Baum-Connes assembly map}
$$\mu_X^{\Ga}:\lim_{d\to\infty}K_n^{\Ga}(P_d(X))\to K_n(C^*(X)^{\Ga}) , \quad \text{ for } n \in \{0,1\} ,$$
where the left-hand side makes use of the Rips complexes of $X$ as a technical tool to circumvent set-theoretical nuisance in the construction of the inductive limit. 
Despite looking more complicated, the left-hand side (also called the topological side) of the assembly map is easier to compute with the help of classical homological algebraic methods, while the right-hand side (also called the analytic side) is far more difficult to handle. 
There are two special cases: 
\begin{itemize}
	\item If $\Ga\car X$ is cobounded, i.e., the quotient space $X/\Ga$ is bounded (e.g., when $X = \Ga$), then the equivariant Roe algebra $C^*(X)^{\Ga}$ is Morita equivalent to the reduced group $C^*$-algebra $C^*_r\Ga$, the topological side is isomorphic to the representable equivariant $K$-homology of $\EG$, the classifying space for proper $\Ga$-actions, 
	and, in this way, $\mu_X^{\Ga}$ coincides with the assembly map 
	$$\mu^{\Ga}:RK^{\Ga}_n(\EG)\to K_n(C^*_r\Ga),$$
	introduced by P.~Baum, A.~Connes, and N.~Higson in \cite{BCH1994} (see also \cite{BC1988} and \cite[Section 4]{KS2003}).

	\item When $\Ga$ is the trivial group, 
	$\mu_X^{\Ga}$ matches the coarse assembly map (\cite{HR1993, YuCBC95})
	\[
	\mu_{X}:KX_n(X)\to K_n(C^*(X))
	\]
	from the coarse $K$-homology group $KX_*(X)$ to the $K$-theory of the Roe algebra of $X$.  
\end{itemize}
The injectivity of the assembly map $\mu_X^{\Ga}$ is known as the \emph{equivariant coarse strong Novikov conjecture} for the $\Ga$-space $X$. 
The adjective ``rational'' is attached at beginning if $\mu_X^{\Ga}$ is only required to be injective after tensoring both sides by $\mathbb{Q}$ (which allows us to ignore torsion on both sides). 
Even in its rational form, this conjecture leads to some far-reaching consequences: the cobounded case alone implies both the (classical) Novikov conjecture and the Gromov-Lawson conjecture, while in general, rational injectivity of $\mu_X^{\Ga}$ for a smooth manifold $X$ with a free and proper $\Ga$-action obstructs the existence of Riemannian metrics on $X / \Ga$ with uniformly positive scalar curvature. 
We refer the reader to \cite[Chapter~7]{HIT2020} for more details. 

Notice that in order to verify the Novikov conjecture and the Gromov-Lawson conjecture, it suffices to show rational injectivity of $\mu_X^{\Ga}$ on a ``smaller'' domain, namely $RK^{\Ga}_*(\E\Ga)$, the representable equivariant $K$-homology of $\E\Ga$, the classifying space for \emph{free} and proper $\Ga$-actions (that is, the ``original'' classifying space, or a $K(\Ga,1)$ space; see \cite[Section 1.B]{Hatcher}; note that $RK^{\Ga}_*(\E\Ga) \cong RK_*(B\Ga)$, where $B\Ga = \E\Ga / \Ga$). 
More precisely, composing the Baum-Connes assembly map $\mu^{\Ga}$ with the natural map $\pi_* : RK^{\Ga}_*(\E\Ga) \to RK^{\Ga}_*(\EG)$ (which is rationally injective by \cite[Section~7]{BCH1994}), we obtain the so-called \emph{Miščenko-Kasparov assembly map}
$$\nu^{\Ga}:RK^{\Ga}_n(\E\Ga)  \xrightarrow{\pi_*} RK^{\Ga}_n(\EG) \xrightarrow{\mu^{\Ga}} K_n(C^*_r\Ga)$$
whose rational injectivity is known to imply the Novikov conjecture and the Gromov-Lawson conjecture \cite{Lusztig1972,Misc1974,Kas83,Rosenberg1983}. 

We generalize this construction beyond the cobounded case in Section~\ref{sec: The Milnor-Rips complex} and introduce the \emph{equivariant coarse Miščenko-Kasparov assembly map} $\nu_X^{\Ga}$, whose domain is now the ``inductive limit'' of all \emph{free} and proper $\Ga$-spaces equivariantly coarsely equivalent to $X$, and it factors as a composition of $\mu_X^{\Ga}$ and a natural map that ``forgets'' the freeness requirement. 
There is however one subtlety: 
in order to avoid certain pathological behaviors, 
we need to assume a condition that we call \emph{equivariant bounded geometry} for the action $\Ga \curvearrowright X$ (see \Cref{ebg}), which, for uniformly locally finite $X$, amounts to the existence of a uniform bound on the sizes of stabilizers $\Ga_x$ for all $x \in X$; equivalently, $X$ is equivariantly coarsely equivalent to a free and proper $\Ga$-space that is uniformly locally finite. 
Now for a $\Ga$-space $X$ with equivariant bounded geometry, in order to rigorously construct the inductive limit, we replace $X$ by a free, proper, uniformly locally finite $\Ga$-space and make use of the \emph{Milnor-Rips complexes} introduced in \cite{Yu1995}, denoted by $\wt{P_{d,k}}(X)$ (see Definition \ref{Milnor Rips}). Hence the equivariant coarse Miščenko-Kasparov assembly map takes the following concrete form: 
\begin{equation}\label{MK assembly map}\nu_X^{\Ga}:\lim_{d,k\to\infty}K^{\Ga}_n(\wt{P_{d,k}}(X))\xrightarrow{\pi_*}\lim_{d\to\infty}K^{\Ga}_n(P_d(X))\xrightarrow{\mu_X^{\Ga}} K_n(C^*(X)^{\Ga}) ,
\end{equation}
where $\pi_*$ is an isomorphism if $\Ga$ is torsion-free, and rationally injective in general (see \Cref{inj between classifying}). 
Hence this assembly map agrees with the Miščenko-Kasparov assembly map in the cobounded case, and also with the equivariant coarse Baum-Connes assembly map when $\Ga$ is torsion-free. 
To distinguish it from the case of $\mu_X^{\Ga}$ above, we shall term the (rational) injectivity of $\nu_X^{\Ga}$ the \emph{(rational) analytic equivariant coarse Novikov conjecture}. This is weaker than the (rational) equivariant coarse strong Novikov conjecture, but still yields the same applications regarding the Novikov conjecture and positive scalar curvature as above.

The main result of this paper is as follows.

\begin{Thm}\label{main result with torsion}
	Let $\Gamma$ be a countable discrete group and let $X$ be a proper $\Gamma$-space (with $\alpha$ denoting the action).
    If $X$ admits a $\Ga$-equivariant coarse embedding into an admissible Hilbert-Hadamard space and $(X, \Ga, \alpha)$ has equivariant bounded geometry, then the rational analytic equivariant coarse Novikov conjecture holds for $(X, \Ga, \alpha)$, i.e., the rationalized equivariant coarse Miščenko-Kasparov assembly map
	$$\nu_X^{\Ga}\ox\operatorname{id}_\IQ:\lim_{d,k\to\infty}K^{\Ga}_*(\wt{P_{d,k}}(X))\ox\IQ\xrightarrow{\pi_*}\lim_{d\to\infty}K^{\Ga}_*(P_d(X))\ox\IQ\xrightarrow{\mu_X^{\Ga}} K_*(C^*(X)^{\Ga})\ox\IQ$$
	is injective.
\end{Thm}

When $\alpha$ is cobounded, Theorem~\ref{main result with torsion} recovers \cite[Theorem~1.1]{GWY2018} for $\Ga$ (and, in particular, implies the classical Novikov conjecture for such groups). 
On the other hand, the ability to deal with the non-cobounded situation (and for that matter, the purely ``coarse case'', where $X$ is unbounded but $\Ga$ is trivial) enables
Theorem~\ref{main result with torsion} to give wider applications to the positive scalar curvature problem for noncompact manifolds. 

\begin{Cor}\label{cor:psc}
	Let $N$ be a complete Riemannian manifold. If its universal cover $\widetilde{N}$ is uniformly contractible and 
	admits a coarse embedding into an admissible Hilbert-Hadamard space, 
	then $N$ does not have uniformly positive scalar curvature.
\end{Cor}


Let us give a few remarks on our proof of \Cref{main result with torsion}. On the one hand, the broad strategy of our proof is inspired by ideas that have appeared in \cite{GWY2018}, in particular, in the use of a noncommutative ``coefficient'' $C^*$-algebra $\A(M)$ associated with a Hilbert-Hadamard space $M$ and the adoption of a deformation trick to complement the classical Dirac-dual-Dirac method\footnote{Note that \cite{GWY2018} makes use of $KK$-theory with real coefficients (\cite{AAS2020}) to deal with groups with torsion, we shall use a localization algebraic technique in its place.}. 
On the other hand, however, due to the new challenges presented by the equivariant coarse setting, 
our implementation of the strategy differs significantly from \cite{GWY2018}. 
The most evident among the differences is the fact that we employ, in place of equivariant $KK$-theory, the toolbox of equivariant localization algebras originally introduced by the last author. 

In order to overcome difficulties that arise in the equivariant coarse setting, we develop a number of new constructions and techniques. For example, we introduce equivariant twisted Roe algebras and equivariant twisted localization algebras with coefficients in the aforementioned $C^*$-algebra $\A(M)$. 
Similar algebras have been constructed in \cite{FW2016} in the special case when $M$ is a Hilbert space, but even in this special case, our construction is much simpler and more intrinsic, since we do not (and, for more general Hilbert-Hadamard spaces, probably cannot; see \cite[Remark~7.7]{GWY2018}) rely on writing $\A(M)$ as an inductive limit of $C^*$-algebras $\A(M_n)$ associated with finite-dimensional convex subspaces $M_n$. 
We also bridge a subtle gap left in \cite{FW2016} in the equivariant cutting-and-pasting argument for this type of algebras (see Remark \ref{difference with FW}). Moreover, in order to implement the deformation trick, we introduce the notion of annihilation ideals in (equivariant twisted) localization algebras; by modding out these ideals, we are able to upgrade the functoriality of the $K$-theory of these algebras from using Lipschitz maps to using more general continuous maps (see Appendix \ref{Appendix A}). 

More precisely, the proof centers around the following commuting diagram:
\begin{equation}\begin{tikzcd}\label{main diag}
K^{\Ga}_{*+1}(\MRdnX) \arrow[drr, "\nu_{X,d,n}^\Ga",bend left=10]  \arrow[d, "\cong"] &&  \\
K_{*+1}(C^*_L(\MRdnX)^{\Ga}) \arrow[r, "\pi_*", "(2)"']  \arrow[d, "(\widetilde{\beta_L})_*", "(5)"'] & K_{*+1}(C^*_L(P_d(X))^{\Ga}) \arrow[r, "ev_*","(1)"'] \arrow[d, "(\beta_L)_*", "(4)"'] & K_{*+1}(C^*(P_d(X))^{\Ga}) \arrow[d, "\beta_*", "(3)"']   \\
K_*(C^*_L(\MRdnX,\A_{[0,1]}(M))^{\Ga}) \arrow[r, "\pi_*","(10)"'] \arrow[d, "(\widetilde{ev_1})_*", "(8)"']  & K_*(C^*_L(P_d(X),\A_{[0,1]}(M))^{\Ga})\arrow[r, "ev_*","(9)"'] \arrow[d, "(ev_1)_*","(7)"']& K_*(C^*(P_d(X),\A_{[0,1]}(M))^{\Ga}) \arrow[d, "(ev_1)_*","(6)"']\\
K_*(C^*_L(\MRdnX,\A(M^{[0,1]}))^{\Ga,(1)}) \arrow[r, "\pi_*","(12)"']& K_*(C^*_L(P_d(X),\A(M^{[0,1]}))^{\Ga,(1)})\arrow[r, "ev_*","(11)"'] & K_*(C^*(P_d(X),\A(M^{[0,1]}))^{\Ga,(1)})
\end{tikzcd}\end{equation}
By taking inductive limits with regard to $d$ and $n$, the top horizontal map becomes the Miščenko-Kasparov assembly map $\nu_X^{\Ga}$, which we need to show is rationally injective. To do this, 
we shall, after taking inductive limits of all items, trace along the leftmost column and then the bottom row, and show their composition is rationally injective. 
Since the entire proof is lengthy, for the convenience of the reader, we shall split the proof into two parts: 
\begin{itemize}
    \item Part~I includes Sections~\ref{sec: Hilbert-Hadamard spaces and their K-theory}-\ref{sec: The twisted assembly map and proof of Theorem 1.3} and shows that in diagram~\eqref{main diag}, the maps in the left column and map~(11) are rationally injective (after taking inductive limits). 
    \item Part~II includes Section~\ref{sec: The Milnor-Rips complex}-\ref{sec: Proof of the main theorem} and shows that map~(12) is rationally injective. 
\end{itemize}

In the special case when $\Ga$ is torsion-free, 
since the Milnor-Rips complexes $\MRdnX$ agree with the usual Rips complexes $P_d(X)$ (for large enough $n$), Part~I alone is sufficient to prove Theorem \ref{main result with torsion} with a strengthened conclusion: the \emph{equivariant strong coarse Novikov conjecture} holds for $(X,\Ga,\alpha)$. So, as a half-range target, we shall show the following theorem as a direct consequence.

\begin{Thm}\label{main result torsion-free}
	Let $\Gamma$ be a countable discrete group and let $X$ be a proper $\Gamma$-space (with bounded geometry, and with $\alpha$ denoting the action). 
    If $X$ admits a $\Ga$-equivariant coarse embedding into an admissible Hilbert-Hadamard space and $\Ga$ is torsion-free, 
	then the equivariant coarse assembly map is rationally injective, i.e.,
$$\mu_{X}^{\Ga}\ox_{\IZ}\operatorname{id}_{\IQ}:\lim_{d\to\infty}K^{\Ga}_*(P_d(X))\ox_{\IZ}\IQ\to K_*(C^*(X)^{\Ga})\ox_{\IZ}\IQ$$
is injective.
\end{Thm}  

We give a bit more details regarding Part~I and  
our proof of Theorem~\ref{main result torsion-free}. We start by following the ``dual-Dirac'' half of the Dirac-dual-Dirac method, 
one of the most powerful methods in dealing with Novikov-type conjectures, first introduced by G.~Kasparov (\cite{Kas88}) in the language of $KK$-theory and further developed by several authors (see, for example, \cite{KS1991, Tu1999, KS2003, GHT2000, HK2001, CE2001}). 
A coarse version of this method, based on twisted localization algebras, was developed by the last author (\cite{Yu2000}) and subsequently applied in a range of works, including some in the \emph{equivariant} coarse setting (see \cite{FW2016,FWY2020}). 
We refer the reader to \cite{GW2022} for a comparison of this (equivariant) coarse version with the original $KK$-theoretic version. 

Following the equivariant coarse version of the Dirac-dual-Dirac method, we divide the proof into two steps: (1) reduce the problem to a "twisted" version of the equivariant coarse Baum-Connes conjecture with "nice" coefficients; (2) show the twisted assembly map (map~(11) in diagram~\eqref{main diag}, after taking inductive limits) is indeed an isomorphism. 
By our assumption in Theorem~\ref{main result with torsion}, $X$ equivariantly and coarsely embeds into an admissible Hilbert-Hadamard space $M$ with a proper isometric $\Ga$-action, which gives rise to the $C^*$-algebra $\A(M)$ with a proper $\Ga$-action (see \cite[Section~5]{GWY2018}). Unfortunately, we are not able to compute the $K$-theory of $\A(M)$, for which reason the usual cutting-and-pasting argument fails for step~(1). To circumvent this problem, we adapt an idea in \cite{GWY2018} and introduce a deformation trick for twisted localization algebras (see Section~5) to "trivialize" the $\Ga$-action on the twisted algebras. Then, step~(1) can be realized by using a K\"unneth formula for the $K$-theory of twisted localization algebras (see Proposition \ref{left column}).

In step~(2), we develop a cutting-and-pasting argument by combining the ideas in \cite{Yu2000} and \cite{GHT2000}. 
The condition of equivariant bounded geometry will be used here, in addition to its use in the construction of the Milnor-Rips complexes $\MRdnX$. 
We shall first "cut" the twisted algebras based on $\A(M)$ into "small" algebras such that the assembly map for these "small" algebras can be reduced to the coarse Baum-Connes conjecture for a "transverse" of the action. Then the twisted assembly map is proved to be an isomorphism via a Mayer-Vietoris argument to "paste" all these "small" isomorphisms. 
A dimension control extracted from the equivariant bounded geometry assumption on $X$ will guarantee that the Mayer-Vietoris argument will end in finitely many steps.

As for Part II, it remains to show the map (12) in diagram \eqref{main diag} is a rational injection. Note that an analogous statement was shown in the $\Ga$-compact case in \cite{GWY2018}, using the machinery of $KK$-theory with real coefficients developed in \cite{AAS2020}. In this paper, we shall only use the language of localization algebras to do this. 
To this end, we devise an analogue of $KK$-products in the framework of localization algebras (see Section \ref{KK-product for localization alg}). 
With this tool at our disposal, we adapt an elegant idea from \cite{AAS2020}: Considering the Cartesian product of $P_d(X)$ with a compact $\Gamma$-space $\Delta$ with Property TAF (indicating ``torsion (elements) act freely'') \cite[Definition 5.2]{AAS2020}, we see that $P_d(X)\times\Delta$ retains its equivariant coarse equivalence with $X$ due to the compactness of $\Delta$, while the diagonal $\Gamma$-action becomes not only proper, but also free, whence by the universal property of the Milnor-Rips complexes (see Section \ref{universal property for MR complexes}), there exists a $\Gamma$-equivariant continuous map from $P_d(X)\times\Delta$ to $\widetilde{P_{d',n}}(X)$ for sufficiently large $d'$ and $n$, enabling the construction of a partial inverse to $\pi_*$ in diagram~\eqref{main diag}.

Actually, we can do better: we introduce a concrete model of a $\Ga$-space $\Omega_X$ with Property TAF, which consists of all linear orders on $X$ (see Section \ref{sec: A Concrete model for space with Property TAF}). Using this model, an ``inverse map'' from $P_d(X)\times\Delta$ to $\MRdnX$ (with the same $d$) can be constructed concretely, simply by ``ordering" vertices in $P_d(X)$ according to a linear order in $\Omega_X$ (see formula~\eqref{eq:ordering_map} in the proof of Lemma~\ref{pi is injective}).
We then proceed, by using a K\"unneth formula and a canonical $\Gamma$-invariant trace on $C(\Delta)$, to construct the desired partial inverse to $\pi_*$ (see Lemma~\ref{pi is injective}). 
Along the way of showing the constructed map is indeed a partial inverse, we settle in the positive a fundamental question about $\A(M)$, namely whether it satisfies the UCT (see Section~\ref{Sec: Approximating AofM}). In fact, we show $\A(M)$ is a direct limit of type I $C^*$-algebras.  

This paper is organized as follows. For Part I, Section \ref{sec: Hilbert-Hadamard spaces and their K-theory} briefly recalls the definition of Hilbert-Hadamard spaces and a $C^*$-algebra associated with a Hilbert-Hadamard space. In Section \ref{sec: The equivariant coarse Novikov conjecture}, we briefly recall the equivariant Roe algebra and Yu's localization algebra and state the equivariant coarse strong Novikov conjecture. In Section \ref{sec: Twisted algebras and the Bott homomorphisms}, we construct several twisted algebras associated with the coarse embedding. We also introduce a version of twisted algebras to formulate a deformation trick with a larger coefficient algebra. Then we construct Bott maps on both sides of the assembly map and build the commutative diagram \eqref{com diag}. In Section \ref{sec: A deformation trick}, we formulate a deformation trick and prove the Bott map for the topological side is a rational injection. This requires a \emph{K\"unneth formula} for the $K$-theory of twisted localization algebras. To complete a proof, we discuss a more flexible version of twisted localization algebras and show the K\"unneth formula in Appendix \ref{Appendix A}. In Section \ref{sec: The twisted assembly map and proof of Theorem 1.3}, we show the twisted assembly map is an isomorphism and complete the proof of Theorem \ref{main result torsion-free}.

For Part II, in Section \ref{sec: The Milnor-Rips complex}, we briefly recall the universal property for the equivariant coarse Baum-Connes assembly map. Then we introduce the Milnor-Rips complexes for a proper $\Ga$-space. The inductive limit $\lim_{d,n\to\infty}\MRdnX$ forms a concrete model for the ``universal" group of all proper and free $\Ga$-spaces which are equivariantly coarsely equivalent to $X$. In Section \ref{sec: A Concrete model for space with Property TAF}, we construct a concrete model of a $\Ga$-space with Property TAF by using all linear orders on $X$ for any proper and free $\Ga$-space $X$. In Section \ref{Sec: Approximating AofM}, we prove that the $C^*$-algebra $\A(M)$ associated with a separable Hilbert-Hadamard space is a limit of Type I $C^*$-algebra, thus in the Bootstrap class. This result will be used to show the K\"unneth formula in the next section. In Section \ref{sec: KK-theoretic constructions for twisted localization algebras}, we introduce two constructions for twisted localization algebras which are inspired by Kasparov's equivariant $KK$-theory, including a $KK$ product and a K\"unneth formula in the language of the twisted localization algebras. In Section \ref{sec: Proof of the main theorem}, we construct the ``inverse map" to the map (12) and complete the proof of Theorem \ref{main result with torsion}.

\section{Hilbert-Hadamard spaces and their $K$-theory}\label{sec: Hilbert-Hadamard spaces and their K-theory}

In this section, we recall a class of metric spaces called Hilbert-Hadamard spaces introduced in \cite{GWY2018}. 
Intuitively speaking, these metric spaces may be understood as (typically infinite-dimensional) nonpositively curved manifolds. 
Despite their lack of local compactness, a method was developed in \cite{GWY2018} to study the $K$-theoretic properties of such a space $M$ by introducing a noncommutative $C^*$-algebra $\A(M)$.

\subsection{Hilbert-Hadamard spaces}\label{sec:Hilbert-Hadamard-spaces}

Let $(M,d)$ be a metric space. A continuous map from the interval $[0,1]$ to $M$ is called a \emph{path}. The \emph{length} of $\gamma$, denoted by $\|\gamma\|$, is defined to be
$$\|\gamma\|=\sup\sum_{i=1}^nd(\gamma(t_{i-1},t_i))$$
where the supremum is taken over the set of all partitions $0=t_0< \cdots < t_n = 1$ of the interval $[0,1]$, with an arbitrary $n\in\mathbb{N}$. We always use the symbol $\gamma_t=\gamma(t)$.

A path $\gamma:[0,1]\to M$ is called a \emph{geodesic} if $d(\gamma_s,\gamma_t)=d(\gamma_0,\gamma_1)|s-t|$ for every $s,t\in [0,1]$. We denoted the image of a geodesic joining $x$ and $y$ by $[x,y]$.  For a point $z\in[x,y]$, we often use the notation $z=(1-t)x+ty$, where $t=\frac{d(x,z)}{d(x,y)}$, and say that $z$ is a convex combination of $x$ and $y$. If $x_0,x_1\in M$ and $t\in [0,1]$, then the point $(1-t)x_0+tx_1$ on a geodesic $[x_0,x_1]$ is usually denoted by $x_t$. The point $x_{\frac 12}$ is called a \emph{midpoint} of $x_0$ and $x_1$, which means that $d(x_0,x_{\frac12})=d(x_{\frac12},x_1)=\frac12d(x_0,x_1)$.

Let $(M,d)$ be a geodesic metric space, i.e., any two points in $M$ can be connected by a geodesic path. A \emph{geodesic triangle} with vertices $p,q,r$ consists of three geodesics $[p,q], [p,r], [q,r]$, denoted by $\triangle(p,q,r)$. Given such a geodesic triangle in $M$, there exists a \emph{comparison triangle} in $\mathbb{R}^2$, that is, three line segments $[\bar p,\bar q]$, $[\bar p,\bar r]$ and $[\bar q,\bar r]$ in $\mathbb{R}^2$, such that
$$d(p,q)=\|\bar p-\bar q\|,\quad d(p,r)=\|\bar q-\bar r\|,\quad d(q,r)=\|\bar q-\bar r\|$$
Notice that the comparison triangle is unique up to isometries on $\IR^2$ and we shall denote it $\triangle(\bar p,\bar q,\bar r)$. There is a natural correspondence between $\triangle(p,q,r)$ and $\triangle(\bar p,\bar q,\bar r)$. Assume that $x=tp+(1-t)q$ for some $t\in[0,1]$. Then the \emph{comparison point} for $x$ in the comparison triangle is the point $\bar x=t\bar p+(1-t)\bar q$.

\begin{Def}
Let $(M,d)$ be a geodesic space. $M$ is said to be a \emph{CAT(0) space} if for every geodesic triangle with $p,q,r\in X$ and $x\in[p,r],y\in[p,q]$, we have
$$d(x,y)\leq\|\bar x-\bar y\|,$$
where $\bar x$ and $\bar y$ are the corresponding comparison points in the comparison triangle $\triangle(\bar p,\bar q,\bar r)$. A complete CAT(0) space is also called a \emph{Hadamard space}.
\end{Def}

\begin{Rem}
If $(M,d)$ is a CAT(0) space, then every geodesic triangle with vertices $p,q,r\in X$ and $x_t\in[q,r]$, we have
$$d(p,x_t)\leq\|\bar p-\bar{x_t}\|,$$
where $x_t=tq+(1-t)r$ and $t\in[0,1]$. By an elementary calculation, we have
$$d(p,x_t)^2\leq \|\bar p,\bar x_t\|^2=td(p,r)^2+(1-t)d(p,q)^2-t(1-t)d(q,r)^2$$
Take $m=x_{\frac12}$ and we have
\begin{equation}\label{CATEQ}4d(p,m)^2+d(q,r)^2\leq 2d(p,r)^2+2d(p,q)^2\end{equation}
A geodesic space is CAT(0) if and only if the inequality \eqref{CATEQ} holds for any $p,q,r\in X$ and the midpoint $m$ of $q,r$ (cf. \cite{Hadamard1999}). This is called the \emph{semi-parallelogram law}.
\end{Rem}

Typical examples for Hadamard spaces include Hilbert spaces, $\IR$-trees, and Hadamard manifolds (i.e., complete, connected, and simply connected Riemannian manifolds).

\begin{Def}
Let $(M,d)$ be a geodesic space. We define the \emph{angle} between two geodesics $\alpha:[0,1]\to X$ and $\beta:[0,1]\to X$ with $\alpha_0=\beta_0$ by
$$\angle(\alpha,\beta):=\lim_{s,t\to 0+}\angle(\bar{\alpha_t},\bar\alpha_0,\bar{\beta_s})=\lim_{s,t\to 0+}\frac{\langle\bar{\alpha_t}-\bar\alpha_0,\bar{\beta_s}-\bar\beta_0\rangle}{\|\bar{\alpha_t}-\bar\alpha_0\|\cdot\|\bar{\beta_s}-\bar\beta_0\|},$$
where $\bar{\alpha_t},\bar{\beta_s}\in\mathbb{R}^2$ are the comparison points in the comparison triangle.
\end{Def}

For a given point $p\in M$, denote by $S'_p$ the metric space consisting of all equivalence classes of geodesic segments emanating from $p$, where two geodesic segments are identified if they have zero angles. The distance function on $S'_p$ is defined by
$$d_S([\alpha],[\beta])=\angle(\alpha,\beta).$$
It is well-known that $d_S$ is a well-defined metric function (see \cite[Lemma 1.1.4]{Hadamard1999} for the proof of $d_S$ satisfying triangle inequality). Let $S_p$ be the completion of $S'_p$ under the metric function. To simplify the notation, we denote the angle between the geodesics $(p,q)$ and $(p,r)$ by $\angle(q,p,r)$.

\begin{Def}
The \emph{tangent cone} $T_pM$ at a point $p\in M$ is then defined to be 
$$T_pM=S_p\times[0,\infty)/S_p\times\{0\}$$
For two points $[\alpha,t],[\beta,s]\in T_pM$, where $[\alpha],[\beta]\in S_p$, $t,s\in[0,\infty)$ and the square brackets means taking equivalence class, the metric is given by
$$d([\alpha,t],[\beta,s])=(t^2+s^2-2st\cos(d_S([\alpha],[\beta])))^{\frac12}$$
\end{Def}

Now, we are ready to introduce the notion of Hilbert-Hadamard spaces:

\begin{Def}\label{HHspace}
A \emph{Hilbert-Hadamard space} is a complete geodesic CAT(0) space all of whose tangent cones are isometrically embeddable into Hilbert spaces.
\end{Def}

For any point $x$ in a Hilbert-Hadamard space $M$, we define the \emph{tangent Hilbert space} $\H_xM$ to be the smallest Hilbert space that contains $T_xM$ with the tip of the cone being the original. See \cite[Construction 2.14]{GWY2018} for concrete construction.

For each $x\in M$, define the logarithm map $\log_x:M\to T_xM$ by
$$\log_x(y)=\left[[x,y],d(x,y)\right],$$
for each $y\in M$. By CAT(0) condition, one can easily check that $\log_x$ is non-expansive, i.e.,
$$d(\log_x(y),\log_x(z))\leq d(y,z)$$
for each $y,z\in M$. If $y,z$ are on the same geodesic emanating from $x_0$, then $d(\log_x(y),\log_x(z))=d(y,z)$. If $M$ is a Hadamard manifold, then the logarithm map as defined above is the inverse of the exponential map in Riemannian geometry. In this case, $\log_x$ is a diffeomorphism for each $x\in M$ by the famous Cartan-Hadamard theorem (see \cite[Theorem 12.8]{GTM176}).

\begin{Def}
A separable Hilbert-Hadamard space $M$ is said to be \emph{admissible} if there exists an increasing sequence of subspaces $\{M_n\}_{n\in\IN}$ of $M$ such that\begin{itemize}
\item[(1)] each $M_n$ is a convex closed subspace of $M$ and isometric to a finite-dimensional Hadamard manifold;
\item[(2)] $\bigcup_{n\in\IN}M_n$ is dense in $M$.
\end{itemize}\end{Def}

The following construction is called the \emph{continuum product} of a metric space which will play an important role in the following proof. Let $(X,d_X)$ be a metric space, and $(Y,\mu)$ be a compact probability space. Fix a base point $x_0\in X$, define $L^2(Y,\mu, X)$ to be the set of all equivalence classes of measurable functions $\xi: Y\to X$ satisfying
$$\int_Yd_X(\xi(y),x_0)d\mu(y)<\infty.$$
Two functions $\xi$ and $\eta$ are equivalent iff $\mu(\{y\in Y\mid \xi(y)\ne\eta(y)\})=0$. For example, when $Y$ is a finite space and $\mu$ is the counting measure on $Y$, then $L^2(Y,\mu,X)$ is exactly the Cartesian product $X^{\mu(Y)}$ equipped with the $\ell^2$-metric. The following proposition is proved in \cite[Proposition 3.13]{GWY2018}:

\begin{Pro}[\cite{GWY2018}]
Let $(M,m_0)$ be an admissible Hilbert-Hadamard space with base point $m_0$ and $(Y,\mu)$ be a separable finite measure space. Then $L^2(Y,\mu, M)$ is still an admissible Hilbert-Hadamard space.
\end{Pro}

One typical example is when $Y=[0,1]$ and $\mu$ is the Lebesgue measure on $[0,1]$, we shall denote $L^2([0,1],\mu,M)$ by $M^{[0,1]}$. We will view $M$ as a subspace of $M^{[0,1]}$ under the constant embedding, i.e., $m\mapsto \xi_m\in M^{[0,1]}$ for each $m\in M$, where $\xi_m(t)=m$ for each $t\in[0,1]$.

\subsection{A $C^*$-algebra associated with $M$ and its $K$-theory}

Throughout the rest of the paper, we shall always assume $M$ to be an admissible Hilbert-Hadamard space, unless otherwise specified. In this subsection, we sketch the construction and basic properties of a $C^*$-algebra $\A(M)$ associated to $M$. The reader is referred to \cite[Sections~5-7]{GWY2018} for a detailed treatment.

Given a Hilbert space $\H$, denote by $\Cl(\H)$ the complex Clifford algebra associated to $\H$ (for an introduction to Clifford algebras associated to infinite-dimensional Hilbert spaces, see, for example, \cite[Section~3.1]{GLWZ2022}). 
\begin{Def}\label{def:Pi-alg}
Denote $\IR_+=[0,\infty)$, Define
$$\Pi(M)=\prod_{(x,t)\in M\times\IR_+}\Cl(\H_xM\oplus t\IR)$$
where
$$t\IR=\left\{\begin{aligned}&\IR,&&t\ne0\\&\{0\},&&t=0\end{aligned}\right.$$
where $\IR$ carries the canonical inner product independent of $t$. Define the $C^*$-algebra
$$\Pi_b(M)=\{\sigma\in\Pi(M)\mid \sup\|\sigma(x,t)\|<\infty\}$$
equipped with pointwise algebraic operations and the uniform norm.
\end{Def}

For any $\sigma\in\Pi_b(M)$, define the support of $M$ to be
$$\supp(\sigma)=\overline{\{(x,t)\in M\times\IR_+\mid \sigma(x,t)\ne 0\}}.$$

\begin{Def}\label{def:Clifford_generator}
For any $x_0\in M$, we define the \emph{Clifford generator} $C_{x_0}\in\Pi(M)$ by
$$C_{x_0}(x,t)=(-\log_x(x_0),t)\in T_xM\times t\IR\subseteq \Cl(\H_xM\oplus t\IR)$$
for each $x\in M$ and $t\in\IR_+$.
\end{Def}

We denote $C_0(\IR)_{ev}$ (or $C_0(\IR)_{odd}$, respectively) the subset of $C_0(\IR)$ of all even (or odd, respectively) functions. For any Hilbert space $\H$, the functional calculus of $\Cl(\H)$ is defined in the following way: for any $f\in C_0(\IR)_{ev}$ and $v\in\H\subset\Cl(\H)$, we define $f(v)=f(\|v\|)\in\IC\subset\Cl(\H)$ and for $g\in C_0(\IR)_{odd}$, we define the functional calculus by
$$g(v)=\left\{\begin{aligned}&0&&,v=0\\&g(\|v\|)\frac{v}{\|v\|}&&,v\ne0,\end{aligned}\right.$$
where $g(v)$ is viewed as an element of $\Cl(\H)$. We will denote $C_0(\IR)$ by $\S$ in the following.

\begin{Def}\label{Bott}
The \emph{Bott map} $\beta_x:\S\to\Pi_b(M)$ is defined by
$$(\beta_{x_0}(f))(x,t)=f(C_{x_0}(x,t)).$$
The $C^*$-algebra $\A(M)$ is defined to be the $C^*$-subalgebra of $\Pi_b(M)$ generated by
$$\{\beta_{x_0}(f)\mid x_0\in M,f\in\S\}.$$
More generally, for any subset $N\subseteq M$, we define $\A(M,N)$ to be the $C^*$-subalgebra of $\A(M)$  generated by
$$\{\beta_{x_0}(f)\mid x_0\in N,f\in\S\}.$$
\end{Def}

When $M$ is a Hilbert space and $V\subseteq M$ is a dense Euclidean subspace of $M$, then the $C^*$-algebra $\A(M)$ coincide with $\S\C(V)$ introduced by N.~Higson, G.~Kasparov and J.~Trout and used in \cite{Yu2000}, see \cite[Remark 7.7]{GWY2018}. The next theorem summarizes some results we will need in the following of this paper; one can find a proof in Section 5 and Section 7 of \cite{GWY2018}:

\begin{Thm}\label{Bott-inj}
Let $M$ be an admissible Hilbert-Hadamard space. Then\begin{itemize}
\item[(1)] for any $x_0,x_1\in M$, $\beta_{x_0}$ and $\beta_{x_1}$ are homotopic to each other;
\item[(2)] for any $x_0\in M$, the Bott map induces an injection on $K$-theory, i.e.,
$$(\beta_{x_0})_*:K_*(\S)\to K_*(\A(M))$$
is an injection;
\item[(3)] when $M$ is a finite-dimensional Hadamard manifold or a Hilbert space, for any $x_0\in M$, the Bott map induces an isomorphism on $K$-theory, i.e.,
$$(\beta_{x_0})_*:K_*(\S)\to K_*(\A(M))$$
is an isomorphism.\qed
\end{itemize}\end{Thm}

\begin{Rem}
Denote $\A_{ev}(M)$ the $C^*$-subalgebra of $\A(M)$ generated by $\{\beta_{x_0}(f)\mid f\in C_0(\IR)_{ev}\}$. Any element in $\A_{ev}(M)$ is a complex-valued function on $M\times\IR_+$. Thus the Gelfand spectrum of $\A_{ev}(M)$ is $M\times\IR_+$ equipped with the weakest topology such that all functions in $\A_{ev}(M)$ are continuous. The topological base of $M\times\IR_+$ is given by the preimage of $B_{\IC}(z,\varepsilon)$ of $\beta_{x_0}(g)$ for all $x_0\in M$ and $g\in C_0(\IR)_{ev}$, for example, the following family gives a concrete example for a topological base
$$\B=\{B(x_0,t_0,\varepsilon)\mid x_0\in M, t_0,\varepsilon\in\IR_+\},$$
where $B(x_0,t_0,\varepsilon)=\{(x,t)\in M\times\IR_+\mid t_0-\varepsilon<d(x-x_0)^2+t^2<t_0+\varepsilon\}$.
\end{Rem}

At the last of this section, we shall discuss the $\Ga$-action on $\A(M)$ and the continuum product $M^{[0,1]}$. We will recall a deformation trick introduced in \cite[Proposition 3.18]{GWY2018}.

Let $\Ga$ be a countable discrete group that admits a metrically proper and isometric action on $M$. Recall that $\Ga\car M$ is metrically proper means that there are only finite many $\gamma$ in $\Ga$ such that $B(x, R)\cap B(\gamma x, R)\ne\emptyset$ for each $x\in M$, $R>0$. Then $\A(M)$ carries a natural and well-behaved action by $\Ga$. We shall then establish the group action of $\Ga$ on $\A(M)$. Denote $\operatorname{Isom}(M)$ the group of all isometries $\varphi:M\to M$ equipped with the weak topology, i.e., a sequence $\{\varphi_n\}\subseteq \operatorname{Isom}(M)$ converge to $\varphi_0$ if and only if
$$d(\varphi_0(x),\varphi_n(x))\to0\quad\text{as}\quad n\to\infty$$
for each $x\in M$. For any $\varphi:M\to M$, we define the tangent map
$$d_x\varphi:T_xM\to T_{\varphi(x)}M\quad\text{by}\quad d_x\varphi([\alpha,t])=[\varphi\circ\alpha,t].$$
Then $d_x\varphi$ induces a unitary $\H d_x\varphi:\H_xM\to\H_{\varphi(x)}M$, which induces a graded $C^*$-isomorphism between the corresponding Clifford algebras:
$$\Cl(\H d_x\varphi\oplus 1_{t\IR}):\Cl(\H_xM\oplus t\IR)\to\Cl(\H_{\varphi(x)}M\oplus t\IR).$$
We define an automorphism $\varphi_*:\Pi_b(M)\to\Pi_b(M)$ by
$$(\varphi_*(\sigma))(x,t)=\Cl(\H d_{\varphi^{-1}(x)}\oplus 1)\varphi(\sigma(\varphi^{-1}(x),t)).$$
One can check by a calculation that $\varphi_*\circ\beta_{x_0}(f)=\beta_{\varphi(x_0)}(f)$, see \cite[Lemma 6.2]{GWY2018}. Thus for any $\gamma\in\Ga$, $\gamma_*$ will take an element in $\A(M)$ to $\A(M)$. This makes $\A(M)$ into a $\Ga$-$C^*$-algebra. The following proposition is proved in section 6 in \cite{GWY2018}:

\begin{Pro}[\cite{GWY2018}]
Let $\Ga$ be a countable discrete group, $\alpha:\Ga\to \operatorname{Isom}(M)$ a metrically proper action on $M$. Then the induced action of $\Ga$ on $\A(M)$ makes $\A(M)$ into a proper $\Ga$-$C^*$-algebra.\qed
\end{Pro}

There exists a canonical inclusion $\operatorname{Isom}(M)\hookrightarrow \operatorname{Isom}(M^{[0,1]})$ by composition of maps, i.e., for any $\varphi\in \operatorname{Isom}(M)$, we define $\varphi^{[0,1]}\in \operatorname{Isom}(M^{[0,1]})$ by
$$(\varphi^{[0,1]}\xi)(t)=\varphi(\xi(t))$$
for any $\xi\in M^{[0,1]}$ and $t\in[0,1]$. Let $\Ga$ be a discrete subgroup of $\operatorname{Isom}(M)$, $\alpha:\Ga\to \operatorname{Isom}(M)$ an isometric metrically proper action on $M$. Then the canonical embedding $\operatorname{Isom}(M)\hookrightarrow \operatorname{Isom}(M^{[0,1]})$ gives rise to an isometric, metrically proper action of $\Ga$ on $M^{[0,1]}$, see \cite[Proposition 3.19]{GWY2018}.

Let $M$ be a Hilbert-Hadamard space. Define
$$H:\operatorname{Isom}(M)\times [0,1]\to \operatorname{Isom}(M^{[0,1]})$$
by
\begin{equation}\label{eqH}
(H(\varphi,t)\xi)(s)=\left\{\begin{aligned}&\varphi(\xi(s))&&,s\in[0,t]\\&\xi(s)&&s\in(t,1]\end{aligned}\right.
\end{equation}
for any $\varphi\in \operatorname{Isom}(M)$, $\xi\in M^{[0,1]}$ and $t,s\in[0,1]$. The following Proposition comes from \cite[Proposition 3.18]{GWY2018}:

\begin{Pro}\label{H-homotopy}
The map $H: \operatorname{Isom}(M)\times [0,1]\to \operatorname{Isom}(M^{[0,1]})$ defined above is a continuous homotopy connecting the canonical embedding $\operatorname{Isom}(M)\hookrightarrow \operatorname{Isom}(M^{[0,1]})$ to the trivial action $\operatorname{Isom}(M)\to\{Id\}\leq \operatorname{Isom}(M^{[0,1]})$, i.e., $H$ satisfies that:
\begin{itemize}
\item[(1)]for each $t\in[0,1]$, $\varphi\mapsto H(\varphi,t)$ is a group homomorphism;
\item[(2)]$H(\varphi,0)=Id$, for any $\varphi\in \operatorname{Isom}(M)$;
\item[(3)]$H(\varphi,1)=\varphi^{[0,1]}$, for any $\varphi\in \operatorname{Isom}(M)$.\qed
\end{itemize}\end{Pro}

\begin{Def}\label{big algebra}
Let $M$ be a Hilbert-Hadamard space, $\Ga$ a discrete countable subgroup of $\operatorname{Isom}(M)$ such that $\Ga\car M$ is proper. Define the $C^*$-algebra:
$$\A_{[0,1]}(M)=C\left([0,1],\A(M^{[0,1]})\right).$$
Equip $\A_{[0,1]}(M)$ with the action of $\Ga$ defined by
$$(\gamma(f))(t)=H(\gamma,t)_*(f(t)),$$
where $H$ is the homotopy defined as above and the lower $*$ means the induced group action on $\A(M^{[0,1]})$.

For each $t\in [0,1]$, we define the evaluation map
$$ev_t:\A_{[0,1]}(M)\to\A(M^{[0,1]}),\quad f\mapsto f(t)$$
which intertwines the $\Ga$-action on $\A_{[0,1]}(M)$ and the $\Ga$-action on $\A(M^{[0,1]})$ associated with $t$, denoted by $\alpha_t$, defined by
$$(\alpha_t(\gamma))(a)=H(\gamma,t)_*(a)$$
for all $a\in \A(M^{[0,1]})$.
\end{Def}

\section{The equivariant coarse Novikov conjecture}\label{sec: The equivariant coarse Novikov conjecture}

In this section, we shall recall some basic facts about equivariant Roe algebras, whose $K$-theory forms the target of the equivariant coarse index map (and thus also the equivariant coarse Baum-Connes assembly map). Moreover, we shall also recall a localization approach to equivariant $K$-homology introduced in \cite{Yu1997}, which greatly simplifies the study of the equivariant coarse index map.

Let $X$ be a proper metric space, and let $\Gamma$ be a countable discrete group acting on $X$ (on the left) properly by isometries. Recall that the action $\Gamma\curvearrowright X$ is {\em proper} if the map $\Gamma\times X\to X$ defined by $(\gamma, x)\mapsto \gamma x$ is a proper map with respect to the product topology on $\Gamma\times X$.
Let $C_0(X)$ be the $C^*$-algebra of all complex-valued continuous functions on $X$ vanishing at infinity. Then the action $\Ga\car X$ induces a $\Ga$-action $\alpha:\Ga\to\text{Aut}(C_0(X))$ by
$$(\alpha_{\gamma}f)(x)=f(\gamma^{-1}x),$$
where the group of $*$-automorphisms of $C_0(X)$ is equipped with the topology of pointwise norm-convergence. To simplify the notation, we shall also denote this action by $\gamma\cdot f$.

The following definition is due to John Roe \cite{Roe1993}:

\begin{Def}
Let $H_X$ be a Hilbert space, and let $\varphi:C_0(X)\to\B(H_X)$ be a $*$-homomorphism from $C_0(X)$ to $\B(H_X)$, the $C^*$-algebra of all bounded operators on $H_X$. Let $T\in\B(H_X)$ be an operator on $H_X$.
\begin{itemize}
\item[(1)] The {\em support} of $T$ is defined to be the complement in $X\times X$ of the set of all points $(x, y)\in X\times X$ for which there exist $f, g\in C_c(X)$ such that $fTg=0$ but $f(x)\not= 0, g(y)\not=0$.
\item[(2)] The {\em propagation} of $T$ is defined to be
$$\prop(T)=\sup\big\{ d(x,y) \mid (x, y)\in \supp(T) \big\}.$$
\item[(3)] The operator $T: H_X\to H_X$ is said to be {\em locally compact} if $fT$ and $Tf$ belong to $\K(H_X)$ for all $f\in C_c(X)$, where $\K(H_X)$ is defined to be the algebra of all compact operators on $H_X$.
\end{itemize}\end{Def}

Let $H_X$ be a $\Gamma$-Hilbert space, i.e., there exists a representation $U$ of $\Gamma$ on $H_X$ which is {\em unitary} in the sense that
$$\langle U_\gamma \xi, \, U_\gamma \eta \rangle = \langle \xi, \eta \rangle$$
for all $\gamma\in \Gamma$, and $\xi, \eta\in H_X$. Let $\varphi: C_0(X)\to \B(H_X)$ be a $*$-homomorphism which is {\em covariant} with $U$ in the sense that
$$Ad_{U_\gamma} (\varphi(f)) := U_\gamma \varphi(f) U_\gamma^{-1} = \varphi(\gamma\cdot f),$$
or simply written as
$$ad_\gamma (f)=\gamma f \gamma^{-1} =\gamma\cdot f,$$
where $(\gamma\cdot f)(x)=f(\gamma^{-1}x)$ for all $\gamma\in \Gamma$, $f\in C_0(X)$ and $x\in X$. Such a quadruple $(C_0(X), \Gamma, \varphi, U)$ is said to be a {\em covariant system} for the system $\Gamma\curvearrowright X$.
\par
An operator $T\in\B(H_X)$ is said to be {\em $\Gamma$-invariant} if $\gamma T \gamma^{-1} =T$ for all $\gamma\in \Gamma$.
\par
\begin{Def}(\cite{KY2012})
A covariant system $(C_0(X), \Gamma, \varphi, U)$ is said to be {\em admissible} if \begin{itemize}
\item[(1)] the action $\Ga\car X$ is proper; 
\item[(2)] there exist $\Gamma$-Hilbert spaces $\H_X$ and $\H_\Gamma$, and a separable and infinite dimensional Hilbert space $\H_0$,
such that $H_X$ is isomorphic to $\H_X\otimes \H_\Gamma \otimes \H_0$ as a $\Gamma$-Hilbert space; 
\item[(3)] the representation $\varphi: C_0(X)\to\B(H_X)$ admits a decomposition $\varphi=\varphi_X\otimes I$ for some $\Gamma$-equivariant $*$-homomorphism $\varphi_X: C_0(X)\to\B(\H_X)$ such that
\\(i) $\varphi_X(f)\notin \K(\H_X)$ for all nonzero $f\in C_0(X)$, and
\\(ii) $\varphi_X$ is non-degenerate in that sense that $\{\varphi_X(f) \xi: f\in C_0(X), \xi\in \H_X\}$ is dense in $\H_X$; 
\item[(4)] $H_X$ is locally free (cf. \cite{HIT2020}), i.e., for any finite subgroup $F$ of $\Gamma$ and any $F$-invariant Borel subset $E$ of $X$ there is a Hilbert space $\H_E$ equipped with the trivial representation of $F$ such that $\chi_E H_X$ is isomorphic to $\H_E\otimes \ell^2(F)$ as $F$-representions.
\end{itemize}\end{Def}

In the above definition, $I$ is the identity operator on $\H_\Gamma\otimes \H_0$ and the $F$-action on $\ell^2(F)$ is the left regular representation, that is, $(\gamma \xi) (z)=\xi(\gamma^{-1} z)$ for all $\xi\in \ell^2(F)$ and $\gamma,z\in F$. And the representation $\varphi:C_0(X)\to\B(H_X)$ is also called an \emph{ample (or absorbing)} representation (see \cite[Definition 4.5.2]{HIT2020}).

\begin{Rem}(\cite{KY2012})
Let $N\subset X$ be a $\Gamma$-invariant countable dense subset of $X$, and let $\H_0$ be an infinite dimensional separable Hilbert space. Then the Hilbert space
$$H_X:=\ell^2(N)\otimes \ell^2(\Gamma)\otimes \H_0$$
equipped with the diagonal action of $\Gamma$ defined by
$$U_\gamma: \delta_z\otimes \delta_h\otimes \xi \mapsto \delta_{\gamma z}\otimes \delta_{\gamma h}\otimes \xi$$
and an action of $C_0(X)$ by pointwise multiplications
$$f: \delta_z\otimes \delta_h\otimes \xi \mapsto f(z)\delta_z\otimes \delta_h\otimes \xi$$
is an admissible covariant system.
\end{Rem}

\begin{Def}
Let $(C_0(X), \Gamma, \varphi, U)$ be an admissible covariant system.\begin{itemize}
\item[(1)]The {\em algebraic equivariant Roe algebra},
denoted by $\IC[X]^{\Ga}$, is defined to be the $*$-algebra of all $\Gamma$-invariant, locally compact, finite propagation operators in $\B(H_X)$.
\item[(2)]The {\em equivariant Roe algebra}, denoted by $C^*(X)^{\Ga}$, is defined to be the $C^*$-algebraic completion of $\mathbb{C}[X]^{\Ga}$ in $\B(H_X)$.
\end{itemize}\end{Def}

The equivariant Roe algebra has a natural "X-by-X" matrix picture in the following sense:

\begin{Rem}[\cite{Yu2000}]\label{Roe finite}
Let $N\subset X$ be a $\Gamma$-invariant countable dense subset of $X$, and let $\H_0$ be an infinite dimensional separable Hilbert space. Let $\IC_{\rm f}[X]^{\Ga}$ be the set of all bounded functions on $N\times N$ such that\begin{itemize}
\item[(1)] $T(x,y)\in \K(\ell^2(\Ga)\ox\H_0)$ for all $x,y\in N$, where $\ell^2(\Ga)$ is equipped with the left regular representation of $\Ga$ and $\H_0$ is equipped with the trivial $\Ga$-action;
\item[(2)] there exists $r>0$ such that $T(x,y)=0$ for all $x,y\in N$ with $d(x,y)>r$;
\item[(3)] for each bounded set $B\subseteq X$, the set
$$\{(x,y)\in (N\times N)\cap (B\times B)\mid T(x,y)\ne0\}$$
is finite;
\item[(4)]there exists $L>0$ such that
$$\#\{y\in N\mid T(x,y)\ne0\}\leq L,\quad\#\{y\in Z\mid T(y,x)\ne0\}\leq L$$
for all $x\in N$;
\item[(5)] for each $\gamma\in\Ga$,
$$T(\gamma x,\gamma y)=\gamma T(x,y)$$
where $\Ga\car\K(\ell^2(\Ga)\ox\H_0)$ is given by $ad_{u_{\gamma}\ox 1}$ and $u_{\gamma}$ is the left regular action on $\ell^2(\Ga)$.
\end{itemize}
The algebraic structure on $\IC_{\rm f}[X]^{\Ga}$ is given by
$$(TS)(x,y)=\sum_{z\in N}T(x,z)S(z,y)$$
$$(T^*)(x,y)=(T(y,x))^*.$$
Notice that $T$ acts on $\delta_x\ox\delta_{\gamma}\ox v\in\ell^{2}(N)\ox\ell^2(\Ga)\ox\H_0$ given by
$$T(\delta_x\ox\delta_{\gamma}\ox v)=\sum_{y\in N}\delta_y\ox T(y,x)(\delta_{\gamma}\ox v),$$
which fits the $*$-structure on $\IC_{\rm f}[X]^{\Ga}$. In the sequel, $\IC_{\rm f}[X]^{\Ga}$ is a dense subalgebra of $C^*(X)^{\Ga}$ and we will use $\IC_{\rm f}[X]^{\Ga}$ to replace $\IC[X]^{\Ga}$ to define the equivariant Roe algebra of $X$.
\end{Rem}

\begin{Def}
Let $X, Y$ be proper metric spaces with proper isometric $\Ga$-actions, and $f: X\to Y$ be a map. Then\begin{itemize}
\item[(1)] $f$ is \emph{metrically proper} if the preimage of a bounded set is bounded;
\item[(2)] $f$ is \emph{bornologous} if there exists a function $\rho_+:[0,\infty)\to[0,\infty)$ such that for all $x_1,x_2\in X$,
$$d(f(x_1),f(x_2))\leq \rho_+(d(x_1,x_2));$$
\item[(3)] $f$ is said to be \emph{coarse} if it is bornologous and metrically proper;
\item[(4)] $f$ is said to be a \emph{coarse equivalence} if there exists $R>0$ and a coarse map $g:Y\to X$ such that $d(g(f(y)),y)\leq R$ and $d(f(g(x)),x)\leq R$ for all $x\in X$ and $y\in Y$;
\item[(5)] $f$ is said to be \emph{$\Ga$-equivariant} if $f(\gamma x)=\gamma f(x)$ for any $\gamma\in\Ga$ and $x\in X$.
\end{itemize}\end{Def}

Let $X, Y$ be proper metric spaces with proper $\Gamma$-actions. If there exists a $\Gamma$-equivariant coarse equivalence $f:X\to Y$, then there exists an equivariant unitary operator $V\in\B(H_X,H_Y)$ associated with $f$ such that the map $ad(V): T\mapsto VTV^*$ gives rise to an isomorphism (cf. \cite{HIT2020})
$$ad(V):C^*(X)^{\Gamma}\xrightarrow{\cong} C^*(Y)^{\Gamma}.$$
As a result, since the identity map $id: X\to X$ is an equivariant coarse equivalence, the equivariant Roe algebra $C^*(X)^\Gamma$ does not depend on the choice of the admissible covariant system. The following result is essentially due to J.~Roe.

\begin{Thm}\label{Roe}
Let $(C_0(X), \Gamma, \varphi, U)$ be an admissible covariant system. If the $\Gamma$-action on $X$ is cocompact, then
$$C^*(X)^\Gamma\cong C^*_r\Gamma\otimes\K$$
where $C^*_r\Gamma$ is the reduced group $C^*$-algebra, and $\K$ is the algebra of all compact operators on a separable infinite dimensional Hilbert space.
\end{Thm}

Let $X$ be a proper metric space with a proper $\Ga$-action. The equivariant $K$-homology group $K_*^{\Ga}(X)$ is defined to be $KK^{\Ga}_*(C_0(X),\IC)$ by using Kasparov's equivariant $KK$-theory (cf. \cite{Kas88}). In the following, we shall recall a localization approach to $K^{\Ga}_*(X)$ which is introduced by G.~Kasparov and G.~Yu in \cite{KY2012}. One is also referred to \cite{Yu1997, QR2010, HIT2020} for backgrounds on localization algebras.

A $KK$-cycle in $KK^{\Ga}_0(C_0(X),\IC)$ is a triple $(H_X,\varphi,F)$, where $H_X$ is a $\Ga$-Hilbert space, $\varphi:C_0(X)\to\B(H_X)$ is a covariant $*$-homomorphism and $F\in\B(H_X)$ satisfies
$$[\varphi(f),F],\,\varphi(f)(FF^*-1),\,\varphi(f)(F^*F-1),\,\varphi(F)(\gamma\cdot F-F)\in\K(H_X)$$
for all $f\in C_0(X)$ and $\gamma\in\Ga$.

For any $R> 0$, one can always take a locally finite, $\Ga$-equivariant open cover $\{U_i\}_{i\in I}$ of $X$ such that the diameter of each $U_i$ is no more than $R$. Let $\{\phi_i\}_{i\in I}$ be a continuous, $\Ga$-equivarinat partition of unity subordinate to the open cover $\{U_i\}_{i\in I}$. Define
$$F'=\sum_{i}\varphi(\phi_i^{\frac 12})F\varphi(\phi_i^{\frac 12}),$$
where the sum converges in the strong operator topology. It is not hard to see that $(H_X,\varphi,F)$ and $(H_X,\varphi,F')$ are equivalent via $(H_X,\varphi,(1-t)F+tF')$, where $t\in[0,1]$. Note that both $F'$ and $F'^2-1$ have finite propagation, so $F'$ is a multiplier of $C^*(X)^{\Ga}$ and $F'$ is invertible modulo $C^*(X)^{\Ga}$. Hence $F'$ gives rise to an element, denoted by $\partial([F'])$ in $K_0(C^*(X)^{\Ga})$, where
$$\partial:K_1\left(M\left(C^*(X)^{\Ga}\right)/C^*(X)^{\Ga}\right)\to K_0\left(C^*(X)^{\Ga}\right)$$
is the boundary map in $K$-theory, and $M\left(C^*(X)^{\Ga}\right)$ is the multiplier algebra of $C^*(X)^{\Ga}$. We define the index map
$$\ind:KK^{\Ga}_0(C_0(X),\IC)\to K_0(C^*(X)^{\Ga})$$
by
$$\ind([H_X,\varphi,F])=\partial([F']).$$
Similarly, one can also define $\ind:KK^{\Ga}_1(C_0(X),\IC)\to K_1(C^*(X)^{\Ga})$.

\begin{Def}\label{equivariant localization algebra}
Let $(C_0(X), \Gamma, \varphi, U)$ be an admissible covariant system.\begin{itemize}
\item[(1)]The {\em algebraic equivariant localization algebra},
denoted by $\mathbb{C}_L[X]^{\Ga} $, is defined to be the $*$-algebra of all bounded and uniformly continuous functions $f:[0,\infty)\to\IC_{\rm f}[X]^{\Ga}$ such that
$$\prop(f(t))\to0$$
as $t\to\infty$.
\item[(2)]The {\em equivariant localization algebra}, denoted by $C^*_L(X)^{\Ga}$, is defined to be the $C^*$-algebraic completion of $\mathbb{C}_L[X]^{\Ga}$ under the norm
$$\|f\|=\sup_{t\in[0,\infty)}\|f(t)\|.$$
\end{itemize}\end{Def}

It is proved that $C^*_L(X)^{\Ga}$ does not depend on the choice of the admissible covariant system \cite{KY2012}. There exists a local index map $\ind_L:KK_*^{\Ga}(C_0(X),\IC)\to K_*(C^*_L(X)^{\Ga})$ as in \cite{KY2012} (see \cite{Yu1997} for a non-equivariant version). We have the following result:

\begin{Pro}[\cite{KY2012}]\label{localization}
If $X$ is a finite-dimensional simplicial polyhedron with a $\Ga$-invariant metric. Then the local index map
$$\ind_L:KK_*^{\Ga}(C_0(X),\IC)\to K_*(C^*_L(X)^{\Ga})$$
is an isomorphism.\qed
\end{Pro}

Define the evaluation map $ev:C^*_L(X)^{\Ga}\to C^*(X)^{\Ga}$ by $ev(f)=f(0)$ for each $f\in C^*_L(X)^{\Ga}$. Then we have the following commuting diagram:
\[\begin{tikzcd}
                                                                       && {K_*(C^*_L(X)^{\Gamma})} \arrow[d, "ev_*"] \\
{KK_*^{\Gamma}(C_0(X),\IC)} \arrow[rr, "\text{Ind}"] \arrow[rru, "\text{Ind}_L"] && {K_*(C^*(X)^{\Ga}).}                    
\end{tikzcd}\]

To formulate the equivariant coarse Novikov conjecture, we shall recall the Rips complexes associated with $\Ga\car X$ which will play a role of classifying space. We should mention that we can only consider $X$ to be discrete.  Recall that a discrete metric space $X$ is said to have \emph{bounded geometry} if for any $R>0$, there exists $N>0$ such that $\sup_{x\in X}\#B(x, R)\leq N$. For any $r>0$, a subset $W\subseteq X$ is said to be a \emph{$r$-net} if for any $x\in X$, there exists $w\in W$ such that $d(w,x)\leq r$. A non-discrete metric space $X'$ with proper $\Ga$-action is said to have \emph{bounded geometry} if there exists a $\Ga$-invariant $r$-net $W$ of $X'$ such that $W$ has bounded geometry as a discrete metric space.

Fix $\Ga$-space $X$. For any $r>0$, one can always choose $W\subseteq X$ to be a $\Ga$-invariant set such that the distance of any two orbits in $W$ is larger than $r$ by using Zorn's Lemma. Then one can see that $W$ must be a $r$-net by the maximal assumption and $W$ is locally finite by the properness of the $\Ga$-action. Since the equivariant coarse Baum-Connes conjecture is invariant under equivariant coarse equivalence, we can replace $\Ga\car X$ by $\Ga\car W$. From now on. we shall always assume $X$ to have bounded geometry.

For each $n\in\IN$, let $\mathcal{S}(\mathbb{R}^n)$ be the sphere of radius one in an $n$-dimensional Euclidean space. The metric on $\mathcal{S}(\mathbb{R}^d)$ is defined by
$$d_{S(\IR^n)}(x,y)=\arccos(\langle x,y\rangle),$$
where $\langle x,y\rangle$ means the standard inner product of $\IR^n$. It is the usual Riemannian metric on $\mathcal{S}(\mathbb{R}^d)$. Before further discussion, let us recall some basic definitions of the Rips complex.

\begin{Def}[Rips complex]
For each $d\geq 0$, the \emph{Rips complex} of $X$ at scale $d$, denoted $P_d(X)$, consists as a set of all formal sums
$$z=\sum_{z\in X}t_xx$$
such that each $t_x$ is in $[0,1]$, such that $\sum_{x\in X}t_x=1$, and such that the support of $x$ defined by $\supp(z):=\{x\in X\,|\,t_x\ne0\}$ has diameter at most $d$.
\end{Def}

Let $F\subseteq X$ be a finite set. Let $\mathbb{R}^F=\{\sum_{z\in F}t_zz\mid t_z\in\IR\text{ and }z\in F\}$ be finite-dimensional real Euclidean space spanned by $F$. Denote $\Delta(F)$ the simplex generated by $F$, i.e.
$$\Delta(F)=\left\{\sum_{z\in F}t_zz\;\Big|\;t_z\geq 0\text{ and }\sum_{z\in F}t_z=1\right\}\subseteq \IR^F.$$
Define 
$$f:\Delta(F)\to\mathcal{S}(\mathbb{R}^F)\quad\text{by}\quad \sum_{z\in F}t_zz\mapsto\sum_{z\in F}\frac{t_z}{\sqrt{\sum_{z\in F}t^2_z}} z$$
The spherical metric on $\Delta(F)$ is the metric defined by
$$d_S(x,y)=\frac{2}{\pi}d_{S(\IR^F)}(f(x),f(y))$$
and $\Delta(F)$ equipped with this metric is called the \emph{spherical simplex} on $F$. The factor $2/\pi$ is chosen so that $\Delta(F)$ has diameter one.

For any $x,y\in P_d(X)$, a semi-simplicial path $\delta$ between $x$ and $y$ is a finite sequence of points
$$x=x_0,y_0,x_1,y_1,\cdots,x_n,y_n=y$$
where $x_i$ and $y_{i-1}$ are in $X$ for each $i=\{1,\cdots,n\}$. The length of $\delta$ is defined to be
$$l(\delta)=\sum_{i=0}^{n}d_X(x_i,y_i)+\sum_{i=0}^{n-1}d_S(y_i,x_{i+1})$$
We define the \emph{semi-spherical metric} $d_{P_d}$ on $P_d(X)$ by
$$d_{P_d}(x,y)=\inf\{l(\delta)\mid\text{ $\delta$ is a semi-simplicial path between $x$ and $y$}\}.$$
One can show that $d_{P_d}$ is the largest metric such that $d_{P_d}(x,y)\leq d_S(x,y)$ and $d_{P_d}(z_1,z_2)\leq d_X(z_1,z_2)$ for any $x,y\in \Delta(F)\subseteq P_d(X)$ for some finite $F\subseteq X$ and $z_1,z_2\in X$.

The Rips complex is made into a metric space under the metric $d_{P_d}$. One can see that $(P_0(X),d_{P_0})$ is identified with the metric space $(X,d_X)$. The $\Ga$-action on $P_d(X)$ is defined by
$$\gamma\cdot\left(\sum_{x\in X}t_xx\right)=\sum_{x\in X}t_x\gamma x.$$
It is proved in \cite{HIT2020} that the canonical inclusion $i: X\to P_d(X)$ is an equivariant coarse equivalence for each $d\geq0$. Then the equivariant coarse Novikov conjecture states as follows:

\begin{ECNcon}
Let $X$ be a proper metric space with bounded geometry, $\Ga$ a countable discrete group which acts properly and isometrically on $X$. Then the equivariant coarse assembly map induced by the evaluation map is rational injective, i.e.,
$$\mu_X^{\Ga}:\lim_{d\to\infty}K_*(C^*_L(P_d(X))^{\Ga})\to\lim_{d\to\infty}K_*(C^*(P_d(X)^{\Ga}))\cong K_*(C^*(X)^{\Ga})$$
is a (rational) injection.
\end{ECNcon}

The conjecture above is strengthened by the \emph{equivariant coarse Baum-Connes conjecture}, which claims the equivariant coarse assembly map is (rational) bijective. It is worth mentioning that when the action $\Ga\car X$ is cocompact, then the topological side of the assembly map is identified with the representable $K$-homology group of the classifying space for proper $\Ga$-actions, which is always denoted by $RK_*^{\Ga}(\EG)$. Combining Theorem \ref{Roe}, the equivariant coarse assembly is identified with the Baum-Connes assembly map for the group $\Ga$, i.e.,
$$\mu^{\Ga}:RK^{\Ga}_*(\EG)\to K_*(C^*_r\Ga).$$
When $\Ga$ is the trivial group, the map
$$\mu_X:\lim_{d\to\infty}K_*(C^*_L(P_d(X)))\to K_*(C^*(X))$$
is called the coarse assembly map.

\section{Twisted algebras and the Bott homomorphisms}\label{sec: Twisted algebras and the Bott homomorphisms}

In this section, we introduce the (equivariant) twisted Roe algebras and the (equivariant) twisted localization algebras associated to (equivariant) coarse embeddings into Hilbert-Hadamard spaces. The constructions presented here are natural generalizations of those from \cite{Yu2000}, which tackles embeddings into Hilbert spaces. 
However, in preparation for a \emph{deformation trick} (see \Cref{sec: A deformation trick}), pivotal in our method of proof, we shall need to deviate from the standard construction to equip our (equivariant) twisted Roe (respectively, localization) algebras with an extra ``deformation parameter'' indicated notationally by the unit interval $[0,1]$. 

\begin{Def}
Let $X$ be a proper metric space, $\Ga$ a countable discrete group with a proper isometric action on $X$. A \emph{(Roe algebraic) coefficient system in $M$} for $\Ga\car X$ is a triple $(N,M,f)$ consisting of the following data:\begin{itemize}
\item[(1)] a $\Ga$-invariant countable dense subset $N\subseteq X$;
\item[(2)] an admissible Hilbert-Hadamard space $M$ with a $\Ga$-action, $\alpha:\Ga\to\operatorname{Isom}(M)$;
\item[(3)] an $\Ga$-equivariant map $f:N\to M$.
\end{itemize}\end{Def}

\begin{Def}\label{twisted Roe small}
The \emph{algebraic equivariant twisted Roe algebra} for $\Ga\car X$ associated with a coefficient system $(N,M,f)$, denoted by $\IC[X,\A(M)]^{\Ga}$, is the set of all bounded functions $T$ on $N\times N$ such that:\begin{itemize}
\item[(1)] $T(x,y)\in\K(\ell^2(\Ga)\ox\H_0)\ox\A(M)$ for each $x,y\in N$ and there exists $R>0$ such that
$$\supp(T(x,y))\subseteq B(f(x),R)\subseteq M\times\IR_+,$$
where $f:N\to M$ is the equivariant map;
\item[(2)] there exists $r>0$ such that $T(x,y)=0$ for all $x,y\in N$ with $d(x,y)>r$;
\item[(3)] for each bounded set $B\subseteq X$, the set
$$\{(x,y)\in (N\times N)\cap (B\times B)\mid T(x,y)\ne0\}$$
is finite;
\item[(4)]there exists $L>0$ such that
$$\#\{y\in N\mid T(x,y)\ne0\}\leq L,\quad\#\{y\in N\mid T(y,x)\ne0\}\leq L$$
for all $x\in N$;
\item[(5)] for each $\gamma\in\Ga$,
$$T(\gamma x,\gamma y)=\gamma T(x,y),$$
where $\Ga\car\K(\ell^2(\Ga)\ox\H_0)\ox\A(M)$ is given by $ad_{u_{\gamma}\ox 1}\ox \gamma_*$ for each $\gamma\in\Ga$, where $u_{\gamma}$ is the left regular action on $\ell^2(\Ga)$ and $\gamma_*$ is the induced action on $\A(M)$.
\end{itemize}\end{Def}

Compared to Remark \ref{Roe finite} the original version, we add a "twisted condition" (Condition (1) above) in this new algebra to do the cutting and pasting trick in Section \ref{sec: The twisted assembly map and proof of Theorem 1.3}.

The algebraic structure on $\IC[X,\A(M)]^{\Ga}$ is given by
$$(TS)(x,y)=\sum_{z\in Z}T(x,z)S(z,y)$$
$$(T^*)(x,y)=(T(y,x))^*.$$
Let $E_{\A(M)}=\ell^2(N)\ox\ell^2(\Ga)\ox\H_0\ox\A(M)$ be a $\Ga$-Hilbert $C^*$-module over $\A(M)$. The $\Ga$-action on $E_{\A(M)}$ is given by
$$\gamma(\delta_z\ox\delta_{\eta}\ox v\ox a)=\delta_{\gamma z}\ox\delta_{\gamma\eta}\ox v\ox \gamma_*a,$$
where $z\in N$, $v\in\H_0$, $a\in\A(M)$ and $\gamma,\eta\in\Ga$. It is clear that an element $T\in\IC[X,\A(M)]^{\Ga}$ defines an adjointable module homomorphism on $E_{\A(M)}$:
$$T(\delta_x\ox\delta_{\gamma}\ox v\ox a)=\sum_{y\in N}\delta_y\ox T(y,x)(\delta_{\gamma}\ox v\ox a),$$
which fits the $*$-structure and the $\Ga$-action on $\IC[X,\A(M)]^{\Ga}$. This gives a $C^*$-representation of $\IC[X,\A(M)]^{\Ga}$ on $E_{\A(M)}$.

\begin{Def}
The \emph{twisted equivariant Roe algebra}, denoted by $C^*(X,\A(M))^{\Ga}$, is the norm closure of $\IC[X,\A(M)]^{\Ga}$ in $\L\left(E_{\A(M)}\right)$, the $C^*$-algebra of all adjointable module homomorphisms.
\end{Def}

\begin{Def}\label{twisted localization small}
Define $\IC_L[X,\A(M)]^{\Ga}$ to be the set of all bounded, uniformly norm-continuous functions
$$g:\IR_+\to\IC[X,\A(M)]^{\Ga}$$
such that\begin{itemize}
\item[(1)]there exists a bounded function $r:\IR_+\to\IR_+$ with $\lim_{s\to\infty}r(s)=0$ such that $(g(s))(x,y)=0$ whenever $d(x,y)>r(s)$;
\item[(2)]there exists $R>0$ such that $\supp((g(s))(x,y))\subseteq B(f(x),R)\subseteq M\times\IR_+$ for all $s\in\IR_+$.
\end{itemize}\end{Def}

Condition (1) in the above definition implies that the propagation of $g(t)$ will tend to 0 as $t$ tends to $\infty$. Condition (2) is the twisted condition for the localization algebra.

\begin{Def}
The \emph{equivariant twisted localization algebra} for $\Ga\car X$ associated to the coefficient system $(N,M,f)$, denoted by $C^*_L(X,\A(M))^{\Ga}$, is the completion of $\IC_L[X,\A(M)]^{\Ga}$ with respect to the norm
$$\|g\|=\sup_{s\in\IR_+}\|g(s)\|.$$
\end{Def}

Let $X$ be a metric space with bounded geometry, $(M,m)$ an admissible Hilbert-Hadamard space with base point $m\in M$. Let $\Ga$ be a countable discrete group that admits a proper isometric action on $X$ and $M$. Assume that $f: X\to M$ is an equivariant coarse embedding. For each $d>0$ and $x_0\in X$, denote
$$Star(x_0)=\left\{\sum_{x\in X}t_xx\in P_d(X)\;\Big|\;t_{x_0}>t_x\text{ for all }x\ne x_0\right\}$$
the open star of $x$ in the barycentric subdivision of $P_d(X)$. It is clear that $\bigsqcup_{x\in X}Star(x)$ is dense in $P_d(X)$. As $\Ga\car X$ is isometric, we have $\gamma\cdot Star(x)\subseteq Star(\gamma x)$ by definition. Take a countable $\Ga$-invariant dense subset $N_d\subseteq P_d(X)$ such that
$$N_d\subseteq\bigsqcup_{x\in X} Star(x)\quad\text{ and }\quad N_d\subseteq N_d'\ \text{ for any }d< d'.$$
For each $z\in N_d$, there exists a unique $x_z\in X$ such that $z\in Star(x_z)$. Then the map $i: N_d\to X$ defined by $z\mapsto x_z$ is an equivariant coarse equivalence. To simplify the notation, we still denote it $f: N_d\to M$ the composition of $f: X\to M$ and $i: N_d\to X$. It is clear that $f$ is an equivariant coarse embedding. Then the triple $(N_d,M,f)$ forms a Roe algebraic coefficient system for $\Ga\car P_d(X)$.

However, since the $K$-theory of $\A(M)$ is not fully clear to us, we are not able to, even partially, calculate the $K$-theory of the twisted algebras for $P_d(X)$ associated with $(N_d, M,f)$. Inspired by \cite{GWY2018}, this problem actually can be solved by amplifying the Hilbert-Hadamard space $M$ to a bigger one, $M^{[0,1]}$ defined as in Section 2. Since the canonical inclusion $M\to M^{[0,1]}$ is an equivariant isometry, $N_d$ still equivariantly coarsely embeds into $M^{[0,1]}$. The $K$-theory of the twisted algebras associated with $M^{[0,1]}$ is more computable since one can use the homotopy in \eqref{eqH} to kill the group action. To formalize this trick, we shall also need the following coefficient systems and even larger versions of twisted algebras. 

Fix a base point $m\in M$. Let $(N, M,f)$ be a Roe algebraic coefficient system for $\Ga\car X$. For each $t\in[0,1]$, define $f_t:N\to M^{[0,1]}$ by
\begin{equation}\label{eqf}
(f_t(z))(s)=\left\{\begin{aligned}&f(z),&&s\in[0,t];\\&m,&&s\in(t,1].\end{aligned}\right.
\end{equation}
Then $f_t$ is an equivariant map associated with the $\Ga$-action on $M^{[0,1]}$ given by $\alpha_t$ defined as in Definition \ref{big algebra} and \eqref{eqH}, and $f_1$ is the composition of $f:N\to M$ and the constant inclusion $M\hookrightarrow M^{[0,1]}$. Then for each $t\in[0,1]$, the triple $(N,M^{[0,1]},f_t)$ is also a Roe algebraic coefficient system for $\Ga\car X$, where $M^{[0,1]}$ is equipped with the $\Ga$-action defined by $\alpha_t$.

{\bf Notation:} Let $(N_d, M,f)$ be the coefficient system for $\Ga\car P_d(X)$. For each $t\in [0,1]$, we can then define the coefficient systems $(N_d, M^{[0,1]},f_t)$ as above. We denote the twisted equivariant Roe algebra (respectively, twisted equivariant localization algebra) for $\Ga\car P_d(X)$ associated with $(N_d, M^{[0,1]},f_t)$ by $C^*(P_d(X),\A(M^{[0,1]}))^{\Ga,(t)}$ (respectively, $C^*_L(P_d(X),\A(M^{[0,1]}))^{\Ga,(t)}$).

We still need the following larger version of twisted algebras. We enlarge the coefficient algebra from $\A(M^{[0,1]})$ to $\A_{[0,1]}(M)$ (see Definition \ref{big algebra}). The algebra $\A_{[0,1]}(M)$ will build a bridge between the proper $\Ga$-action on $\A(M^{[0,1]})$ to the trivial action on $\A(M^{[0,1]})$. The idea of the following constructions comes from \cite[Lemma 8.3]{GWY2018}.

\begin{Def}\label{twist Roe big}
For each $t\in[0,1]$, let $(N,M,f)$ be a Roe alegbraic coefficient system for $\Ga\car X$ with $f$ an equivariant coarse embedding,  $(N,M^{[0,1]},f_t)$ as above. Define $\IC[X,\A_{[0,1]}(M)]^{\Ga}$ to be the set of all bounded functions $T$ on $N\times N$ such that\begin{itemize}
\item[(1)] $T(x,y)\in\K(\ell^2(\Ga)\ox\H_0)\ox\A_{[0,1]}(M)$ for each $x,y\in N$ and there exists $R>0$ such that
$$\supp\Big((1\ox ev_t)\Big(T(x,y)\Big)\Big)\subseteq B(f_t(x),R)\subseteq M^{[0,1]}\times\IR_+,$$
where $t\in[0,1]$, $ev_t:\A_{[0,1]}(M)\to\A(M^{[0,1]})$ is the evaluation map introduced in Definition \ref{big algebra}, $f_t:N\to M^{[0,1]}$ is the equivariant map associated to the $\Ga$-action on $M^{[0,1]}$ given by $H(\cdot,t)$.
\item[(2)] $T$ satisfies condition (2)-(5) in Definition \ref{twisted Roe small}, where the $\Ga$-action on $\A_{[0,1]}(M)$ is defined as in Definition \ref{big algebra}.
\item[(3)] the sequence of functions $\{t\mapsto(T(x,y))(s_x\cdot t)\}_{x,y\in N}$ is equi-continuous, where we view $T(x,y)$ as a continuous function from $[0,1]$ to $\K(\ell^2(\Ga)\ox\H_0)\ox\A(M^{[0,1]})$ and $s_x=d(f(x),m)^{-1}$.
\end{itemize}\end{Def}

Condition (3) in the definition above is a technical condition that will be used in Lemma \ref{deformation} for an Eilenberg swindle argument.
The algebraic structure of $\IC[X,\A_{[0,1]}(M)]^{\Ga}$ is defined similarly with $\IC[X,\A(M^{[0,1]})]^{\Ga}$. Let $E_{\A_{[0,1]}(M)}=\ell^2(N)\ox\ell^2(\Ga)\ox\H_0\ox\A_{[0,1]}(M)$ be the $\Ga$-Hilbert $\A_{[0,1]}(M)$-module. The twisted algebra $\IC[X,\A_{[0,1]}(M)]^{\Ga}$ admits a faithful representation on $E_{\A_{[0,1]}(M)}$. Define $C^*(X,\A_{[0,1]}(M))^{\Ga}$ to be the norm closure of $\IC[X,\A_{[0,1]}(M)]^{\Ga}$ in $\L\left(E_{\A_{[0,1]}(M)}\right)$.

\begin{Def}\label{twist localization big}
Define $\IC_L[X,\A_{[0,1]}(M)]^{\Ga}$ to be the set of all bounded, uniformly norm-continuous functions
$$g:\IR_+\to\IC[X,\A_{[0,1]}(M)]^{\Ga}$$
such that $g$ satisfies condition (1) in Definition \ref{twisted localization small} and there exists $R>0$ such that
$$\supp\Big((1\ox ev_s)\Big((g(t))(x,y)\Big)\Big)\subseteq B(f_s(x),R)$$
for all $t\in\IR_+$ and $s\in[0,1]$.

Define $C^*_L(X,\A_{[0,1]}(M))^{\Ga}$ to be the completion of $\IC_L[X,\A_{[0,1]}(M)]^{\Ga}$ with respect to the norm
$$\|g\|=\sup_{t\in\IR_+}\|g(t)\|.$$
\end{Def}

To simplify the description, we shall omit $\Ga$ on the right-top corner for all these twisted algebras if the action on $X$ is trivial. For example, if $X$ is equipped with a trivial $\Ga$-action and $(N,M^{[0,1]},f_t)$ is a coefficient system, the corresponding twisted Roe algebras will be denoted by $C^*(X,\A(M^{[0,1]}))$ and $C^*(X,\A_{[0,1]}(M))$, similar to twisted localization algebras.

Our main purpose in this section is to build the following commuting diagram:
\begin{equation}\label{com diag}\begin{tikzcd}
& K_{*+1}\Big(C^*_L(P_d(X))^{\Ga}\Big) \arrow[r, "ev_*"] \arrow[d, "(\beta_L)_*"] & K_{*+1}\Big(C^*(P_d(X))^{\Ga}\Big) \arrow[d, "\beta_*"]   \\
& K_*\Big(C^*_L(P_d(X),\A_{[0,1]}(M))^{\Ga}\Big)\arrow[r, "ev_*"] \arrow[d, "(ev_1)_*"]& K_*\Big(C^*(P_d(X),\A_{[0,1]}(M))^{\Ga}\Big) \arrow[d, "(ev_1)_*"]\\
& K_*\Big(C^*_L(P_d(X),\A(M^{[0,1]}))^{\Ga,(1)}\Big)\arrow[r, "ev_*"] & K_*\Big(C^*(P_d(X),\A(M^{[0,1]}))^{\Ga,(1)}\Big)
\end{tikzcd}\end{equation}

The horizontal maps are induced by the evaluation at zero maps. For example, the map
$$ev:C^*_L(P_d(X),\A(M^{[0,1]}))^{\Ga,(1)}\to C^*(P_d(X),\A(M^{[0,1]}))^{\Ga,(1)}$$
is defined by $g\mapsto g(0)$. The induced map on $K$-theory
\begin{equation}\label{eq:twisted-assembly}
    ev_*:\lim_{d\to\infty}K_*\Big(C^*_L(P_d(X),\A(M^{[0,1]}))^{\Ga,(1)}\Big)\to\lim_{d\to\infty}K_*\Big(C^*(P_d(X),\A(M^{[0,1]}))^{\Ga,(1)}\Big)
\end{equation}
is called the \emph{twisted assembly map}. The lower vertical maps are induced by the evaluation map on the coefficients as in Definition \ref{big algebra}, i.e.,
$$ev_1:\A_{[0,1]}(M)\to\A(M^{[0,1]})\quad f\mapsto f(1).$$
For example, the map
$$ev_1:C^*(P_d(X),\A_{[0,1]}(M))^{\Ga}\to C^*(P_d(X),\A(M^{[0,1]}))^{\Ga,(1)}$$
is defined by $T\mapsto ev_1(T)$, where $(ev_1(T))(x,y)=(1\ox ev_1)(T(x,y))$ and $(1\ox ev_1)$ is defined on $\K(\ell^2(\Ga)\ox\H_0)\ox\A_{[0,1]}(M)$.

It remains to define the upper vertical maps, i.e., the Bott maps. For any $g\in\S$, we define $\sigma:[1,\infty)\car\S$ to be the rescaling action given by
$$(\sigma_s(g))(t)=g(s^{-1}t)$$ 
for any $s\in[1,\infty)$, $g\in\S$ and $t\in\IR$. For each $g\in \S$ and $x\in N_d$, define
$$\beta^{[0,1]}(x):\S\to\A_{[0,1]}(M)$$
by using the formula:
$$\Big(\Big(\beta^{[0,1]}(x)\Big)(g)\Big)(s)=\beta^{M^{[0,1]}}_{f_s(x)}(g)\in\A(M^{[0,1]}),$$
where $\beta^{M^{[0,1]}}_{f_s(x)}$ is the Bott map for $M^{[0,1]}$ with base point $f_s(x)\in M^{[0,1]}$. By \cite[Proposition 5.11]{GWY2018}, the map $[0,1]\to\A(M^{[0,1]})$ defined by $s\mapsto\beta^{M^{[0,1]}}_{f_s(x)}(g)$ is continuous as $s\mapsto f_s(x)$ is continuous for fixed $x\in N_d$. Notice that $d(f_s(x),f_{s'}(x))=|s-s'|\cdot d(x,m)$. Let $t_x=d(x,m)^{-1}$ for each $x\in N_d$. Then the sequence of functions
$$\left\{s\mapsto\Big(\Big(\beta^{[0,1]}(x)\Big)(g)\Big)(t_x\cdot s)\right\}_{x\in N}$$
is clearly equi-continuous.

Since $f_s$ is equivariant associated to the $\Ga$-action on $M^{[0,1]}$ given by $H(\cdot,s)$, then one can check that
\[\Big(\Big(\beta^{[0,1]}(\gamma x)\Big)(g)\Big)(s)=\beta^{M^{[0,1]}}_{f_s(\gamma x)}(g)=H(\gamma,s)_*(\beta^{M^{[0,1]}}_{f_s(x)}(g)).\]
This shows that $\beta^{[0,1]}(\gamma x)=\gamma\cdot\beta^{[0,1]}(x)$.

\begin{Def}\label{Big Bott map}
For each $t\in[1,\infty)$ and $d\geq 0$, the \emph{Bott map for the Roe algebras}
$$\beta_t:\S\ox\IC_{\rm f}[P_d(X)]^{\Ga}\to\IC[P_d(X),\A_{[0,1]}(M)]^{\Ga}$$
is defined by the formula
$$(\beta_t(g\ox T))(x,y)=T(x,y)\ox(\beta^{[0,1]}(x))(\sigma_t(g))$$
for all $g\in\S$ and $T\in\IC_{\rm f}[P_d(X)]^{\Ga}$.

The \emph{Bott map for the localization algebras}
$$(\beta_L)_t:\S\ox\IC_L[P_d(X)]^{\Ga}\to\IC_L[P_d(X),\A_{[0,1]}(M)]^{\Ga}$$
is defined by the formula
$$\Big((\beta_L)_t(g\ox h)\Big)(r)=\beta_t(g\ox h(r))$$
for all $g\in\S$ and $h\in\IC_{L}[P_d(X)]^{\Ga}$.
\end{Def}

Notice that
\[\begin{split}\Big(\beta_t(g\ox T)\Big)(\gamma x,\gamma y)&=T(\gamma x,\gamma y)\ox\Big(\beta^{[0,1]}(\gamma x)\Big)(\sigma_t(g))\\
&=ad_{u_{\gamma}\ox 1}\cdot T(x,y)\ox \gamma\cdot(\beta^{[0,1]}(x))(\sigma_t(g)).\end{split}\]
This shows that $\beta_t(g\ox T)$ is $\Ga$-invariant. It also satisfies condition (3) in Definition \ref{twist Roe big}. Recall that an equivariant asymptotic morphism from $A$ to $B$ is an equivariant $C^*$-homomorphism
$$\alpha:A\to\frac{C_b([1,\infty),B)}{C_0([1,\infty),B)}.$$
Since $t\mapsto\beta_t(g\ox T)$ is obviously continuous, it gives a linear map
$$\S\ox \IC_{\rm f}[P_d(X)]^{\Ga}\to C_b([1,\infty),C^*(P_d(X),\A_{[0,1]}(M))^{\Ga}).$$

\begin{Lem}
For each $d\geq 0$, the maps $\beta_t$ and $(\beta_L)_t$ extend to asymptotic morphisms
$$\beta:\S\ox C^*(P_d(X))^{\Ga}\to C^*(P_d(X),\A_{[0,1]}(M))^{\Ga}$$
$$\beta_L:\S\ox C^*_L(P_d(X))^{\Ga}\to C^*_L(P_d(X),\A_{[0,1]}(M))^{\Ga}.$$
\end{Lem}

\begin{proof}
It suffices to show that
$$\|\beta_t(TS)-\beta_t(T)\beta_t(S)\|\to 0$$
as $t\to\infty$ for any $S,T\in\IC_{\rm f}[P_d(X)]$. By a easy calculation and using \cite[Lemma 3.4]{GWY2008}, it suffices to show that for any $r>0$ and $\varepsilon>0$, there exists $N>0$ such that
$$\|\beta^{[0,1]}(x)(\sigma_t(g))-\beta^{[0,1]}(y)(\sigma_t(g))\|\leq \varepsilon$$
for all $t>N$ and $x,y\in N_d$ with $d(x,y)<r$. Since $f: N_d\to M$ is a coarse embedding, there exists $r'>0$ such that $d(f(x),f(y))<r'$ for all $x,y$ with $d(x,y)<r$.

For any $s\in[0,1]$, one has
\[\begin{split}&\left\|\left(\beta^{[0,1]}(x)(\sigma_t(g))-\beta^{[0,1]}(y)(\sigma_t(g))\right)(s)\right\|=\|\beta^{M^{[0,1]}}_{f_s(x)}(\sigma_t(g))-\beta^{M^{[0,1]}}_{f_s(y)}(\sigma_t(g))\|\\
=&\sup_{(z,r)\in M^{[0,1]}\times\IR_+}\left\|\sigma_t(g)\left(C^{M^{[0,1]}}_{f_s(x)}(z,r)\right)-\sigma_t(g)\left(C^{M^{[0,1]}}_{f_s(y)}(z,r)\right)\right\|\\
=&\sup_{(z,r)\in M^{[0,1]}\times\IR_+}\left\|g\left(-\frac{\log_z(f_s(x))}t,\frac rt\right)-g\left(-\frac{\log_z(f_s(y))}t,\frac rt\right)\right\|.
\end{split}\]
Set $z_{x,t,s}=\exp_{z}\left(\frac{\log_z(f_s(x))}t\right)\in M^{[0,1]}$ and $z_{y,t,s}=\exp_{z}\left(\frac{\log_z(f_s(y))}t\right)\in M^{[0,1]}$. Then it follows that
\[\begin{split}&\sup_{(z,r)\in M^{[0,1]}\times\IR_+}\left\|g\left(-\frac{\log_z(f_s(x))}t,\frac rt)\right)-g\left(-\frac{\log_z(f_s(y))}t,\frac rt)\right)\right\|\\
=&\sup_{(z,r)\in M^{[0,1]}\times\IR_+}\left\|g\left(-\log_z(z_{x,t,s}),r\right)-g\left(-\log_z(z_{y,t,s}),r\right)\right\|\\
=&\|\beta^{M^{[0,1]}}_{z_{x,t,s}}(\sigma_t(g))-\beta^{M^{[0,1]}}_{z_{y,t,s}}(g)\|.
\end{split}\]
By the CAT(0) condition, one has that
$$d(z_{x,t,s},z_{y,t,s})\leq \frac 1td(f_s(x),f_s(y))\leq \frac std(f(x),f(y))\leq \frac {sr'}t.$$
By using \cite[Proposition 5.11]{GWY2018}, there exists a large $N>0$ such that 
$$\|\beta^{M^{[0,1]}}_{z_{x,t,s}}(\sigma_t(g))-\beta^{M^{[0,1]}}_{z_{y,t,s}}(\sigma_t(g))\|\leq\varepsilon$$
for all $t\geq N$ and $s\in[0,1]$, which also shows that
$$\|\beta^{[0,1]}(x)(\sigma_t(g))-\beta^{[0,1]}(y)(\sigma_t(g))\|\leq \varepsilon$$
for all $t\geq N$. This completes the proof.
\end{proof}

The asymptotic morphisms induce homomorphisms on $K$-theory:
$$\beta_*:K_{*+1}(C^*(P_d(X))^{\Ga})\to K_*(C^*(P_d(X),\A_{[0,1]}(M))^{\Ga})$$
$$(\beta_L)_*:K_{*+1}(C^*_L(P_d(X))^{\Ga})\to K_*(C^*_L(P_d(X),\A_{[0,1]}(M))^{\Ga}).$$
Thus, we finally complete the commutative diagram \eqref{com diag}.

\section{A deformation trick}\label{sec: A deformation trick}

In this section, we show, under the assumption that $\Gamma$ is torsion-free, that the composition of the maps in the leftmost column in diagram~\eqref{main diag} is rationally injective, i.e., the composition
$$K_{*+1}(C^*_L(P_d(X))^{\Ga})\xrightarrow{(\beta_L)_*} K_*(C^*_L(P_d(X),\A_{[0,1]}(M))^{\Ga})\xrightarrow{(ev_1)_*} K_*(C^*_L(P_d(X),\A(M^{[0,1]}))^{\Ga,(1)})$$
is rationally injective. The strategy of proof is inspired by \cite{GWY2018}, in which a deformation trick was introduced to show that for a $\Ga$-compact free and proper $\Ga$-space $Z$, the equivariant $KK$-group $KK^{\Ga}(Z,\A)$ stays the same when we continuously deform the $\Gamma$-action on the coefficient algebra $\A$. In particular, $KK^{\Ga,\alpha_1}(X,\A(M^{[0,1]}))$ is isomorphic to $KK^{\Ga,\alpha_0}(X,\A(M^{[0,1]}))$, where the superscript $\alpha_t$ means that $\A(M^{[0,1]})$ is equipped with the $\Ga$-action given by $H(\cdot,t)$ as in Definition \ref{big algebra}. In fact, the sole purpose of introducing in diagram~\eqref{main diag} a middle row involving $\A_{[0,1]}(M)$ is to implement this deformation trick. On a technical level, while in \cite{GWY2018} the deformation trick was described using the language of equivariant $KK$-theory, we shall provide a proof here using the language of equivariant localization algebras \emph{without the $\Ga$-compactness assumption}.

Before we start, we need a K\"unneth formula for $K$-theory of twisted localization algebras. The basic idea of the proof is to follow an outline in \cite{Duke1987}. Thus we shall need some of the functorial properties of localization algebras. However, the proofs of these properties are quite technical and may distract the reader from the main theme of the current article. For this reason, we place these technical arguments in Appendix \ref{Appendix A} and merely summarize the main results here.

\begin{Thm}\label{Kunneth}
Let $X$ be a $CW$-complex with a proper metric and a \emph{trivial group action}. Let $(N,M^{[0,1]},f_0)$ be a coefficient system as we discussed in Section \ref{sec: Twisted algebras and the Bott homomorphisms}, where $f_0:N\to M^{[0,1]}$ is the constant map to a base point $m$ of $M$ viewed as a subspace of $M^{[0,1]}$. Then
$$K_i(C^*_L(X,\A(M^{[0,1]})))\ox\IQ\cong \bigoplus_{j\in\IZ_2\IZ}K_j(X)\ox K_{i-j}(\A(M^{[0,1]}))\ox\IQ,$$
where $C^*_L(X,\A(M^{[0,1]}))$ is the twisted localization algebra defined as in Definition \ref{twisted localization small} with the trivial group action.
\end{Thm}

\begin{proof}
See Appendix \ref{Appendix A} and Theorem \ref{kunneth}.
\end{proof}

We are now ready to introduce a deformation trick for twisted localization algebras. To begin with, we shall need the following definitions.

\begin{Def}\label{Def ulf cover}
Let $\Ga$ be a countable discrete group, $X$ a proper metric space with a proper and isometric $\Ga$-action (\emph{not necessarily free}). Let $\F$ be the family of finite subgroups of $\Ga$.
\begin{itemize}
\item[(1)]A cover $\{U_i\}_{i\in I}$ of $X$ is said to be \emph{equivariant} if for each $\gamma\in\Ga$ and $U_j$, the set $\gamma\cdot U_j$ is still an element in $\{U_i\}_{i\in I}$;
\item[(2)]An equivariant cover $\{U_i\}_{i\in I}$ is called a \emph{$\F$-cover} if for each $U_i$, there exists $\Ga_i\in\F$ such that $U_i$ is $\Ga_i$-invariant and $\Ga\cdot\overline{U_i}$ is homeomorphic to $\Ga\times_{\Ga_i}\overline{U_i}$, where $\Ga\times_{\Ga_i}\overline{U_i}$ means the balanced product (see \cite[Example A.2.6]{HIT2020});\footnote{The notation of $\F$-cover comes from \cite[Definition 1.1]{BLR2008}}
\item[(3)]An $\F$-cover $\{U_i\}_{i\in I}$ is said to be \emph{uniformly locally finite}, if there exists $N>0$ such that $\forall x\in X$,
$$\#\{i\in I\mid x\in\overline{U_i}\}\leq N.$$
\end{itemize}\end{Def}

\begin{Pro}\label{ulf cover}
If $\Ga$ is a countable discrete group, $X$ is a metric space with bounded geometry, and $\Ga\car X$ is proper and isometric, then for each $d\geq 0$, there exists a uniformly locally finite $\F$-cover $\U$ for $P_d(X)$.
\end{Pro}

\begin{proof}
Fix $d\geq 0$. Since $X$ has bounded geometry, the Rips complex $P_d(X)$ is finite-dimensional, say $dim(P_d(X))=N$. Set $\varepsilon=\frac 2{(N+2)^2}>0$.

For each $x_0\in X$, define
$$U_{x_0}=\left\{\sum t_xx\in P_d(X)\mid t_{x_0}>t_x+\varepsilon\text{ for any }x\ne x_0\right\}.$$
We should mention that $U_{x_0}$ is a subset of the open star of $x_0$ in the barycentric subdivision of $P_d(X)$. Set $\U^0=\{U_x\}_{x\in X}$. For each $x_0,x_1\in X$ with $d(x_0,x_1)\leq d$, denote $[x_0,x_1]$ the sub-simplex in $P_d(X)$ generated by $x_0,x_1$. Define
$$U_{[x_0,x_1]}=\left\{\sum t_xx\in P_d(X)\mid \min\{t_{x_0},t_{x_1}\}>t_x+\varepsilon\text{ for any }x\ne x_0,x_1 \right\}.$$
Set $\U^1=\{U_{[x,y]}\}_{[x,y]\subseteq P_d(X)}$. Similarly, for any $n$-subsimplex $\Delta\subseteq P_d(X)$ whose vertex set is $\{x_0,x_1,\cdots,x_n\}\subseteq X$, we define
$$U_{\Delta}=\left\{\sum t_xx\in P_d(X)\mid \min\{t_{x_0},\cdots,t_{x_n}\}>t_x+\varepsilon\text{ for any }x\ne x_0,\cdots,x_{n-1} \right\}.$$
Set $\U^n=\{U_{\Delta}\}_{\Delta\text{ is a $n$-subsimplex of }P_d(X)}$.  Thus we can construct $\U^0,\cdots,\U^N$ by repeating the procedure above.

Set $\U=\U^0\cup\cdots\cup\U^N$. We claim that $\U$ is a uniformly locally finite $\F$-cover for $P_d(X)$. First, for any $z=\sum t_xx\in P_d(X)$, it can be written as
$$z=t_0x_0+t_1x_1+\cdots+t_nx_n+\cdots$$
such that $t_0\geq t_1\geq \cdots\geq t_n> 0=t_{n+1}=\cdots$ and $t_0+\cdots+t_n=1$, where $n\leq N$. Then there exists $m\leq n$ such that $t_m-t_{m+1}>\varepsilon$. Indeed, if $t_m-t_{m+1}\leq\varepsilon$ for all $m\leq n$, then $t_m\leq(n+1-m)\varepsilon$. Then
$$\sum_{m=0}^{n}t_m\leq \frac{(n+1)(n+2)}2\cdot\varepsilon<1$$
which makes a contradiction. Denote $\Delta$ the subsimplex of $P_d(X)$ generated by $\{x_0,\cdots,x_m\}$. Then $z\in U_{\Delta}$ which shows that $\U$ is an open cover.

 For any two different $n$-subsimplexes $\Delta_1,\Delta_2\subseteq P_d(X)$, it is clear that $U_{\Delta_1}\cap U_{\Delta_2}=\emptyset$ since $z=\sum t_xx\in U_{\Delta}$ means that the first $n$ largest coefficients $t_x$ must satisfy that $x$ belongs to the vertex set of $\Delta$. Thus for any $z$, there is at most one element in $\U$ which covers $z$.

For any $U_{\Delta}\in\U$ and $\gamma\in\Ga$, one has $\gamma\cdot U_{\Delta}=U_{\gamma\cdot\Delta}$ by definition. This shows that $U_{\Delta}=\gamma\cdot U_{\Delta}$ if and only if $\Delta=\gamma\cdot\Delta$, which means that $\gamma$ is a torsion element. All such elements form a finite subgroup $\Ga_{\Delta}\leq \Ga$. This means that
$$\Gamma\cdot U_{\Delta}=\Ga\times_{\Ga_{\Delta}} U_{\Delta}.$$
Notice that
$$\overline{U_{\Delta}}=\left\{\sum t_xx\in P_d(X)\mid \min\{t_{x_0},\cdots,t_{x_n}\}\geq t_x+\varepsilon\text{ for any }x\ne x_0,\cdots,x_{n-1} \right\}.$$
where $\{x_0,\cdots,x_n\}$ is the vertex set of $\Delta$. One can use the same argument to show that 
$$\Gamma\cdot\overline{U_{\Delta}}=\Ga\times_{\Ga_{\Delta}}\overline{U_{\Delta}}.$$
As there is at most one element in $\U^m$ which covers $z$ for each $z\in P_d(X)$, thus there are at most $N+1$ elements in $\U$ which cover $z$ since $\U=\U^1\cup\cdots\cup\U^N$. This shows the claim.
\end{proof}

\begin{Rem}
In the above, any element in $\U$ is contractible. For any $U_{\Delta}\in\U$, the barycentric of $\Delta$, denoted by $z_{\Delta}$, clearly belongs to $U_{\Delta}$. One also has that $U_{\Delta}$ is $\Ga_{\Delta}$-equivariantly homotopic to $z_{\Delta}$ by using a linear homotopy. This means that each $U_{\Delta}$ (or $\overline{U_{\Delta}}$) is contractible.
\end{Rem}

\begin{Cons}\label{finite pieces}
Assume that $\Ga\car P_d(X)$ is free for any $d\geq 0$ (this only happens when $\Ga$ is torsion-free). Since $P_d(X)$ admits a uniformly locally finite $\F$-cover, there exists a decomposition $P_d(X)=K_1\cup \cdots\cup K_N$ such that for each $1\leq i\leq N$, there exists a closed subset $K_{i,0}$ such that $K_i=\Ga\cdot K_{i,0}$ and $K_{i,0}$ is homotopic to the disjoint union of at most countable single points. The proof is similar to \cite[Lemma 12.2.3]{HIT2020}. Let $\U$ be a uniformly locally finite $\F$-cover as above and $\overline{\U}=\{\overline{U}\}_{U\in\U}$. Let $C_1$ be a maximal subset of $\U$ which is $\Ga$-equivariant and the orbits of any two elements in $C_1$ have empty intersections. By induction, having defined $C_1,\cdots, C_n$, define $C_{n+1}$ to be a maximal equivariant subset of $\U\backslash(C_1\cup\cdots\cup C_n)$ such that the orbits of any two elements in $C_n$ have empty intersection. This procedure will end for finite steps since $\U$ is uniformly locally finite. Set
$$K_i=\bigcup_{\overline{U}\in C_i}\overline{U}$$
and $K_{i,0}$ to be a fundamental domain for $K_i$. Since $K_i=\bigsqcup_{\Ga\cdot\overline{U}\subseteq C_i} \Ga\cdot\overline{U}$, we can choose a component for each orbit and their union forms a model of $K_{i,0}$. Each $\overline{U_i}$ in $\overline{\U}$ is contractible; it is clear that $K_{i,0}$ is homotopic to the disjoint union of at most countable points.

Actually, for any free $\Ga$-space $X$ which admits a uniformly locally finite $\F$-cover $\U$ such that every element in $\U$ is contractible, the process above still works for $X$. Then $X$ can be viewed as $X=K_1\cup\cdots\cup K_N$ such that each $K_i=\Ga\cdot K_{i,0}$ for some fundamental domain and each $K_{i,0}$ is homotopic to a disjoint union of at most countable points.
\end{Cons}

We should remark that the following lemma relies on the condition that the action $\Ga\car X$ is free.

\begin{Lem}\label{deformation}
Let $X$ be a metric space with a \emph{free} and proper $\Ga$-action which admits a uniformly locally finite $\F$-cover $\U$ such that each element in $\U$ is contractible, $(N,M,f)$ a coefficient system. Then for each $t\in[0,1]$, the $K$-theortic map
$$(ev_t)_*:K_*(C^*_L(X,\A_{[0,1]}(M))^{\Ga})\to K_*(C^*_L(X,\A(M^{[0,1]}))^{\Ga,(t)})$$
is an isomorphism.
\end{Lem}

\begin{proof}
Let $K_1\cdots,K_N\subseteq X$ be as in Construction \ref{finite pieces}, where $K_i=\Ga\cdot K_{i,0}$ is closed for some fundamental domain $K_{i,0}=\bigsqcup_{m\in\IN} K^m_{i,0}$. Here each $K_{i,0}^m$ is contractible to the point $z^m\in K_{i,0}^m$.

Let $(N,M^{[0,1]},f_t)$ be the coefficient system induced by $(N,M,f)$. Then $(N\cap K_i,M^{[0,1]},f_t)$ and $(N\cap K_{i,0},M^{[0,1]},f_t)$ form coefficient systems for $\Ga\car K_i$ and $K_{i,0}$, respectively. Here $K_{i,0}$ is equipped with the trivial group action. Then we can define the twisted localization algebras for $\Ga\car K_i$ and $K_{i,0}$ by Definition \ref{twist localization big}.

{\bf Claim 1.} $K_*(C^*_L(K_i,\A_{[0,1]}(M))^{\Ga})$ is isomorphic to $K_*(C^*_L(K_{i,0},\A_{[0,1]}(M)))$ for each $i$.

One can compare this claim with \cite[Lemma 2.3]{Yu1995}. For any $[g]\in K_1(C^*_L(K_i,\A_{[0,1]}(M))^{\Ga})$ with $g\in(C^*_L(K_i,\A_{[0,1]}(M))^{\Ga})^+$, let $T>0$ be such that $\prop(g(t))\leq\frac{r_0}2$ for any $t>T$, where $r_0=\min_{\gamma\in\Ga} d(K_{i,0},\gamma\cdot K_{i,0})>0$ since $K_{i,0}$ is closed. Then the function $t\mapsto\chi_{K_{i,0}}\cdot g(t+T)\cdot\chi_{K_{i,0}}$ clearly defines an element of $(C^*_L(K_{i,0},\A_{[0,1]}(M)))^+$. Then the construction gives rise to a homomorphism:
$$c:K_1(C^*_L(K_i,\A_{[0,1]}(M))^{\Ga})\to K_1(C^*_L(K_{i,0},\A_{[0,1]}(M))).$$
On the other hand, for any $[h]\in K_1(C^*_L(K_{i,0},\A_{[0,1]}(M)))$, $h^{\Ga}\in (C^*_L(K_i,\A_{[0,1]}(M))^{\Ga})^+$ is defined as follows. For each $t\in\IR_+$ and $(x,y)\in K_i\times K_i$, if $x,y$ are not in the same $\gamma\cdot K_{i,0}$, we define $(h^{\Ga}(t))(x,y)=0$. If $x,y$ both belong to $\gamma\cdot K_{i,0}$ for some $\gamma\in\Ga$, we define
$$(h^{\Ga}(t))(x,y)=\gamma\cdot\left((h(t))(\gamma^{-1}x,\gamma^{-1}y) \right).$$
It is clear that $[h^{\Ga}]$ defines an element in $K_1(C^*_L(K_i,\A_{[0,1]}(M))^{\Ga})$. The correspondence $[h]\mapsto [h^{\Ga}]$ gives the inverse homomorphism of $c$. The case for $K_0$ follows from a suspension argument.

Similarly, one can also prove that
$$K_*(C^*_L(K_i,\A(M^{[0,1]}))^{\Ga,(t)})\cong K_*(C^*_L(K_{i,0},\A(M^{[0,1]}))^{(t)})$$
for each $t\in[0,1]$. Here we omit "$\Ga$" on the right-top corner on the right side since the action on $K_{i,0}$ is trivial. The corner mark $(t)$ only stands for the coarse embedding given by $f_t$.

{\bf Claim 2.} $(ev_t)_*:K_*(C^*_L(K_i,\A_{[0,1]}(M))^{\Ga})\to K_*(C^*_L(K_{i},\A(M^{[0,1]}))^{\Ga,(t)})$ is an isomorphism for each $i$ and $t\in[0,1]$.

By using Claim 1, we have the following commuting diagram:
\[\begin{tikzcd}
K_*(C^*_L(K_i,\A_{[0,1]}(M))^{\Ga}) \arrow[r, "(ev_t)_*"] \arrow[d, "\cong"'] & K_*(C^*_L(K_i,\A(M^{[0,1]}))^{\Ga,(t)}) \arrow[d, "\cong"] \\
K_*(C^*_L(K_{i,0},\A_{[0,1]}(M))) \arrow[r, "(ev_t)_*"]                     & K_*(C^*_L(K_{i,0},\A(M^{[0,1]}))^{(t)}).                   
\end{tikzcd}\]
It suffices to prove the bottom evaluation map for the non-equivariant case is an isomorphism. Since both twisted localization algebras are invariant under homotopy equivalence on the level of $K$-theory (see Appendix \ref{Appendix A}), it suffices to prove the evaluation map is an isomorphism for $\sqcup_{\IN}\{pt\}$. We then have that
\begin{equation}\label{comm5}\begin{tikzcd}
K_*(C^*_L(\sqcup_{m\in\IN}\{z^m\},\A_{[0,1]}(M))) \arrow[r, "(ev_t)_*"] \arrow[d, "(1)","\cong"'] & K_*(C^*_L(\sqcup_{m\in\IN}\{z^m\},\A(M^{[0,1]}))^{(t)}) \arrow[d, "\cong","(2)"'] \\
K_*(C_{ub}(\IR_+,\prod^u_{\IN}\K\ox\A_{[0,1]}(M))) \arrow[r, "(ev_t)_*"] \arrow[d, "ev_*"']  & K_*(C_{ub}(\IR_+,\prod^u_{\IN}\K\ox\A(M^{[0,1]}))) \arrow[d, "ev_*"]  \\
K_*(\prod^u_{\IN}\K\ox\A_{[0,1]}(M)) \arrow[r, "(ev_t)_*"]                     & K_*(\prod^u_{\IN}\K\ox\A(M^{[0,1]})).                  
\end{tikzcd}\end{equation}
Here are some explanations. Let $I:\sqcup_{m\in\IN}\{z^m\}\to K_{i,0}$ be the canonical inclusion where each $K_{i,0}^m\subseteq K_{i,0}$ is contractible to $\{z^m\}$. Then the triple $(\sqcup_{\IN}\{z^m\},M^{[0,1]},f_t\circ I)$ forms a coefficient system for each $t$. Thus the algebras are well-defined. The isomorphism (1) is given by
\[K_*(C^*_L(\sqcup_{m\in\IN}\{z^m\},\A_{[0,1]}(M)))\cong K_*\left(\prod_{m\in\IN}^uC^*_L(\{z^m\},\A_{[0,1]}(M))\right).\]
For any $[g]\in K_*(C^*_L(\sqcup_{m\in\IN}\{z^m\},\A_{[0,1]}(M)))$, since the propagation of $g$ tends to $0$ as $t$ tends to infinity, $g$ can be chosen to be of the form $g=(g_0,g_1,\cdots)$, where each $g_i=\chi_i g\chi_i$ is an element of $C^*_L(\{z^i\},\A_{[0,1]}(M))$ and $\chi_i$ is the characteristic function on $\sqcup\{z^m\}$ associated with $z^i$. The letter $``u"$ in $\prod_{\IN}^u$ means the sequence $(g_0,g_1,\cdots)$ are required to satisfy that there exists $R>0$ such that
$$\supp(ev_t(g_i(s)))\subseteq B((f_t(I(z^i)),0),R)\subseteq M^{[0,1]}\times\IR_+$$
for each $s\in\IR_+$ and $t\in[0,1]$ (i.e., uniformly bounded supported). Since the twisted localization algebra for a single point is isomorphic to the $C^*$-algebra of all bounded and uniformly continuous functions from $\IR_+$ to $\K\ox\A_{[0,1]}(M)$, where $\K$ is the $C^*$-algebra of all compact operators, this shows that the map (1) in the diagram \eqref{comm5} is an isomorphism. Similarly, map (2) in diagram \eqref{comm5} is also an isomorphism.

It is clear that the evaluation-at-zero map from $C_{ub}(\IR_+,\prod^u_{\IN}\K\ox\A_{[0,1]}(M))$ to $\prod^u_{\IN}\K\ox\A_{[0,1]}(M)$ induces an isomorphism on $K$-theory (see \cite[Lemma 12.4.3]{HIT2020}). It suffices to prove the map
$$ev_t:\prod^u_{\IN}\K\ox\A_{[0,1]}(M)\to\prod^u_{\IN}\K\ox\A(M^{[0,1]})$$
induces an isomorphism on $K$-theory for each $t\in[0,1]$. Let $u\in (ker(ev_t))^+$ be an invertible element and we may write $u=(u_0,\cdots,u_m,\cdots)$ such that each $u_m$ is in $(\K\ox\A_{[0,1]}(M))^+$ viewed as a function whose value at $t$ is $1$. Notice that $s\mapsto (u_0(s),\cdots, u_m(s),\cdots)$ may be \emph{not} continuous since the sequence $\{u_m\}_{m\in\IN}$ is \emph{not} equi-continuous, thus we can not directly use a homotopy to show $[u]=0$ in $K_1(ker(ev_t))$. However, we can use a ``stacking argument" and an Eilenberg swindle argument to get around this problem (see \cite[Proposition 12.6.3]{HIT2020} for an example of a stacking argument).

For any invertible element $u=(u_0,\cdots,u_m,\cdots)\in (ker(ev_t))^+$, by condition (3) in Definition \ref{twist Roe big}, set $N_m=[d(f(z_m),m)]+1$ and the sequence of functions 
$$\left\{s\mapsto u_m(N_m^{-1}\cdot s)\right\}_{m\in\IN}$$
is equi-continuous. Extend each $u_m$ to a continuous function on $\IR$ by defining $u_m(s)=u_m(0)$ for $t\leq0$ and $u_m(s)=u_m(1)$ for $t\geq 1$. For each $k\geq 0$, we define
$$u_m^{(k)}(s)=\left\{\begin{aligned}&1,&&s\in\left[t-k\cdot{N_m^{-1}},t+k\cdot{N_m^{-1}}\right];\\
&u_m\left(s+k\cdot{N_m^{-1}}\right),&&s\in\left[0,t-k\cdot{N_m^{-1}}\right];\\
&u_m\left(s-k\cdot{N_m^{-1}}\right),&&s\in\left[t+k\cdot{N_m^{-1}},1\right].\end{aligned}\right.$$
Notice that $u_m^{(0)}=u_m$ and $u_m^{(k)}=1$ for sufficiently large $k$.

Define $\wt{u}=(\wt{u_0},\cdots,\wt{u_m},\cdots)\in (ker(ev_t))^+$, where $\wt{u_m}\in (\K\ox\K\ox\A_{[0,1]}(M))^+\cong(\K\ox\A_{[0,1]}(M))^+$ is given by
$$\wt{u_m}(s)=\bigoplus_{k=1}^{\infty}u_m^{(k)}(s).$$
Write $v=u\oplus\wt u=(v_0,\cdots,v_m,\cdots)\in (ker(ev_t))^+$ and one can see that
$$v_m(s)=\bigoplus_{k=0}^{\infty}u_m^{(k)}(s).$$
It is not hard to see that $v$ is homotopic to $\wt u$ by the map $w:[0,1]\to(ker(ev_t))^+$ defined by $w(r)=(w_1(r),\cdots,w_m(r),\cdots)$, where
$$w_m(r)=\bigoplus_{k=0}^{\infty}u_m^{(k+r)}(s).$$
It is clear that this map is continuous and $w(0)=v$, $w(1)=\wt u$. Thus one then has that $[u]+[\wt u]=[\wt u]$ in $K_1(ker(ev_t))$ for any $[u]\in K_1(ker(ev_t))$. This shows that $K_1(ker(ev_t))=0$ for any $t\in[0,1]$. By using a suspension argument one can see that $K_*(ker(ev_t))=0$. Thus $(ev_t)_*$ is an isomorphism for any $t\in[0,1]$.
This completes the proof of Claim 2.
Since $X=K_1\cup \cdots\cup K_N$, the lemma follows from a Mayer-Vietoris argument and the Five Lemma.
\end{proof}

\begin{Pro}\label{left column}
Let $X$ be a CW-complex with a proper metric and a free and proper $\Ga$-action that admits a uniformly locally finite $\F$-cover $\U$ such that each element in $\U$ is contractible. Then the Bott map
$$(\beta_L)_*:K_{*+1}(C^*_L(X)^{\Ga})\ox\IQ\to K_*(C^*_L(X,\A_{[0,1]}(M))^{\Ga})\ox\IQ$$
is an injection.
\end{Pro}

\begin{proof}
It suffices to prove $(ev_0)_*\circ(\beta_L)_*$ is an injection, where
$$(ev_0)_*:K_*(C^*_L(X,\A_{[0,1]}(M))^{\Ga})\ox\IQ\to K_*(C^*_L(X,\A(M^{[0,1]}))^{\Ga,(0)})\ox\IQ$$
is an isomorphism by Lemma \ref{deformation}. For the case when $t=0$, the group action on $\A(M^{[0,1]})$ is trivial. Thus one can show that
$$K_*(C^*_L(X/\Ga,\A(M^{[0,1]}))^{(0)})\cong K_*(C^*_L(X,\A(M^{[0,1]}))^{\Ga,(0)})$$
by using an argument similar to \cite[Theorem 6.5.15]{HIT2020}, where the coefficient system is given by $(N/\Ga,M^{[0,1]},f^Q_0)$ and $f_0^Q$ is the constant map from $N/\Ga$ to the base point $m$ of $M$.

Combining the K\"unneth formula in Theorem \ref{Kunneth}, we have the following commuting diagram:
\[\begin{tikzcd}
K_{i+1}(C^*_L(X/\Ga))\ox\IQ \arrow[r, "(ev_0)_*\circ(\beta_L)_*"] \arrow[d, "\cong"'] & K_i(C^*_L(X/\Ga,\A(M^{[0,1]}))^{(0)})\ox\IQ\arrow[d, "\cong"] \\
\bigoplus_{j\in\IZ_2\IZ}K_{i-j+1}(C^*_L(X/\Ga))\ox K_j(\S)\ox\IQ  \arrow[r, "1\ox(\beta_m)_*"] & \bigoplus_{j\in\IZ_2\IZ}K_{i-j}(C^*_L(X/\Ga))\ox K_j(\A(M^{[0,1]}))\ox\IQ,                
\end{tikzcd}\]
where the vertical maps are given by the K\"unneth formula, the bottom horizontal map is given by the Bott map $\beta_m:\S\to\A(M^{[0,1]})$ associated with the base point $m\in M\subset M^{[0,1]}$ as in Definition \ref{Bott}. One is referred to \cite{Duke1987} for the K\"unneth formula for $KK$-theory, also see \cite[Proposition 2.11]{TopforCII} for a lovely proof for special cases.
Since $M$ is admissible, we conclude that the bottom map is injective by Theorem \ref{Bott-inj}. This implies that $(ev_0)_*\circ(\beta_L)_*$ is an injection, which completes the proof.
\end{proof}

\section{The twisted assembly map and proof of Theorem \ref{main result torsion-free}}\label{sec: The twisted assembly map and proof of Theorem 1.3}

In this section, we shall complete our half-range target, \Cref{main result torsion-free}, which deals with the equivariant strong coarse Novikov conjecture in the case of torsion-free acting groups. In view of diagram~\eqref{main diag} and what we have done in \Cref{sec: A deformation trick}, it suffices to show map~(11), namely the twisted assembly map defined in \eqref{eq:twisted-assembly}, is an isomorphism. This turns out to be true regardless of the torsion-free assumption. 

To this end, it is natural to combine the ideas in \cite[Theorem 6.8]{Yu2000} for the purely coarse case (i.e., actions by the trivial group) and those in \cite[Theorem 13.1]{GHT2000} for cocompact group actions. More precisely, we shall use a Mayer-Vietoris argument to reduce the problem to verifying the equivariant coarse Baum-Connes conjecture for a sequence of singletons as proper $F$-spaces, where $F$ are finite subgroups of $\Gamma$. In doing so, we need a dimension control condition for $\Ga\car X$ to ensure that the Mayer-Vietoris argument ends in finitely many steps. In \cite{Yu2000} and \cite{GHT2000}, dimension control is guaranteed by the bounded geometry assumption of $X$ and co-compactness of $\Ga\car X$, respectively. In our situation, we shall need the following condition.

\begin{Def}\label{ebg}
Let $X$ be a space with bounded geometry.
The proper group action $\Ga\car X$ is said to have \emph{equivariant bounded geometry} if for any $R>0$ there exists $N>0$ such that 
$$\#\{\gamma\in\Ga\mid B(\gamma\cdot x,R)\cap B(x,R)\ne\emptyset\}<N$$
for all $x\in X$.
\end{Def}

\begin{Rem}\label{Rem ebg}
For a discrete space $X$ with bounded geometry, the following conditions are equivalent:\begin{itemize}
\item[(1)] $\Ga\car X$ has equivariant bounded geometry;
\item[(2)] there exists $N>0$ such that $\#\Ga_x<N$ for all $x\in X$, where $\Ga_x$ is the stabilizer of $x$;
\item[(3)] there exists a free and proper $\Ga$-space $Y$ with bounded geometry and an equivariant coarse equivalence $\phi: Y\to X$.
\end{itemize}
$(1)\Rightarrow(2)$ holds clearly by taking $R=0$. To see $(2)\Rightarrow(1)$, for any $x,y\in X$, if there exists $\gamma\in\Ga$ such that $\gamma x=y$, for any $\gamma'$ satisfying $\gamma'x=y$, there must be $g=\gamma^{-1}\gamma'\in\Ga_x$ such that $\gamma'=\gamma g$. Thus, for any $x,y$,
$$\#\{\gamma\in\Ga\mid \gamma x=y\}\leq N.$$
Since $X$ has bounded geometry, for each $R>0$, there exists $N'>0$ such that $\#B(x,2R)<N'$. Thus one has that
$$\#\{\gamma\in\Ga\mid B(\gamma\cdot x,R)\cap B(x,R)\ne\emptyset\}\leq (N')^2\cdot N.$$
For any $x\in X$, one should notice that $\#\Ga_x=\#\phi^{-1}(x)$. Then the implication $(3)\Rightarrow(2)$ follows from that $\phi$ is an equivariant coarse equivalent and $Y$ has bounded geoemtry. We will show $(2)\Rightarrow (3)$ in Section \ref{sec: A Concrete model for space with Property TAF} (the construction of the model space $M_X$). As a result, if $\Ga$ is torsion-free and $X$ has bounded geometry, $\Ga\car X$ must have equivariant bounded geometry.
\end{Rem}

For the convenience of description, we shall call it a \emph{$\Ga$-slice} an open subset $U\subseteq M^{[0,1]}\times\IR_+$ with the form of a balanced product $U=\Ga\times_{F}U_0$ for some finite subgroup $F\leq \Ga$ and an $F$-invariant open subset $U_0\subseteq M^{[0,1]}$. The following proposition shows the dimension control condition brought by equivariant bounded geometry.

\begin{Pro}\label{equi-bg}
Assume that $\Ga\car X$ has equivariant bounded geometry, and $f_1: X\to M^{[0,1]}$ is the equivariant coarse embedding. For any $R>0$, there exists $N>0$ such that $\Ga\cdot B(f_1(x), R)$ can be written as the union of at most $N$ $\Ga$-slices for all $x\in X$, where
$$B(f_1(x),r)=\{(y,t)\in M^{[0,1]}\times\IR_+\mid d(y,f_1(x))^2+t^2<r^2\}.$$
\end{Pro}

Before starting the proof, we shall recall a notion of \emph{nerve complex}. Let $Y$ be a topological space, $\U=\{U_i\}_{i\in I}$ be a family of locally finite open cover of $Y$. The nerve simplicial complex for $\U$ is a simplicial complex $N(\U)$ whose vertex set is $I$, and a finite subset $J\subseteq I$ spans a simplex if and only if $\bigcap_{j\in J}U_j\ne\emptyset$. Let $\{h_i\}_{i\in I}$ be a partition of unity associated with $\U$. Then there is a canonical continuous map $Y\to N(\U)$ defined by
$$y\mapsto\sum_{i\in I}f_i(y)\cdot i.$$
One can check that this map does not depend on the choice of the partition of unity up to homotopic equivalence. Nerve theorems assert that the homotopy type of a sufficiently nice topological space is encoded in the nerve complex of a good open cover. One is referred to \cite{Hatcher} for more information on nerve complex

Moreover, if $\U=\{U_i\}_{i\in I}$ is a $\Ga$-equivariant open cover, i.e., $\gamma U_i=U_j$ for some $j\in I$, then $I$ is also an $\Ga$-set. This makes $N(\U)$ a $\Ga$-simplicial complex. In this case, if we choose an equivariant partition of unity $\{h_i\}_{i\in I}$ associated with $\U$, the canonical map $Y\to N(\U)$ associated with $\{h_i\}_{i\in I}$ is also an equivariant map.

\begin{proof}[Proof of Proposition \ref{equi-bg}]
Fix $x\in X$. For any $R>0$, the family of open sets
$$\B=\{B(f_1(\gamma\cdot x),R)\mid \gamma\in\Ga\}$$
forms a $\Ga$-equivariant open cover $\Ga\cdot B(f_1(x),R)$, where $B(f_1(x),R)$ is the open set in $M^{[0,1]}\times\IR_+$ defined as above. Let $N(\B)$ be the nerve simplicial complex of $\B$. Since $\Ga\car X$ has equivariant bounded geometry and $f_1$ is a coarse embedding, the nerve complex $N(\B)$ is finite-dimensional, and the dimension of $N(\B)$ is only dependent on $R$.

Notice that the vertex set of $N(\B)$ is the homogeneous space $\Ga/\Ga_x$ of $\Ga$ associated with $\Ga_x$, where $\Ga_x$ is the stabilizer of $x$. The $\Ga$-action on the vertex set of $N(\B)$ is exactly the quasi-regular representation of $\Ga$ on $\Ga/\Ga_x$. Thus $\Ga\car N(\B)$ is cocompact. Indeed, fix a vertex $z_0\in N(\B)$, let $N(\B)_0$ be the subcomplex of $N(\B)$ containing all simplicies in $N(\B)$ with $z_0$ as a vertex.  Since $\Ga\car X$ has equivariant bounded geometry, $N(\B)_0$ is a finite complex and the cardinal number of simplicies in $N(\B)_0$ is only dependent on $R>0$. Then $\Ga\cdot N(\B)_0=N(\B)$. Apply the construction in Proposition \ref{ulf cover} to $N(\B)_0$, one can construct a finite open cover for $N(\B)_0$. As a result, there exists $N>0$ which is only dependent on the dimension of $N(\B)_0$ and the cardinal number of simplicies in $N(\B)_0$ such that $N(\B)_0$ is covered by at most $N$ open sets $\U=\{U_i\}$. Thus $N$ is only dependent on $R$. 

By Proposition \ref{ulf cover}, $\Ga\cdot U_i$ is in the form of a balanced product for each $i$. Since $\Ga\cdot N(B)_0=N(\B)$, $\Ga\cdot\U=\{\Ga\cdot U_i\}$ forms an open cover of $N(\B)$, where each $\Ga\cdot U_i$ takes the form of a balanced product. Denote by $j_{\B}$ the canonical $\Ga$-equivariant continuous map
$$j_{\B}:\Ga\cdot B(f_1(x),R)\to N(\B).$$
Then the family $\{j^{-1}_{\B}(\gamma\cdot U_i)\}_{U_i\in\U,\gamma\in\Ga}$ forms an open cover for $\Ga\cdot B(f_1(x),R)$. Since $j_{\B}$ is $\Ga$-equivariant, this cover is also equivariant satisfying that each $j^{-1}_{\B}(\Ga\cdot U_i)$ being the form of a balanced product for each $i$. This completes the proof.
\end{proof}

The following theorem can be viewed as a twisted version of the coarse Baum-Connes conjecture, which can also be viewed as a non-cocompact version of the generalized Green-Julg Theorem.

\begin{Thm}\label{twist assembly}
Let $(N, M,f)$ be a Roe algebraic coefficient system for $\Ga\car X$, where $\Ga$ acts on $X$ properly by isometries. If $f$ is an equivariant coarse embedding and $\Ga\car X$ has equivariant bounded geometry, then the twisted assembly map is an isomorphism, i.e.,
$$ev_*:\lim_{d\to\infty}K_*(C^*_L(P_d(X),\A(M^{[0,1]}))^{\Ga,(1)})\to\lim_{d\to\infty}K_*(C^*(P_d(X),\A(M^{[0,1]}))^{\Ga,(1)})$$
is an isomorphism.
\end{Thm}

\begin{Rem}\label{difference with FW}
It is studied in \cite{FW2016} for the case when $M$ is a separable Hilbert space by using a similar method to ours. Theorem \ref{twist assembly} above is similar to \cite[Theorem 4.11]{FW2016}. However, there is a subtle gap in their proof, more concretely, the gap is in the proof of \cite[Lemma 4.20]{FW2016}. One should need the family $\{\Delta_{i,j}\}_{j\in J}$ to be disjoint to make sure the family is homotopy equivalent to $\{y_{i,j}\}_{j\in J}$ in their paper. This does not happen when the parameter $d$ in $P_d(X)$ tends large enough. To solve this problem, we shall need to cut the twisted algebras by using a much finer family of open sets.
\end{Rem}

To begin with, we shall first introduce an ideal of twisted algebras associated with a given open set in $M^{[0,1]}$.

\begin{Def}
Let $O$ be a $\Ga$-invariant open subset of $M^{[0,1]}\times\IR_+$. Define $C^*(P_d(X),\A(M^{[0,1]}))^{\Ga,(1)}_O$ to be the subalgebra of $C^*(P_d(X),\A(M^{[0,1]}))^{\Ga,(1)}$ generated by all $T\in\IC[P_d(X),\A(M^{[0,1]})]^{\Ga,(1)}$ such that
$$\supp(T(x,y))\subseteq O$$
for all $x,y\in N_d$.

Define $C^*_L(P_d(X),\A(M^{[0,1]}))^{\Ga,(1)}_O$ to be the subalgebra of $C^*_L(P_d(X),\A(M^{[0,1]}))^{\Ga,(1)}$ generated by all $g\in\IC_L[P_d(X),\A(M^{[0,1]})]^{\Ga,(1)}$ such that
$$\supp((g(t))(x,y))\subseteq O$$
for all $x,y\in N_d$ and $t\in\IR_+$.
\end{Def}

Notice that $C^*(P_d(X),\A(M^{[0,1]}))^{\Ga,(1)}_O$ and $C^*_L(P_d(X),\A(M^{[0,1]}))^{\Ga,(1)}_O$ are two-sided closed ideals of $C^*(P_d(X),\A(M^{[0,1]}))^{\Ga,(1)}$ and $C^*_L(P_d(X),\A(M^{[0,1]}))^{\Ga,(1)}$, respectively. There still is an evaluation map
$$ev:C^*_L(P_d(X),\A(M^{[0,1]}))^{\Ga,(1)}_O\to C^*(P_d(X),\A(M^{[0,1]}))^{\Ga,(1)}_O$$
defined by
$$ev(g)=g(0).$$

\begin{Lem}\label{iso for slice}
If $O=\bigsqcup_{i\in J}O_i$ is a disjoint union of $\{O_i\}$ in $M^{[0,1]}\times\IR_+$ such that\begin{itemize}
\item[(1)] $O_i$ is a $\Ga$-slice, i.e., $O_i=\Ga\times_{F_i}O_{i,0}$ for some finite subgroup $F_i\leq\Ga$ and an $F_i$-invariant open set $O_{i,0}\subseteq M^{[0,1]}\times\IR_+$;
\item[(2)] for each $i$, $\exists r>0$, $x_i\in X$ such that
$$O_{i,0}\subseteq B(f_1(x_i),r)=\{(y,t)\in M^{[0,1]}\times\IR_+\mid d(y,f_1(x_i))^2+t^2<r^2\}$$
where $f_1:X\to M^{[0,1]}$ is the coarse embedding.
\end{itemize}
Then
$$ev_*:\lim_{d\to\infty}K_*\left(C^*_L(P_d(X),\A(M^{[0,1]}))^{\Ga,(1)}_O\right)\to\lim_{d\to\infty}K_*\left(C^*(P_d(X),\A(M^{[0,1]}))^{\Ga,(1)}_O\right)$$
is an isomorphism.
\end{Lem}

One should compare condition (1) in Lemma \ref{iso for slice} with the proof of \cite[Theorem 13.1]{GHT2000} for the cocompact case (also see \cite[Proposition 3.3]{GW2022} for our motivation). We still need some preparations before we can prove Lemma \ref{iso for slice}. Let $O=\bigsqcup_{i\in J} O_i$ and $\{x_i\}_{i\in J}$ be as in Lemma \ref{iso for slice}. Let $\{X_i\}_{i\in J}$ be a family of subsets of $P_d(X)$ such that $X_i$ is $F_i$-invariant and $\exists R>0$ such that $X_i\subseteq B(x_i,R)$ for all $i\in J$. Typical examples for $\{X_i\}$ include the following cases:\begin{itemize}
\item[(1)] $X_i=F_i\cdot B_{P_d(X)}(x_i,R)=\{\gamma\cdot z\in P_d(X)\mid d(z,x_i)\leq R, \gamma\in F_i\}$ for some common $R>0$ for all $i\in J$;
\item[(2)] $X_i=\Delta_i(R)$, a simplex in $P_d(X)$ with vertex set $\{\gamma\cdot x\mid \gamma\in F_i,x\in X, d(x,x_i)\leq R\}$ for sufficient large $d$.
\end{itemize}
Notice that $\Ga\car X$ has equivariant bounded geometry, the family of subgroup $\{F_i\}_{i\in J}$ is uniformly finite, i.e., $\sup_{i\in J}\#F_i<\infty$. Indeed, since $\gamma\cdot O_{i,0}=O_{i,0}\subseteq B(f_1(x_i).r)$ for all $\gamma\in F_i$, there exists $L>0$ decided by $r$ and the coarse embedding $f_1$ such that 
$$F_i\subseteq \{\gamma\in\Ga\mid B(\gamma\cdot x_i,L)\cap B(x_i,L)\ne\emptyset\}.$$
Thus, for sufficiently large $d\geq 0$, the simplex $\Delta_i(R)\subseteq P_d(X)$  for each $i\in J$.

For an open set $O\subseteq M^{[0,1]}\times\IR_+$, let $\A(M^{[0,1]})_{O}$ be the $C^*$-subalgebra of $\A(M^{[0,1]})$ generated by elements whose supports are contained in $O$. Notice that $\A(M^{[0,1]})_{O_{i,0}}$ is an $F_i$-algebra instead of a $\Ga$-algebra for each $i\in J$.

\begin{Def}
Define
\[A^*(X_i,i\in J)=\prod_{i\in J} C^*(X_i)^{F_i}\ox\A(M^{[0,1]})_{O_{i,0}}.\]
Similarly, define $A_L^*(X_i,i\in J)$ to be the $C^*$-subalgebra of $\prod_{i\in J}C^*_L(X_i)\ox\A(M^{[0,1]})_{O_{i,0}}$ generated by all elements $\oplus_{i\in J}b_i$ such that\begin{itemize}
\item[(1)] the sequence $\{b_i\}_{i\in J}$ is equicontinuous by norm, i.e., the map $t\mapsto\oplus b_i(t)\in A^*(X_i,i\in J)$ is uniformly norm-continuous;
\item[(2)] there exists a bounded function $c:\IR_+\to\IR_+$ with $\lim_{t\to\infty}c(t)=0$ such that $(b_i(t))(x,y)=0$ whenever $d(x,y)>c(t)$ for all $i\in J$, $x,y\in X_i$ and $t\in\IR_+$.
\end{itemize}\end{Def}

\begin{Lem}\label{to simplex}
Let $O=\bigsqcup_{i\in J}O_i$ be the disjoint union of the family of $\Ga$-equivariant open subsets of $M^{[0,1]}\times\IR_+$ as in Lemma \ref{iso for slice}. Then we have the following isomorphisms:\begin{itemize}
\item[(1)] $C^*(P_d(X),\A(M^{[0,1]}))^{\Ga,(1)}_O\cong\lim\limits_{R\to\infty} A^*(F_i\cdot B_{P_d(X)}(x_i,R):i\in I)$;
\item[(2)] $C^*_L(P_d(X),\A(M^{[0,1]}))^{\Ga,(1)}_O\cong\lim\limits_{R\to\infty} A^*_L(F_i\cdot B_{P_d(X)}(x_i,R):i\in I)$;
\item[(3)] $\lim\limits_{d\to\infty}C^*(P_d(X),\A(M^{[0,1]}))^{\Ga,(1)}_O\cong\lim\limits_{R\to\infty} A^*(\Delta_i(R):i\in I)$;
\item[(4)] $\lim\limits_{d\to\infty}C^*_L(P_d(X),\A(M^{[0,1]}))^{\Ga,(1)}_O\cong\lim\limits_{R\to\infty} A^*_L(\Delta_i(R):i\in I)$;
\end{itemize}\end{Lem}

\begin{proof}
We shall only prove the first term and (2)-(4) follow from (1) directly. One is referred to \cite[Lemma 6.4]{Yu2000} for the non-equivariant case, \cite[Proposition 12.9]{GHT2000} and \cite[Lemma 3.1]{GW2022} for cocompact action case.

Define $E_O=\ell^2(N_d)\ox\ell^2(\Ga)\ox\H_0\ox\A(M^{[0,1]})_O$ which is a Hilbert $\A(M^{[0,1]})_O$-module and also a submodule of $E_{\A(M^{[0,1]})}$ as a Hilbert $\A(M^{[0,1]})$-module. Clearly, $C^*(P_d(X),\A(M^{[0,1]}))^{\Ga,1}_{O}$ has a faithful representation on $E_O$. By condition (1) in Lemma \ref{iso for slice}, one has that
$$\A(M^{[0,1]})_{O}\cong\prod_{i\in J}\A(M^{[0,1]})_{O_i}\cong\prod_{\gamma\cdot O_{i,0}\subseteq O}\A(M^{[0,1]})_{\gamma\cdot O_{i,0}}.$$
As a result, we also have a decomposition
$$E_O=\bigoplus_{\gamma\cdot O_{i,0}\subseteq O}E_{\gamma\cdot O_{i,0}}.$$
Thus, for each $T\in C^*(P_d(X),\A(M^{[0,1]}))^{\Ga}$, there is a corresponding decomposition
$$T=\bigoplus_{\gamma\cdot O_{i,0}\subseteq O}T_{\gamma,i},$$
where
$$T_{\gamma,i}(x,y)=T(x,y)|_{\gamma\cdot O_{i,0}}\in\A(M^{[0,1]})_{\gamma\cdot O_{i,0}}.$$
Since $T$ is $\Ga$-equivariant, for any $\gamma\in\Ga$ and $x,y\in O_{i,0}$, one has that
$$T_{\gamma,i}(\gamma x,\gamma y)=\gamma\cdot T_{e,i}(x,y).$$

We claim that $\bigoplus_{i\in J}T_{e,i}\in A^*(F_i\cdot B_{P_d(X)}(x_i,R):i\in I)$ for sufficiently large $R>0$. First of all, since $T$ is $\Ga$-equivariant and $O_{i,0}$ is $F_i$-invariant, $T_{e,i}(\gamma x,\gamma y)=\gamma\cdot T_{e,i}(x,y)$ for all $\gamma\in F_i$. This means that $T_{e,i}$ is $F_i$-invariant.
By condition (1) of Definition \ref{twisted Roe small}, there exists $R_1>0$ such that $T_{e, i}(x,y)=0$ if $B(f_1(x), R_1)\cap O_{i,0}=\emptyset$ for each $i\in J$. Since $f_1$ is a coarse embedding, there exists $R>0$ such that $T_{e,i}(x,y)=0$ for all $x,y\notin B_{P_d(X)}(x_i,R)$. Thus $T_{e,i}$ can be viewed as an element of $C^*(F_i\cdot B_{P_d(X)}(x_i,R))^{F_i}\ox\A(M^{[0,1]})_{O_{i,0}}$ since $C^*(F_i\cdot B_{P_d(X)}(x_i,R))^{F_i}\ox\A(M^{[0,1]})_{O_{i,0}}$ also has a faithful representation on $E_{e,i}$ and fits the action of $T_{e,i}$ on $E_{e,i}$. This shows the claim. Then the map
$$T\mapsto \oplus_{i\in J}T_{e,i}$$
gives a $C^*$-homomorphism.

For the inverse, let $\oplus_{i\in J}T_i\in A^*(F_i\cdot B_{P_d(X)}(x_i,R):i\in I)$. One can first define $T_{\gamma,i}$ by 
$$T_{\gamma,i}(x,y)=\gamma\cdot T_i(\gamma^{-1}x,\gamma^{-1}y)$$
for each $[\gamma]\in\Ga/F_i$. Define $\widetilde{T_i}=\oplus_{\gamma\in \Ga}T_{\gamma,i}$. Then
$$\oplus_{i\in J}T_i\mapsto \oplus_{i\in J}\widetilde{T_i}$$
gives the inverse homomorphism of the map above, thus we have the isomorphism
$$C^*(P_d(X),\A(M^{[0,1]}))^{\Ga,(1)}_O\cong\lim\limits_{R\to\infty} A^*(F_i\cdot B_{P_d(X)}(x_i,R):i\in I)$$
as desired in (1).
\end{proof}

It is known that the equivariant localization algebras are invariant under equivariant strongly Lipschitz homotopy equivalence on the level of $K$-theory (see \cite[Lemma 6.5]{Yu2000}, also see \cite[Lemma 4.19]{FW2016}). One can follow a similar argument to show the algebras $A^*_L(X_i: i\in I)$ we defined above are also invariant under a similar setting as follows. We say that a map $g:\bigsqcup_{i\in I} X_i\to\bigsqcup_{i\in I} Y_i$ is said to be equivariant if $g|_{X_i}$ is $F_i$-equivariant. An equivariant coarse map $g:\bigsqcup_{i\in I} X_i\to\bigsqcup_{i\in I} Y_i$ is said to be Lipschitz if $g(X_i)\subset Y_i$ and $d(g(x),g(y))\leq c d(x,y)$. Two equivariant Lipschitz maps $g,h$ from $\bigsqcup_{i\in I} X_i\to\bigsqcup_{i\in I} Y_i$ are said to be equivariantly strongly Lipschitz homotopy equivalent if there exists a continuous homotopy $H(x,t)(t\in [0,1])$ satisfying that
\begin{itemize}
\item[(1)]$d(H(x,t),H(y,t))\leq Cd(x,y)$ for all $x,y\in X_i$, $i\in I$ and $t\in[0,1]$, where $C$ is a constant (called the Lipschitz constant of $H$);
\item[(2)]$H$ is equi-continuous in $t$, i.e., for any $\varepsilon>0$, there exists $\delta>0$ such that $d(H(x,t_1),H(x,t_2))\leq\varepsilon$ for all $x\in \bigsqcup_{i\in I}X_i$ if $|t_1-t_2|<\delta$;
\item[(3)]$H(x,0)=g(x),H(x,1)=h(x)$ for all $x\in \bigsqcup_{i\in I}X_i$ and $H(x,t)$ is equivariant and coarse on $x$ for each $t\in[0,1]$.
\end{itemize}
The space $\bigsqcup_{i\in I}X_i$ is equivariantly strongly Lipschitz homotopy equivalent to $(Y_i)_{i\in I}$ if there exist equivariant Lipschitz maps
$$g:\bigsqcup_{i\in I}X_i\to\bigsqcup_{i\in I}Y_i$$
and
$$h:\bigsqcup_{i\in I}Y_i\to\bigsqcup_{i\in I}X_i$$
such that $g\circ h$ and $h\circ g$ are equivariantly strongly Lipschitz homotopy equivalent to the identity maps, respectively.

\begin{Lem}\label{to single point}
Let $O=\{O\}_{i\in J}$ be the family of open sets as in Lemma \ref{iso for slice}, $\Delta_i(R)\subseteq P_d(X)$ as above. Then the evaluation map induces an isomorphism on $K$-theory
$$ev_*:K_*(A^*_L(\Delta_i(R),i\in I))\to K_*(A^*(\Delta_i(R),i\in I))$$
for each $R>0$.
\end{Lem}

\begin{proof}
The family $\{\Delta_i(R)\}_{i\in J}$ is equivariantly strongly Lipschitz homotopy equivalent to $\{z_i\}_{i\in J}$ by using a linear homotopy, where $z_i$ is the barycenter of $\Delta_i(R)$. On the other hand, the canonical inclusions $\{z_i\to \Delta_i(R)\}_{i\in J}$ are uniformly equivariant coarse equivalences for each $i\in J$. Thus, we have the following commuting diagram
$$\begin{tikzcd}
{K_*(A^*_L(\Delta_i(R),i\in I))}  \arrow[r, "ev_*"] & {K_*(A^*(\Delta_i(R),i\in I))} \\
{K_*(A^*_L(\{z_i\},i\in I))} \arrow[u, "\cong"]\arrow[r, "ev_*"]    & {K_*(A^*(\{z_i\},i\in I))\arrow[u, "\cong"'] .}                       
\end{tikzcd}$$
Notice that $A^*_L(\{z_i\},i\in I)\cong C_{ub}(\IR_+,A^*(\{z_i\},i\in I))$, the $C^*$-algebra of all bounded, uniformly continuous functions from $\IR_+$ to $A^*(\{z_i\},i\in I)$. It suffices to show the evaluation-at-zero map
$$ev:C_{ub}(\IR_+,A^*(\{z_i\},i\in I))\to A^*(\{z_i\},i\in I)$$
induces an isomorphism on $K$-theory. Then the lemma follows from \cite[Lemma 12.4.3]{HIT2020}.
\end{proof}

\begin{proof}[Proof of Lemma \ref{iso for slice}]
It follows from Lemma \ref{to simplex} and Lemma \ref{to single point}.
\end{proof}

Before we can prove Theorem \ref{twist assembly}, we still need the following lemma.

\begin{Lem}\label{cutting and pasting}
Equip $M^{[0,1]}\times\IR_+$ with the weak topology given by the center of $\A_{ev}(M^{[0,1]})$. Let $O^{(1)}$ and $O^{(2)}$ be $\Ga$-invariant open sets of $M^{[0,1]}\times\IR_+$. Then
\begin{itemize}
\item[(1)]$C^*(P_d(X),\A(M^{[0,1]}))^{\Ga,(1)}_{O^{(1)}}+C^*(P_d(X),\A(M^{[0,1]}))^{\Ga,(1)}_{O^{(2)}}=C^*(P_d(X),\A(M^{[0,1]}))^{\Ga,(1)}_{O^{(1)}\cup O^{(2)}}$
\item[(2)]$C^*(P_d(X),\A(M^{[0,1]}))^{\Ga,(1)}_{O^{(1)}}\cap C^*(P_d(X),\A(M^{[0,1]}))^{\Ga,(1)}_{O^{(2)}}=C^*(P_d(X),\A(M^{[0,1]}))^{\Ga,(1)}_{O^{(1)}\cap O^{(2)}}$
\item[(3)]$C^*_L(P_d(X),\A(M^{[0,1]}))^{\Ga,(1)}_{O^{(1)}}+C^*_L(P_d(X),\A(M^{[0,1]}))^{\Ga,(1)}_{O^{(2)}}=C^*_L(P_d(X),\A(M^{[0,1]}))^{\Ga,(1)}_{O^{(1)}\cup O^{(2)}}$
\item[(4)]$C^*_L(P_d(X),\A(M^{[0,1]}))^{\Ga,(1)}_{O^{(1)}}\cap C^*_L(P_d(X),\A(M^{[0,1]}))^{\Ga,(1)}_{O^{(2)}}=C^*_L(P_d(X),\A(M^{[0,1]}))^{\Ga,(1)}_{O^{(1)}\cap O^{(2)}}$
\end{itemize}\end{Lem}

\begin{proof}
The proof is similar to \cite[Lemma 6.3]{Yu2000}. There is a slight difference we need to explain. The space $M^{[0,1]}\times\IR_+$ with weak topology is a locally compact, Hausdorff space since it is the spectrum of $\A_{ev}(M^{[0,1]})$. Thus, we can directly use an equivariant partition of unity $\{\phi_1,\phi_2\}$ associated with $\{O^{(1)}, O^{(2)}\}$ instead of $g_x$ and $g_x'$ in \cite[Lemma 6.3]{Yu2000}. We leave the details to the reader.
\end{proof}

\begin{proof}[Proof of Theorem \ref{twist assembly}]
Since $X$ has bounded geometry, for any $s>0$, there exists $N>0$ and a decomposition
$$X=\bigsqcup_{n=1}^NX_n$$
such that each $X_n$ is $\Ga$-invariant and for each $x,y\in X_n$ with $\Ga\cdot x\cap\Ga\cdot y=\emptyset$, $d(\Ga\cdot x,\Ga\cdot y)>s$. One can construct the decomposition by using an induction argument similar to \cite[Lemma 12.2.3]{HIT2020}. As a corollary, for each $R>0$, the family $\{B(f_1(x),R)\}_{x\in X}$ has also a decomposition
$$\bigcup_{x\in X}B(f_1(x),R)=\bigcup_{n=1}^{N_1}\bigsqcup_{\Ga\cdot x\subseteq X_n}\Ga\cdot B(f_1(x),R).$$
For each orbit $\Ga\cdot x$, the open set $\Ga\cdot B(f_1(x),R)$ has a decomposition
$$\Ga\cdot B(f_1(x),R)=\bigcup_{k=1}^{N_2}O_{x,k}(R)$$
such that each $O_{x,k}(R)$ is in the form of a balanced product by using Proposition \ref{equi-bg}. Then for each $k\in\{1,\cdots,N_2\}$, let
$$O^k(R)=\bigsqcup_{x\in X_n}O_{x,k}(R).$$
By Lemma \ref{iso for slice}, the evaluation map induces an isomorphism on $K$-theory:
$$ev_*:\lim_{d\to\infty}K_*\left(C^*_L(P_d(X),\A(M^{[0,1]}))^{\Ga,(1)}_{O^k(R)}\right)\xrightarrow{\cong}\lim_{d\to\infty}K_*\left(C^*(P_d(X),\A(M^{[0,1]}))^{\Ga,(1)}_{O^k(R)}\right).$$
We denote 
$$B_n(R)=\bigsqcup_{\Ga\cdot x\subseteq X_n}\Ga\cdot B(f_1(x),R)\quad\text{and}\quad B(R)=\bigcup_{x\in X}B(f_1(x),R)$$
for convenience. By definition, $B_n(R)=\bigcup_{k=1}^{N_2}O^{k(R)}$. By using Lemma \ref{cutting and pasting} and a Mayer-Vietoris argument, we have that
$$ev_*:\lim_{d\to\infty}K_*\left(C^*_L(P_d(X),\A(M^{[0,1]}))^{\Ga,(1)}_{B_n(R)}\right)\xrightarrow{\cong}\lim_{d\to\infty}K_*\left(C^*(P_d(X),\A(M^{[0,1]}))^{\Ga,(1)}_{B_n(R)}\right)$$
is an isomorphism for each $n$ and $R>0$. Since $B(R)=\bigcup_{i=1}^{N_1}B_n(R)$, by using a Mayer-Vietoris argument again, we have that
$$ev_*:\lim_{d\to\infty}K_*\left(C^*_L(P_d(X),\A(M^{[0,1]}))^{\Ga,(1)}_{B(R)}\right)\xrightarrow{\cong}\lim_{d\to\infty}K_*\left(C^*(P_d(X),\A(M^{[0,1]}))^{\Ga,(1)}_{B(R)}\right)$$
is an isomorphism for each $R>0$. By definition, we have the following commuting diagram:
$$\begin{tikzcd}
\lim\limits_{d\to\infty}K_*\left(C^*_L(P_d(X),\A(M^{[0,1]}))^{\Ga,(1)}\right) \arrow[r,"ev_*"] \arrow[d,"\cong"'] & 
\lim\limits_{d\to\infty}K_*\left(C^*_L(P_d(X),\A(M^{[0,1]}))^{\Ga,(1)}\right) \arrow[d,"\cong"] \\
\lim\limits_{d\to\infty}\lim\limits_{R\to\infty}K_*\left(C^*_L(P_d(X),\A(M^{[0,1]}))^{\Ga,(1)}_{B(R)}\right) \arrow[r,"ev_*"] \arrow[d,"\cong"'] & 
\lim\limits_{d\to\infty}\lim\limits_{R\to\infty}K_*\left(C^*_L(P_d(X),\A(M^{[0,1]}))^{\Ga,(1)}_{B(R)}\right) \arrow[d,"\cong"] \\
\lim\limits_{R\to\infty}\lim\limits_{d\to\infty}K_*\left(C^*_L(P_d(X),\A(M^{[0,1]}))^{\Ga,(1)}_{B(R)}\right) \arrow[r,"ev_*","\cong"']           & 
\lim\limits_{R\to\infty}\lim\limits_{d\to\infty}K_*\left(C^*_L(P_d(X),\A(M^{[0,1]}))^{\Ga,(1)}_{B(R)}\right)         
\end{tikzcd}$$
The upper vertical maps are isomorphisms by definition and the lower vertical maps are isomorphisms by \cite[Theorem 3.3]{SW2007}. The bottom horizontal isomorphism is proved as above. We then complete the proof.
\end{proof}

\subsection*{Proof of Theorem \ref{main result torsion-free}}

Finally, we complete the proof of Theorem \ref{main result torsion-free}.

\begin{proof}[Proof of Theorem \ref{main result torsion-free}]
Consider the commuting diagram \eqref{com diag}:
\[\begin{tikzcd}
& K_{*+1}(C^*_L(P_d(X))^{\Ga}) \arrow[r, "ev_*"] \arrow[d, "(\beta_L)_*"] & K_{*+1}(C^*(P_d(X))^{\Ga}) \arrow[d, "\beta_*"]   \\
& K_*(C^*_L(P_d(X),\A_{[0,1]}(M))^{\Ga})\arrow[r, "ev_*"] \arrow[d, "(ev_1)_*"]& K_*(C^*(P_d(X),\A_{[0,1]}(M))^{\Ga}) \arrow[d, "(ev_1)_*"]\\
& K_*(C^*_L(P_d(X),\A(M^{[0,1]}))^{\Ga,(1)})\arrow[r, "ev_*"] & K_*(C^*(P_d(X),\A(M^{[0,1]}))^{\Ga,(1)})
\end{tikzcd}\]
By Proposition \ref{ulf cover} and Construction \ref{finite pieces}, we have that Lemma \ref{deformation} and Lemma \ref{left column} hold for Rips complex of $X$ at any scale $d>0$. Tracing along the left column and the bottom row, we see that $(ev_1)_*\circ(\beta_L)_*$ is an injection by Lemma \ref{deformation} and Lemma \ref{left column}, the bottom map is an isomorphism by Theorem \ref{twist assembly}. This implies the assembly map is $\IQ$-injective, which completes the proof.
\end{proof}

\section{The Milnor-Rips complexes}\label{sec: The Milnor-Rips complex}

Now we enter Part~II of the paper, where we shift our focus to the rational analytic equivariant coarse Novikov conjecture beyond the case of a torsion-free acting group. 
While it appears difficult to tackle the equivariant strong coarse Novikov conjecture (a stronger statement) with our method because we do not know how to prove a version of \Cref{deformation} for a possibly non-free $\Ga$-space $X$, fortunately, this is not a problem for the rational analytic equivariant coarse Novikov conjecture, since, roughly speaking, the classifying spaces involved in this conjecture are concerned only with free and proper $\Gamma$-spaces (that are equivariantly coarsely equivalent to a given proper $\Gamma$-space). 

In this section, we will construct a concrete model for such a classifying space and discuss the topological side of the analytic equivariant coarse Novikov conjecture. Our presentations here are inspired by Sections~7.1 and~7.2 of \cite{HIT2020}.

\subsection{The Rips complexes and universal $K^{\Ga}$-group}

Let $X$ be a discrete metric space with bounded geometry and proper isometric action of a countable discrete group $\Ga$. Write $\P(X,\Ga)$ for the category with objects including all pairs $(Z,p_Z)$ where
\begin{itemize}
\item[-] $Z$ is a proper metric space with bounded geometry and a proper and isometric $\Ga$-action;
\item[-] $p_Z:Z\to X$ is an equivariant coarse equivalence.
\end{itemize}
A morphisms in $\P(X,\Ga)$ between $(Z,p_Z)$ and $(Y,p_Y)$ is a equivariant continuous coarse map $f:Z\to Y$ such that the following diagram commutes up to closeness:
$$\begin{tikzcd}
Z \arrow[r, "f"] \arrow[d, "p_Z"'] & Y \arrow[d, "p_Y"] \\
X \arrow[r, Rightarrow, no head]   & X.                
\end{tikzcd}$$

\begin{Def}
With notations as above,
\begin{itemize}
\item[(1)]a group $G$ is said to be a \emph{$K^{\Ga}$-group} for $\P(X,\Ga)$ if for any object $Z$ of $\P(X,\Ga)$ there exists a $c_Z:K^{\Ga}_*(Z)\to G$ such that for any morphism $f:Z\to Y$, the following diagram commutes:
$$\begin{tikzcd}
K^{\Ga}_*(Z) \arrow[r, "f_*"] \arrow[d, "c_Z"'] & K_*^{\Ga}(Y) \arrow[d, "c_Y"] \\
G \arrow[r, Rightarrow, no head]   & G .                
\end{tikzcd}$$

\item[(2)]a \emph{universal $K^{\Ga}$-group} for $\P(X,\Ga)$, denoted by $KX^{\Ga}_*(X)$, is a $K^{\Ga}_*$-group satisfying that  for any $K^{\Ga}_*$-group $G$ there exists a group homomorphism $\mu:KX_*(X)\to G$ such that the following diagram commutes:
$$\begin{tikzcd}
K^{\Ga}_*(Z) \arrow[r, "c_Z"] \arrow[d,  Rightarrow, no head] & KX_*^{\Ga}(X) \arrow[d, "\mu"] \\
K_*^{\Ga}(Z) \arrow[r, "d_X"]              & G   
\end{tikzcd}$$
for any $Z\in \P(X,\Ga)$, where $c_Z$ and $d_Z$ come from the definition of $K^{\Ga}$-group.
\end{itemize}\end{Def}

For each $Z\in\P(X,\Ga)$, we have the index map:
$$\ind_Z:K_*^{\Ga}(Z)\to K_*(C^*(Z)^{\Ga})\cong K_*(C^*(X)^{\Ga}),$$
where the index map is defined as in Section \ref{sec: The equivariant coarse Novikov conjecture}. One can check that $K_*(C^*(X)^{\Ga})$ is a $K^{\Ga}_*$-group associated with the family $\{\ind_Z\mid Z\in\P(X,\Ga)\}$. The universal property gives a group homomorphism which we call it assembly map:
$$\mu_{X}^{\Ga}:KX^{\Ga}_*(X)\to K_*(C^*(X)^{\Ga}).$$

The universal group $KX^{\Ga}_*(X)$ is made to assemble all the indices of elliptic operators \cite{BCH1994}. Such a group always exists \cite[Lemma 7.1.6]{HIT2020} and is unique unto canonical isomorphism from the universal property. When $\Ga\car X$ is cocompact, $KX_*^{\Ga}(X)$ is isomorphic to $RK^{\Ga}_*(\EG)$ where $\EG$ is the classifying space for proper $\Ga$-actions. If we moreover assume that $\Ga$ is torsion-free, then $KX_*^{\Ga}(X)$ is isomorphic to $RK^{\Ga}_*(\E\Ga)$ where $\E\Ga$ is the universal principal $\Ga$-bundle (cf. \cite{GTM20}).

\begin{ECBCcon}
The assembly map $\mu_{X}^{\Ga}$ is an isomorphism.
\end{ECBCcon}

In other words, the equivariant coarse Baum-Connes conjecture claims that any element in $K_*(C^*(X)^{\Ga})$ is the index of an elliptic operator, which is unique up to equivalence in $KX^{\Ga}_*(X)$. The injectivity of $\mu_X^{\Ga}$ (i.e., the uniqueness part) implies the Novikov conjecture, so we always call the injectivity part of $\mu_X^{\Ga}$ the equivariant coarse Novikov conjecture:

\begin{ECNcon}
The assembly map $\mu_{X}^{\Ga}$ is an injection.
\end{ECNcon}

There are two ways to approach a concrete model for the universal $K^{\Ga}_*$-group for $\P(X,\Ga)$, which we will introduce in this section. One is to use the \emph{join construction} introduced by J.~Milnor in \cite{Milnor56}. The other one is to use the \emph{Rips complexes} as we introduced in Section \ref{sec: The equivariant coarse Novikov conjecture}. 
The following lemma shows the existence of the equivariant partition of unity that one can similarly prove as \cite[Corollary A.2.8]{HIT2020} or \cite[Proposition 4.12.1]{GTM20}:

\begin{Lem}\label{partition of unity}
Let $\Ga$ be a countable discrete group acting properly by isometries on a proper metric space $X$. Then there exists a countable locally finite $\Ga$-invariant subset $W\subseteq X$ such that:\begin{enumerate}
\item[(1)] $W$ is a $r$-net in $X$ for some $r>0$;
\item[(2)] for each $w\in W$, there exists $0<r_w<r$ and a finite subgroup $\Ga_w\leq\Ga$ such that $\Ga\cdot B(w,r_w)$ is homeomorphic to the balanced product $\Ga\times_{\Ga_w}B(w,r_w)$ and $\bigcup_{w\in W}B(w,r_w)=X$;
\item[(3)] there exists an equivariant partition of unity associated to $\{\phi_w\}_{w\in W}$ such that $\supp(\phi_w)\subseteq B(w,r_w)$ where $r_w$ is defined as in (2) and $\phi_{\gamma w}=\gamma\phi_w$.
\end{enumerate}\end{Lem}

\begin{proof}[Sketch of Proof]
Denote $\pi: X\to X/\Ga$ the quotient map. By using \cite[Lemma A.2.7]{HIT2020}, for each $x\in X/\Ga$, there exists $r_x>0$ such that $\pi^{-1}(B(x,r_x))$ has the form of the balanced product by a finite subgroup. Thus the existence of $W$ is guaranteed by Zorn's Lemma and the fact that $X$ is locally compact and secondly countable.
\end{proof}

The following proposition can be found in \cite{HIT2020}, we provide a sketch of proof here for the convenience of readers.

\begin{Pro}\label{classifying proper}
Let $X$ be a metric space with bounded geometry and proper isometric action of a countable discrete group $\Ga$. Then
for each $Z\in \P(X,\Ga)$, there exists $d\geq0$ and a continuous equivariant coarse equivalence $f:Z\to P_d(X)$. Moreover, $f$ is unique up to homotopy equivalence.
\end{Pro}

\begin{proof}
By using Lemma \ref{partition of unity}, choose a countable discrete $\Ga$-invariant subset $W\subseteq Z$ such that $W$ is a $r$-net in $Z$ and a $\Ga$-partition of unity $\{\phi_w\}_{w\in W}$.

As $p_Z: Z\to X$ is an equivariant coarse equivalence, the restriction $p_Z: W\to Z$ is also a coarse equivalence. We define $f_0:W\to P_d(X)$ by $f_0(w)=p_Z(w)$ for each $w\in W$ for sufficiently large $d$, which is a coarse equivalence. Define $f:Z\to P_d(X)$ by
$$f(x)=\sum_{w\in W}\phi_w(x)p_X(w).$$
To see $f$ is well-defined, $\phi_w(x)\ne 0$ means that $d(x,w)\leq r$. Thus
$$\supp(f(x))=\{p_X(w)\mid d(x,w)\leq r\}.$$
As $p_Z$ is a coarse equivalence, there must exists $d>0$ associated to $r$ such that $\diam(\supp(f(x)))\leq d$. It is clear that $f$ is continuous and $\Ga$-equivariant. $f$ is also a coarse equivalence as $W$ is a net in $Z$ and $f_0:W\to P_d(X)$ is a coarse equivalence. Finally, the uniqueness of $f$ follows from a linear homotopy as in \cite[Lemma 7.2.14]{HIT2020}.
\end{proof}

As a result, one has the following model for $KX^{\Ga}_*(X)$:

\begin{Thm}[\cite{HIT2020}]
Consider the direct system:
$$P_0(X)\to P_1(X)\to \cdots,$$
where the map $P_d(X)\to P_{d+1}(X)$ is given by the canonical inclusion. It gives us a sequence
$$\begin{tikzcd}
K^{\Ga}_*(P_0(X)) \arrow[r, ] \arrow[d,  "\ind"] & K^{\Ga}_*(P_1(X)) \arrow[d, "\ind"] \arrow[r, ]& \cdots \\
K_*(C^*(X)^{\Ga}) \arrow[r, Rightarrow, no head]              & K_*(C^*(X)^{\Ga})\arrow[r, Rightarrow, no head]   & \cdots,
\end{tikzcd}$$
Then the inductive limit group $\lim_{d\to\infty}K^{\Ga}_*(P_d(X))$ is a universal $K^{\Ga}_*$-group for $\P(X,\Ga)$, i.e.,
$$KX^{\Ga}_*(X)\cong\lim_{d\to\infty}K^{\Ga}_*(P_d(X)),$$
and the assembly map is identified with the induced index map from the above sequence
$$\mu_X^{\Ga}:\lim_{d\to\infty}K^{\Ga}_*(P_d(X))\to K_*(C^*(X)^{\Ga}).$$ 
\end{Thm}

Combining Yu's localization technique, we have the following corollary:

\begin{Cor}
The assembly map is identified with the induced evaluation map:
$$ev_*:\lim_{d\to\infty}K_*(C^*_L(P_d(X))^{\Ga})\to K_*(C^*(X)^{\Ga}),$$
where $C^*_L(P_d(X))^{\Ga}$ is the equivariant localization algebra as we defined in Definition \ref{equivariant localization algebra}.
\end{Cor}

\subsection{The Milnor-Rips complexes for free and proper actions}\label{universal property for MR complexes}

In this subsection, we will give a concrete model of classifying spaces for \emph{free and proper} actions $\Ga\curvearrowright X$ by using the Milnor-Rips complexes introduced by G.~Yu in \cite{Yu1995}. The construction is modeled on Milnor's join construction \cite{Milnor56} but slightly different from Milnor's original version. Without loss of generality, we shall moreover assume $\Ga\car X$ is \emph{free} in this subsection (we shall explain detailed reasons in the next section).

Recall that a $\Ga$-space $Y$ is said to admit a uniformly locally finite $\F$-cover if there exists an $\F$-cover $\{U_i\}_{i\in I}$ such that
$$\#\{i\in I\mid y\in U_i\}\leq N$$
for all $y\in Y$ and for some $N>0$ (see Definition \ref{Def ulf cover}). Denote by $\F\P(X,\Ga)$ the subcategory of $\P(X,\Ga)$ which only contains objects $(Z,p_Z)$ where the $\Ga$-action on $Z$ is free and \emph{$Z$ admits a uniformly locallt finite $\F$-cover}. 

 Denote $X^{*\infty}=X*X*\cdots*\cdots$ the infinite join product of $X$, see \cite{Milnor56} or \cite[Definition 11.1]{GTM20}. Then any point in $X^{*\infty}$ can be written as an infinite sequence $\langle t,x\rangle=(t_0x_0,t_1x_1,\cdots,t_kx_k,\cdots)$, where $x_k\in X$, $t_k\in[0,1]$ and only finite $t_k\ne 0$, $\sum_{k}t_k=1$. We denote
$$\supp(\langle t,x\rangle)=\{x_k\mid t_k\ne 0\}.$$
To simplify the notation, we shall omit $x_k$ in the sequence $\langle t,x\rangle$ when $t_k=0$ and we denote $\langle 1,x_0\rangle=(1x_0,0,0,\cdots)$ for each $x_0\in X$.

\begin{Def}\label{Milnor Rips}
For any $d\geq 0$ and $n\in\IN$, the \emph{Milnor-Rips complex}, denoted by $\MRdnX$, is defined to be the set of equivalence classes of infinite sequences $\langle t,x\rangle=(t_0x_0,t_1x_1,\cdots,t_kx_k,\cdots)$, where
\begin{itemize}
\item[(1)]$t_k\geq0$ for each $k$ such that $t_k=0$ for all $i>n$ and $\sum_{k\in\IN}t_k=1$;
\item[(2)]$x_k\in X$ for each $k$ and $\diam(\supp(\langle t,x\rangle))\leq d$;
\item[(3)]two sequences $\langle t,x\rangle=(t_0x_0,t_1x_1,\cdots,t_kx_k,\cdots)$ and $\langle t',x'\rangle=(t'_0x'_0,t'_1x'_1,\cdots,t'_kx'_k,\cdots)$ are equivalent if
\par $(3_a)$ $t_k=t'_k$ for each $k$ and $x_k=x_k'$ for all $k$ such that $t_k\ne 0$, or;
\par $(3_b)$ there exists $k_0$ such that $t_{k_0}=0$, $t_k=t_k'$ and $x_k=x_k'$ for all $k<k_0$ and $t_{k+1}=t_k'$ and $x_{k+1}=x_k'$ for all $k\geq k_0$.
\end{itemize}\end{Def}

In this view, $\MRdnX$ has a natural structure of $CW$-complex since it is a quotient of a subspace of $Z^{*\infty}$. For a fixed sequence $x=(x_0,x_1,\cdots,x_n,\cdots)$ with $\diam(\{x_0,\cdots,x_n\})\leq d$, a \emph{semi-simplex} in $\MRdnX$ associated with $X$ is defined to be the set
$$\Delta_{x}=\{\langle s,w\rangle\in\MRdnX\mid w=x\}.$$
Let $\Delta_n$ be the standard $n$-simplex. Then there exists a canonical attaching map $\varphi_x:\Delta_n\to\Delta_x$ defined by
$$\varphi_x\left(\sum_{i=0}^nt_i[i]\right)=(t_0x_0,t_1x_1,\cdots,t_nx_n,0,\cdots).$$
This gives the CW-complex structure for $\MRdnX$. By definition, one can also see that the dimension of $\MRdnX$ is exactly equal to $n$.

For each $d_1\leq d_2$ and $n_1\leq n_2$, there exists a canonical inclusion $\widetilde{P_{d_1,n_1}}(X)\to \widetilde{P_{d_2,n_2}}(X)$, so we can identify $\widetilde{P_{d_1,n_1}}(X)$ with the image of the inclusion in $\widetilde{P_{d_2,n_2}}(X)$. When $d=0$ and $n=0$, one can see that $X$ is identified with $\widetilde{P_{0,0}}(X)$ by $x_0\mapsto\langle 1,x_0\rangle$, which is exactly the 0-dimensional cell of $\MRdnX$ for each $d\geq 0$ and $n\in\IN$. For any $x_0\in X$, any semi-simplex $\Delta_w$ attaching $x_0$ must satisfy that $x_0=w_i$ for some $i\leq n$ and $\diam(\{w_0,\cdots,w_n\})\leq d$. Thus, there are only finitely many semi-simplexes that attach $x_0$ since $X$ has bounded geometry, i.e., $\MRdnX$ is locally compact.

We shall then equip $\MRdnX$ with a proper metric. Let $\Delta_w$ be a semi-simplex, $\Delta_n$ the standard $n$-simplex equipped with the sphere metric $d_S$, $\varphi_w:\Delta_n\to \Delta_w$ the attaching map defined as above. For any two points $\langle s^{(0)},w\rangle,\langle s^{(1)},w\rangle\in\Delta_w$, a \emph{path} $\delta$ in $\Delta_w$ linking $\langle s^{(0)},w\rangle$ and $\langle s^{(1)},w\rangle$ is a finite sequence
$$\langle t^{(1)},w\rangle,\langle t^{(2)},w\rangle,\cdots,\langle t^{(k)},w\rangle\in\Delta_w.$$
and its length is defined to be
$$l(\delta)=\min\left\{d_S(x_s^{(0)},x_t^{(0)})+d_S(x_s^{(1)},x_t^{(k)})+\sum_{i=1}^{k-1}d_S(x_t^{(i)},x_t^{(i+1)})\right\}.$$
where $x_s^{i}\in\varphi^{-1}_w(\langle s^{(i)},w\rangle),x_t^{i}\in\varphi^{-1}_w(\langle t^{(i)},w\rangle)$ and the minimum is taken over all the choices of $x_s^{i}$ and $x_t^{i}$. The \emph{sphere metric} $d_{\Delta_w}$ on $\Delta_w$ is defined for any two points $\langle s^{(0)},w\rangle,\langle s^{(1)},w\rangle\in\Delta_w$ to be the infimum of the length of all path in $\Delta_w$ linking $\langle s^{(0)},w\rangle$ and $\langle s^{(1)},w\rangle$.

The metric on $\MRdnX$ is defined to be the largest metric satisfying that
$$d(\langle 1,x_1\rangle,\langle1,x_2\rangle)\leq d(x_1,x_2)\quad\text{and}\quad d(\langle t,x\rangle,\langle t',x\rangle)\leq d_{\Delta_x}(\langle t,x\rangle,\langle t',x\rangle)$$
where $x_1,x_2\in X$, $\langle t,x\rangle=(t_0x_0,\cdots,t_kx_k,\cdots)$, $\langle t' ,x\rangle=(t_0'x_0,\cdots,t_k'x_k,\cdots)\in\Delta_x\subseteq \MRdnX$. We define the $\Ga$ action on $\MRdnX$ by
$$\gamma(t_0x_0,t_1x_1,\cdots,t_kx_k,\cdots)=(t_0\gamma x_0,t_1\gamma x_1,\cdots,t_k\gamma x_k,\cdots)$$
As $\Ga\car x$ is free and proper, it is clear that $\Ga\car\MRdnX$ is also free and proper. One can also check by definition that this $\Ga$-action is isometric and the canonical inclusion $X\hookrightarrow\widetilde{P_{0,0}}(X)$ is a coarse equivalence (this metric restricts to the same metric with the metric on $X$ hereditary from $P_d(X)$ for sufficiently large $n$, one is referred to Section 7 of \cite{HIT2020} for more helpful discussions).

There is a canonical map $\pi:\MRdnX\to P_d(X)$ defined by
$$\pi(\langle t,x\rangle)=\sum_{k}t_kx_k.$$
It is clear that $\pi$ is continuous and $\Ga$-equivariant. As $\pi$ maps $\widetilde{P_{0,0}}(X)$ to $P_0(X)\subseteq P_d(X)$, thus $\pi$ is also a coarse equivalence. Then $\bigcup_{d\geq0,n\in\IN}\MRdnX$ is made into a classifying space for the category $\F\P(X,\Ga)$ in the following sense:

\begin{Thm}\label{classify free and proper}
Let $X$ be a bounded geometry metric space with a proper, free, and isometric action of a countable discrete group $\Ga$. Then\begin{enumerate}
\item[(1)]$\MRdnX\in\F\P(X,\Ga)$ for each $d\geq 0$ and $n\in\IN$;
\item[(2)]for any $Z\in\F\P(X,\Ga)$, there exist sufficiently large $d,n$, and an equivariant continuous map $g: Z\to\MRdnX$ such that $\pi\circ g: Z\to P_d(X)$ is a coarse equivalence which makes the following diagram commutes:
$$\begin{tikzcd}
X\arrow[d,  "g"'] \arrow[dr,  "f"]\\
\MRdnX\arrow[r, "\pi"]             & P_d(X)   
\end{tikzcd}$$
where $f$ is defined as in Proposition \ref{classifying proper}. To sum up, the map $f$ we defined in Proposition \ref{classifying proper} factors through $\pi:\MRdnX\to P_d(X)$.
\end{enumerate}\end{Thm}

\begin{proof}
For the first term, it suffices to show $\MRdnX$ admits a uniformly locally finite $\F$-cover for each $d\geq 0$ and $n\in\IN$. The proof is similar to Proposition \ref{ulf cover}. Since $dim(\MRdnX)=n$, we still set $\varepsilon=\frac{2}{(n+2)^2}>0$.

For any $x_0\in X$ and $k\in\{0,\cdots,n\}$, defined
$$U_{x_0,k}=\{(t_0z_0,\cdots,t_nz_n,0,\cdots)\mid z_k=x_0, t_k>t_j+\varepsilon\text{ where }j\ne k\}.$$
Set $\U^0=\{U_{x,k}\}_{x\in X,1\leq k\leq n}$. For a pair $(x_0,x_1)\in X^2$ with $d(x_0,x_1)\leq d$ and $k_0,k_1\in\{1,\cdots,n\}$ with $k_0< k_1$, define
$$U_{(x_0,x_1),k_0,k_1}=\{(t_0z_0,\cdots,t_nz_n,0,\cdots)\mid z_{k_i}=x_i, t_{k_i}>t_j+\varepsilon\text{ where }i\in\{0,1\},j\ne k_0,k_1\}.$$
Set $\U^1=\{U_{(x_0,x_1),k_0,k_1}\}_{(x_0,x_1)\in X^2,1\leq k_0<k_1\leq n}$. Similarly, for a sequence $(x_0,x_1,\cdots,x_m)\in X^{m+1}$ with $d(z_i,z_j)\leq d$ for any $0\leq i<j\leq m$ and a sequence of integers $0\leq k_0<\cdots<k_m\leq n$
\[\begin{split}&U_{(x_0,\cdots,x_m),k_0,\cdots,k_m}\\
=&\{(t_0z_0,\cdots,t_nz_n,0,\cdots)\mid z_{k_i}=x_i, t_{k_i}>t_j+\varepsilon\text{ where }i\in\{0,\cdots,m\},j\ne k_0,\cdots,k_m\}.\end{split}\]
Set $\U^m=\{U_{(x_0,\cdots,x_m),k_0,\cdots,k_m}\}_{(x_0,\cdots,x_m)\in Z^{m+1},1\leq k_0<\cdots<k_m\leq n}$. Since $\MRdnZ$ is finite-dimensional, we can construct $\U^0,\cdots,\U^n$ by repeating the procedure above. Set $\U=\U^0\cup\cdots\cup\U^n$. One can follow a proof similar to Proposition \ref{ulf cover} to show that $\U$ is a uniformly locally finite $\F$-cover since
\[\begin{split}&\overline{U_{(x_0,\cdots,x_m),k_0,\cdots,k_m}}\\
=&\{(t_0z_0,\cdots,t_nz_n,0,\cdots)\mid z_{k_i}=x_i, t_{k_i}\geq t_j+\varepsilon\text{ where }i\in\{0,\cdots,m\},j\ne k_0,\cdots,k_m\}.\end{split}\]
And each set $U_{(x_0,\cdots,x_m),k_0,\cdots,k_m}$ is homotopic to the point $\{(t_0z_0,\cdots,t_nz_n,0,\cdots)\}$ where $z_{k_i}=x_i$ and $t_{k_i}=\frac 1{m+1}$ for each $i$ by using a linear homotopy.

For the second term, by using Lemma \ref{partition of unity}, choose a countable discrete $\Ga$-invariant subset $W\subseteq Z$ such that $W$ is a $r$-net in $Z$ and a $\Ga$-partition of unity $\{\phi_w\}_{w\in W}$. Since $Z\in\F\P(X,\Ga)$, we moreover assume that $\{B(w,r_w)\}_{w\in W}$ is a uniformly locally finite $\F$-cover. As the action is free, we can choose a fundamental domain $W_0$ of $W$ such that $W=\bigsqcup_{\gamma\in\Ga}\gamma W_0$. Without loss of generality, write $W_0=\{w_1,w_2,\cdots\}$. Then $\{\phi_{\gamma w_i}\}_{w_i\in W_0,\gamma\in\Ga}$ is identified with $\{\phi_w\}_{w\in W}$ and $\supp(\phi_{\gamma_1w_i})\cap\supp(\phi_{\gamma_2w_i})=\emptyset$ if $\gamma_1\ne\gamma_2$.

Since $(Z,p_Z)\in\F\P(X,\Ga)$, there exists $n\in\IN$ such that, there are at most $n$ open set in $\{B(w,r)\}_{w\in W}$ that covers $z$ for any $z\in Z$ and there is at most one element in $\{\gamma\cdot B(w,r)\}_{\gamma\in\Ga}$ that covers $z$. Define $g_0:W\to\MRdnX$ by
$$g_0(w)=\langle 1,p_X(w)\rangle$$
We then define $g:Z\to\MRdnX$ by
$$g(z)=\left(\phi_{\gamma_1w_1}(z)p_X(\gamma_1w_1),\cdots,\phi_{\gamma_kw_k}(z)p_X(\gamma_kw_k),\cdots\right)$$
where $\gamma_k$ is the unique element such that $x\in\gamma_kB(w_k,r)$.

It follows from a similar argument in Proposition \ref{classifying proper} that $g$ is well-defined for sufficiently large $d$ and an equivariant continuous map. One can check that
$$\pi\circ g(x)=\sum_{w_i\in W_0}\phi_{\gamma_iw_i}(x)p_X(\gamma_iw_i)=\sum_{w\in W}\phi_w(x)p_X(w).$$
Combining the definition of $f: Z\to P_d(X)$ in Proposition \ref{classifying proper}, this shows that the diagram commutes.
\end{proof}

\begin{con}[Analytic equivariant coarse Novikov conjecture]
Let $\Ga$ be a countable discrete group, $X$ a proper metric space with a free and proper $\Ga$-action. Then the composition of $\pi_*$ and the assembly map is $\IQ$-injective, i.e., the map:
$$\nu_{X}^{\Ga}:\lim_{d,n\to\infty}K^{\Ga}_*\left(\MRdnX\right)\ox\IQ\xrightarrow{\pi_*}\lim_{d\to\infty}K_*^{\Ga}(P_d(X))\ox\IQ\to K_*(C^*(X)^{\Ga})\ox\IQ$$
is injective. We shall call the map $\nu_X^{\Ga}$ the Miščenko-Kasparov assembly map.
\end{con}

\begin{Rem}
\label{rem:Milnor-Ripe-classifying}
When $X$ is $\Ga$-compact, one can always choose $W$ in Lemma \ref{partition of unity} to be of bounded geometry. In this case, the category $\F\P(X,\Ga)$ is exactly all cocompact principal $\Ga$-bundles. Any cocompact principal $\Ga$-boundle $\pi:E\to B$ admits a uniformly locally finite $\F$-cover by using the local trivialization condition and the compactness of $B$. Then the classifying space can be realized by using the universal principal bundle $\pi:\E\Ga\to B\Ga$ (see \cite{GTM20} for a detailed introduction to the principal bundles). Then the conjecture above is exactly the analytic Novikov conjecture stated in \cite{Roe1996}, i.e., the map above is identified with the group homomorphism below:
$$RK_*(B\Ga)\ox\IQ\to RK^{\Ga}_*(\EG)\ox\IQ\to K_*(C^*_r\Ga)\ox\IQ.$$
\end{Rem}

If $\Ga$ is torsion-free, then $\Ga\car P_d(X)$ is free and proper for any $d\geq 0$. Moreover, one also has that  $P_d(X)\in \F\P(X,\Ga)$ since $P_d(X)$ is equivariantly coarsely equivalent to $X$ and $P_d(X)$ admits a uniformly locally finite $\F$-cover as we proved in Proposition \ref{ulf cover}. Thus, $\lim\limits_{d\to\infty}K_*^{\Ga}(P_d(X))$ also plays a role of a universal $K^{\Ga}$-group for $\F\P(X,\Ga)$. One has that
$$\pi_*:\lim_{d,n\to\infty}K^{\Ga}_*\left(\MRdnX\right)\to\lim_{d\to\infty}K_*^{\Ga}(P_d(X))$$
is an isomorphism from the universal property. In particular, when $\Ga$ is a trivial group, omit $\Ga$ from the notation and 
$$\lim_{d,n\to\infty}K_*\left(\MRdnX\right)\xrightarrow{\pi_*}\lim_{d\to\infty}K_*(P_d(X))\xrightarrow{\mu_X}K_*(C^*(X))$$
still gives the coarse Baum-Connes assembly map. This generalizes \cite[Lemma 2.2]{Yu1995} from cocompact action to the general case.

\section{A concrete model for spaces with property TAF}\label{sec: A Concrete model for space with Property TAF}

In this section, we discuss compact Hausdorff spaces with property TAF, as defined in \cite[Definition~5.2]{AAS2020}. 
With their help, it is easy to produce free and proper actions (that is, elements in $\F\P(X,\Ga)$) out of proper ones (that is, elements in $\P(X,\Ga)$), which ultimately will help us prove map~(12) in diagram~\eqref{main diam} is rationally injective. 
While a construction of spaces with property TAF was given in \cite[Definition~5.6]{AAS2020}, we shall present an alternative concrete model for spaces with property TAF which is tailored to our construction of the Milnor-Rips complexes, simplifying proofs in \Cref{sec: Proof of the main theorem}.

Recall that $\Ga\car X$ has \emph{equivariant bounded geometry} means that for any $R>0$, there exists $N>0$ such that 
$$\#\{\gamma\in\Ga\mid B(\gamma\cdot x,R)\cap B(x,R)\ne\emptyset\}<N$$
for all $x\in X$ (see Definition \ref{ebg}).

The following definition comes from \cite{AAS2020}:

\begin{Def}[\cite{AAS2020}]
A compact Hausdorff space $Z$ with an action of $\Ga$ is said to have property TAF (torsion acts freely) if every torsion element acts freely on $Z$.
\end{Def}

It is proved in \cite[Theorem 5.7]{AAS2020} that for any discrete group $\Gamma$,  there exists a compact $\Gamma$-space with property TAF. For the convenience of the reader, we shall briefly recall the proof here.
For any $F\subseteq\Ga$ be a finite subgroup, define
$$Z_F=\prod_{[g]\in\Ga/F}g F=\{s:\Ga/F\to\Ga\mid \pi_F\circ s=id\}$$
where $\pi_F:\Ga\to\Ga/F$ is the quotient map.
Equip $Z_F$ with the Tychonoff topology, which makes $Z_F$ a compact Hausdorff space. 
Define the group action by
$$(\gamma\cdot s)([g])=\gamma\cdot s([\gamma^{-1}g])$$
for all $\gamma\in \Ga$ and $s\in Z_F$. Notice that $F$ acts freely on $Z_F$. Indeed, if $f\in F$ satisfies that $(f\cdot s)([e])=s([e])$, then
$$f\cdot s([f^{-1}])=f\cdot s([e])=s([e]),$$
i.e., $f=e$. Then
$$Z_\Ga=\prod_{F\leq\Ga,\text{$F$ is finite}}Z_F.$$
equipped with the Tychonoff topology has property TAF. Then it is clear that $Z_\Ga$ is a compact Hausdorff space which has property TAF.

Considering the product space $X\times Z_\Ga$, we define the $\Ga$-action on $X\times Z_\Ga$ by the diagonal action, i.e., $\gamma\cdot(x,a)=(\gamma\cdot x,\gamma\cdot a)$ for all $x\in X$ and $a\in Z_{\Ga}$. This action is proper and free. Define $p:X\times Z_{\Ga}\to X$ to be the projection onto $X$. Fix $a_0\in Z_\Ga$. By definition, we can write as the disjoint union of all the orbits $X=\bigsqcup_{x\in D}\Ga\cdot x$, where $D$ is the maximal set of $X$ such that $\Ga\cdot x\cap\Ga\cdot y=\emptyset$ for any $x\ne y\in D$. Then let $M_X=\bigsqcup_{z\in D}\Ga\cdot(z,a_0)$ be a $\Ga$-invariant subspace of $X\times X_{\Ga}$. Define the metric on $M_X$ by
$$d((x,a),(x',a'))=\left\{\begin{aligned}&d(x,x')+1,&&\text{if }a\ne a';\\&d(x,x'),&&\text{if }a=a'.\end{aligned}\right.$$
It is clear that this metric is well-defined and $\Ga\car M_X$ is an isometric action under this metric. The canonical projection $p: M_X\to X$ is an equivariant coarse equivalence.

As $X$ has bounded geometry, there exists $N_0>0$ such that $\#B(z,R)\leq N_0$. Since $\Ga\car x$ has equivariant bounded geometry, there exists $N>0$ such that $\#\Ga_{x}\leq N_1$ for some $N_1>0$, where $\Ga_x$ means the stabilizer of $x$. Then we have that $\#(B((x,a),R)\cap M_X)\leq N_0\cdot N_1$. This shows that $M_X$ has bounded geometry, we conclude that $M_X\in\F\P(X,\Ga)$. Moreover, $\Ga\car M_X$ has equivariant bounded geometry since the action $\Ga\car M_X$ is free (combining Remark \ref{Rem ebg}). Then we define the \emph{Miščenko-Kasparov assembly map for $\Ga\car X$} to be:
$$\nu^{\Ga}_X:\lim_{d,n\to\infty}K_*^{\Ga}(\widetilde{P_{d,n}}(M_X))\ox\IQ\xrightarrow{\pi_*}\lim_{d\to\infty}K_*^{\Ga}(P_d(M_X))\cong\lim_{d\to\infty}K_*^{\Ga}(P_d(X))\ox\IQ\to K_*(C^*(X)^{\Ga})\ox\IQ,$$
the middle isomorphism follows from the fact that $M_X$ is equivariant coarse equivalent to $X$ and the universal property of $KX_*^{\Ga}(X)$. Without loss of generality, we shall assume $\Ga\car X$ is free in the rest of this paper (otherwise just consider $\Ga\car M_X$ instead of $\Ga\car X$).

In the rest of this section, we shall introduce a new model of $\Ga$-space with property TAF by using the action $\Ga\car X$, this model will be easier for us to deal with in the following sections.

\begin{Def}\label{TAF space}
    Given a set $S$, we let $\Omega_S$ be the set of all \emph{linear orders} (or \emph{irreflexive total orders}) on $S$. Recall that an order $\R$ on $S$ is linear if either $x_1<_\R x_2$ or $x_2<_\R x_1$ for any $x_1\ne x_2$. For any subset $S' \subseteq S$, there exists a canonical projection $p_{S,S'}:\Omega_{S}\to\Omega_{S'}$ defined by $\R\mapsto \R|_{S'}$, restricting $\R$ onto $S'$ for any $\R\in\Omega_{S}$. 

    Let $X$ be a free and proper $\Ga$-space with bounded geometry. We write $F \Subset X$ to indicate $F$ is a finite subset of $X$. 
    Note that $\{ \Omega_F \}_{F \Subset X}$ together with the restriction maps $\{ p_{F, F'} \}_{F' \subseteq F \Subset X}$ forms an inverse system of finite sets. 
    Observe that as a set, $\Omega_X$ is canonically identified with the inverse limit of this inverse system; more precisely, $\Omega_X$ is identified, using the maps $p_{X, F}$ for $F \Subset X$, with the subset of $\prod_{F \Subset X} \Omega_F$ consisting of tuples $(\R_F)_{F \Subset X}$ satisfying $(\R_F)|_{F'} = \R_F$ whenever $F' \subseteq F \Subset X$. 

    Now we equip each $\Omega_F$, for $F \Subset X$, with the discrete topology, and topologize $\Omega_X$ as the inverse limit of this inverse system of finite discrete spaces. 
    Applying the Tychonoff theorem or the Alexander subbase theorem, we see that $\Omega_X$ is a compact Hausdorff space. 
\end{Def}

Alternatively, one can take $\Omega_X$ to be the spectrum of the direct limit $C^*$-algebra of the directed system consisting of $\{ C(\Omega_F) \}_{F \Subset X}$ together with the induced $*$-homomorphisms $\{ p^*_{F, F'} \colon  C(\Omega_{F'}) \to C(\Omega_F) \}_{F' \subseteq F \Subset X}$. 

Since $X$ is countable, we may choose an increasing sequence of finite subsets $\{F_n\}_{n\in\IN}\subseteq\F$ such that $\bigcup F_n=X$. Since this forms a cofinal sequence in the directed set of all finite subsets of $X$, we see that $C(\Omega_X) = \lim\limits_{n\to\infty}C(\Omega_{F_n})$ and thus $C(\Omega_X)$ is separable.

The $\Ga$-action on $\Omega_X$ is induced by the $\Ga$-action on $X$. For any $\R\in\Omega_X$ and $\gamma\in\Ga$, we define $\gamma\R$ by
$$x_1<_{\gamma \R} x_2\iff \gamma^{-1}x_1<_{\R}\gamma^{-1}x_2.$$
For any torsion element $\gamma\in\Ga$, there exists $n\in\IN$ such that $\gamma^n=1$. Since $\Ga\car X$ is free, for any $x\in X$, the points $x,\gamma x,\cdots,\gamma^n x$ are distinct from each other. Write $x_m=\gamma^mx$ for simplicity. For any $\R\in\Omega_X$, there must be a ``minimal" element under $\R$, say $x_m$, such that $x_m\leq_\R x_k$ for any $k$. Then $\gamma\R$ can never be equal to $\R$ since $x_m$ can never be the ``minimal" element under $\gamma\R$. This shows that $\Omega_X$ has property TAF.

We now discuss the metric structure of $\Omega_X$. Let $\{F_n\}$ be an increasing sequence of finite subsets $\{F_n\}_{n\in\IN}\subseteq\F$ such that $\bigcup F_n=X$. For each $n\in\IN$, the metric $d_n$ on each $\Omega_{F_n}$ is defined to be
$$d_n(\R,\R')=\left\{\begin{aligned}&2^{-n},&&\R\ne\R';\\&0,&&\R=\R',\end{aligned}\right.$$
for any $\R,\R'\in \Omega_{F_n}$. This metric clearly induces the topology of $\Omega_{F_n}$. Define the metric $d$ on $\Omega_X$ to be
$$d(\R,\R')=\sum_{n\in\IN}d_n(\R|_{F_n},\R'|_{F_n}),\quad\text{ for any }\R,\R'\in\Omega_X.$$
We leave it to the reader to check that $d$ is a well-defined metric. The topological base induced by this metric is given by $\{p^{-1}_{F_n}(\R|_{F_n})\mid n\in\IN,\R\in\Omega_X\}$. Since $\bigcup_{n\in\IN}F_n=X$, the topology induced by this metric coincides with that we defined above. 

A metrization on a $\Ga$-space $X$ is said to be \emph{invariant} if the $\Ga$-action under this metric is isometric. One should notice that the $\Ga$-action on $\Omega_X$ is not isometric under the metric we defined above. It also seems unlikely that an invariant metrization of $\Omega_X$ can be found in general.  However, for any separable metric space $X$ with a proper $\Ga$-space, it is reasonable to expect a metrization on $X\times\Omega_X$ such that the $\Ga$-action is isometric since $X\times\Omega_X$ is separable, metrizable and equipped with a proper $\Ga$-action (cf. \cite[Theorem 4.3.4]{metrizationGaspace}). In the following, we shall construct a concrete invariant metrization on $X\times\Omega_X$.

Let $X$ be a proper $\Ga$-space, one can find a $\Ga$-invariant $r$-net $W$ in $X$ by using Lemma \ref{partition of unity} such that $W$ is locally finite. Set $\{\phi_w\}_{w\in W}$ to be an equivariant partition of unity associated with $\{B(w,r_w)\}_{w\in W}$, where $r_w$ is chosen as in Lemma \ref{partition of unity}. For a fixed point $w\in W$, denote by $\Ga\cdot w=\{\gamma w\mid \gamma\in\Ga\}$ the orbit of $w$ which is a subset of $W$. Since $\Ga\car X$ is proper, the stabilizer $\Ga_w$ is a finite subgroup of $\Ga$. Moreover, it is clear that $\gamma\Ga_w\gamma^{-1}$ is the stabilizer subgroup for $\gamma w$. We define $d_w$ to be the metric on $\Omega_X$ associated with $\{F_n\}_{n\in\IN}$ as above such that each $F_n$ is $\Ga_w$-invariant. For any $\gamma w\in\Ga\cdot w$, we define $d_{\gamma w}$ to be the metric on $\Omega_X$ associated with $\{\gamma F_n\}_{n\in\IN}$ as above. It is not hard to check that $d_w(\R,\R')=d_{\gamma w}(\gamma\R,\gamma\R')$ for any $\R,\R'\in\Omega_X$ and $\gamma\in\Ga$ by definition. For any $x\in X$, we define the metric $d_x$ on $\Omega_{X}$ by
$$d_x(\R,\R')=\sum_{w\in W}\phi_w(z)\cdot d_w(\R,\R').$$
It is clear that $d_x$ is a metric on $\Omega_X$. Then the family of metrics $\{d_x\}_{x\in X}$ on $\Omega_X$ satisfies that
$$d_x(\R,\R')=d_{\gamma x}(\gamma\R,\gamma\R')$$
for any $\R,\R'\in\Omega_X$ and $\gamma\in\Ga$ by definition.

Define the \emph{twisted metric} $\wt d$ on $X\times\Omega_X$ to be the largest metric such that
$$\wt d((x_1,\R),(x_2,\R))=d_X(x_1,x_2)\quad\text{ and }\quad\wt d((x,\R_1),(x,\R_2))\leq d_x(\R_1,\R_2).$$
Such a metric always exists by \cite[Lemma 1.4]{Roe2005}. The idea for this metric is inspired by \cite[Definition 3.6]{SW2021}. The twisted metric we defined can also be realized by the infimum of \emph{mileages} of sequences, see \cite[Proposition 1.6]{Roe2005} and \cite[Proposition 3.7]{SW2021}. It is clear that the $\Ga$-action under this metric is isometric, and the canonical projection from $X\times\Omega_X\to X$ is an equivariant coarse equivalence.

\begin{Rem}\label{metric on product space}
By using \cite[Theorem 4.3.4]{metrizationGaspace}, one can also find an invariant metrization of $X\times\Omega_X$ by embedding $X\times\Omega_X$ equivariantly into a Hilbert space. We should mention that the canonical projection $X\times\Omega_X\to X$ may not be a coarse equivalence under this metric. However, in the following sections, we shall only consider the localization algebra for $X\times\Omega_X$. Since the $K$-theorey of localization algebra of $X\times\Omega_X$ only depends on its topological structure (see Theorem \ref{non-coe comm-prop} or \cite[Theorem 6.6.2]{HIT2020}), it actually does not matter if one uses this invariant metric instead of the metric we defined above.
\end{Rem}

In the last part, we shall construct an \emph{invariant trace} $\tau$ on $C(\Omega_X)$, that is a trace satistying that $\tau(\gamma f)=\tau(f)$ for any $f\in C(\Omega_X)$ and $\gamma\in\Ga$. Such a trace induces a canonical trace on $C(\Omega_X)\rtimes_r\Ga$ defined by
$$\tau\left(\sum_{\gamma\in\Ga}f_\gamma\lambda_\gamma\right)=\tau(f_e),$$
we here abuse the notation and write $\tau$ for both trace on $C(\Omega_X)$ and $C(\Omega_X)\rtimes_r\Ga$, and we write elements in $C_c(\Gamma,C(\Omega_X))\subseteq C(\Omega_X)\rtimes_r\Ga$ as a finite sum $\sum f_\gamma\lambda_\gamma$. Fix an increasing sequence of finite subsets $\{F_n\}_{n\in\IN}$ such that $\bigcup_{n\in\IN} F_n=X$, we define the measure $\mu_n$ on $\Omega_{F_n}$ to be the counting measure. For any positive function $f\in C(\Omega_X)$, we define
$$\tau(f)=\sup\left\{\int_{\Omega_{F_n}}\phi(x)d\mu_n(x)\ \Big|\ 0\leq\phi\in C(\Omega_{F_n})\text{ and }p^*_{F_n}(\phi)\leq f  \right\},$$
where $p^*_{F_n}:C(\Omega_{F_n})\to C(\Omega_X)$ is induced by the canonical projection $\Omega_X\to\Omega_{F_n}$. It is clear that $\tau$ extends to a tracial state. We still need to show that $\tau$ is invariant. Fix $f\in C(\Omega_X)$, for any $\varepsilon>0$, there exists $n\in\IN$ and $\phi\in C(\Omega_X)$ such that $\|p^*_{F_n}(\phi)-f\|\leq\varepsilon$. For any $\gamma\in\Ga$, there must exists $N>n$ such that $F_n\cup\gamma F_n\subseteq F_N$. Then there is a canonical projection $p:\Omega_{F_N}\to\Omega_{\gamma F_n}$ by restriction. For any $\R\in\Omega_{F_n}$, we define $\gamma \R\in\Omega_{\gamma F_n}$ by
$$\gamma x_1<_{\gamma \R}\gamma x_2\iff x_1<_{\R}x_2.$$
Notice that $\R\mapsto\gamma\R$ is a bijection from $\Omega_{F_n}$ to $\Omega_{\gamma F_n}$, and we define $\gamma\phi\in C(\Omega_{\gamma F_n})$ by
$$(\gamma\phi)(\gamma\R)=f(\R).$$
It is direct to see that $\|\gamma f-p^*_{\gamma F_n}(\gamma\phi)\|\leq\varepsilon$.
Denote by $p^*_N(\phi)$ and $p^*_N(\gamma\phi)$ the image of $\phi$ and $\gamma\phi$ under the homomorphism induced by canonical projections from $\Omega_{F_N}$ onto $\Omega_{F_n}$ and $\Omega_{\gamma F_n}$, respectively. By definition, one has that
$$\int_{F_N}p^*_N(\phi)d\mu_N=\int_{F_N}p^*_N(\gamma\phi)d\mu_N.$$
We then conclude that $|\tau(f-\gamma f)|\leq 2\varepsilon$ for any $\varepsilon>0$, this means that $\tau(f)=\tau(\gamma f)$.

\section{$\A(M)$ and the bootstrap class}\label{Sec: Approximating AofM}

In this section, we take a digression to show the $C^*$-algebra $\A(M)$ associated to a separable Hilbert-Hadamard space $M$ is in the bootstrap class (the class $\N$ in \cite{Duke1987}) and thus satisfies the universal coefficient theorem (UCT) in $KK$-theory. 
Establishing this fundamental property of $\A(M)$ will enable us to complete our $K$-theoretic computations in \Cref{subsec:Kunneth}. 

Our approach here is built on the following immediate consequence of Definition~\ref{Bott}: one can write $\A(M)$ as a direct limit of $C^*$-subalgebras $\A(M,F)$, where $F$ ranges over finite subsets of (a dense subset of) $M$ (see also \cite[Lemma~7.2]{GWY2018}). 
The main objective of this section is thus to study the structure of these $C^*$-subalgebras $\A(M, F)$. 
It turns out that each $\A(M,F)$ is a Type I $C^*$-algebra when $F$ is a finite subset, a fact we prove by constructing, for each $n \in \mathbb{N}$, a universal $C^*$-algebra $\B_0(n)$ that is subhomogeneous and maps onto any $\A(M, F)$ with $|F| = n$.


\begin{Def}
For any positive integer $n$, let $M_n(\IR)_+$ be the space of all positive semidefinite real matrices, equipped with the usual topology inherited from $M_n(\IR)$. We define the $*$-algebra
$$\B(n)=\left\{f:M_n(\IR)_+\to\Cl(\IR^n)\,\Big|\, f\text{ is continuous and }f(T)\in\Cl(\operatorname{ran}(T^{1/2}))\right\}$$
with the algebraic operations defined pointwise, where $T^{1/2}\in M_n(\IR)_+$ is the square root matrix of $T$ and $\operatorname{ran}(T^{1/2})$ is the range\footnote{Note that $T^{1/2}$ in fact has the same range as $T$. } of $T^{1/2}$, as a subspace in $\IR^n$. We further define the $C^*$-algebras
$$\B_b(n)=\B(n)\cap C_b(M_n(\IR)_+,\Cl(\IR^n));$$
$$\B_0(n)=\B(n)\cap C_0(M_n(\IR)_+,\Cl(\IR^n)).$$
\end{Def}

Since $C_0(M_n(\IR)_+,\Cl(\IR^n))$ is a Type I $C^*$-algebra \cite{Blackoperator}, as a result, $\B_0(n)$ is Type I. We would like to define $*$-homomorphisms from $\B_0(n)$ to $\A(M)$. To do this, we first state an elementary fact in linear algebra.

\begin{Rem}\label{positive matrix}
Let $v_1, \ldots, v_n$ be vectors in a real Hilbert space $\H$. Let $A_{v_1, \ldots, v_n}$ be the covariance matrix, that is,  
\[A_{v_1, \ldots, v_n} = \left( \left\langle v_i , v_j \right\rangle \right)_{1 \leq i,j \leq n} \; ,\]
which is a positive semidefinite real matrix. Then there is a unique isometry 
\[\rho_{v_1, \ldots, v_n} \colon \operatorname{ran}\left({A_{v_1, \ldots, v_n}}^{1/2}\right) \to \H\]
such that for each $e_i$ in the standard basis of $\IR^n$ (considered as column vectors), we have 
\[\rho_{v_1, \ldots, v_n} \left( {A_{v_1, \ldots, v_n}}^{1/2} e_i \right) = v_i \; .\]
Indeed, the vectors $v_1, \ldots, v_n$ determine a linear map $T \colon \mathbb{R}^n \to \H$ such that $v_i = T e_i$ for $i = 1, \ldots, n$. Note that $A_{v_1, \ldots, v_n} = T^* T$ and $\rho_{v_1, \ldots, v_n}$  can be obtained by a polar decomposition $T = \rho_{v_1, \ldots, v_n} (T^* T)^{1/2}$. 
\end{Rem}

Now, let $M$ be a Hilbert-Hadamard space and let $x_1,\cdots,x_n$ be points in $M$. 
Using the notations in Definition~\ref{def:Clifford_generator} and Remark \ref{positive matrix}, we write, for any $x\in M$ and $t\in\IR_+=[0,\infty)$, 
\[A_{x,t; x_1, \cdots, x_n} = A_{C_{x_1}(x,t), \cdots, C_{x_n}(x,t)} = \left( \left\langle \log_x \left( x_i \right), \log_x \left( x_j \right) \right\rangle + t^2 \right)_{1 \leq i,j \leq n} \in M_n(\IR)_+\]
and 
\[\rho_{x,t; x_1, \cdots, x_n} = \rho_{C_{x_1}(x,t), \cdots, C_{x_n}(x,t)} : \operatorname{ran}\left({A_{x,t; x_1, \cdots, x_n}}^{1/2}\right) \to \H_x M \oplus \IR \; .\]
Note that the range of $\rho_{x,t; x_1, \cdots, x_n}$ lies in the $t$-level set in $\H_xM\oplus\IR$. Let 
\[\Cl\left(\rho_{x,t; x_1, \cdots, x_n}\right): \Cl \left(\operatorname{ran}\left({A_{x,t; x_1, \cdots, x_n}}^{1/2}\right)\right) \to \Cl\left(\H_x M \oplus \IR\right)\]
be the induced $*$-homomorphism between the Clifford algebras. Using these as well as the notations from Definition~\ref{def:Pi-alg}, we can define a $*$-homomorphism 
\[\beta_{x_1, \cdots, x_n} : \B(n) \to \Pi(M)\]
such that for any $f \in \B(n)$, $x \in M$ and $t \in [0,\infty)$, we have
\[\beta_{x_1, \cdots, x_n} (f) (x,t) = \left(\Cl\left(\rho_{x,t; x_1, \cdots, x_n}\right) \right) \left(f\left(A_{x,t; x_1, \cdots, x_n}\right)\right)\]
This restricts to a $*$-homomorphism from $\B_b(n)$ to $\Pi_b(M)$. 

We shall see that $\beta_{x_1, \cdots, x_n}$ restricts to a $*$-homomorphism from $\B_0(n)$ to $\A(M,F)$ where $F=\{x_1, \cdots, x_n\}$. To this end, we define the following unbounded multipliers analogous to the Clifford generators we used in the definition of $\A(M)$. 

\begin{Def}
For any positive integer $n$ and any column vector $v \in \IR^n$, we define an element $C_{n,v}$ in $\B(n)$ such that 
\[C_{n,v}(T) = T^{1/2} v\]
for any $T \in M_n(\IR)_{+}$. Moreover, for any $f \in C_0(\IR)$, we define $f \left(C_{n,v}\right) \in \B(n)$ such that at any  $T \in M_n(\IR)_{+}$, its value is derived by applying functional calculus, with $f$, to the self-adjoint element $T^{1/2} v$ in $\Cl\left(\operatorname{ran}\left({T}^{1/2}\right)\right)$.
\end{Def}

The following basic properties of $C_{n,v}$ are immediate. Recall from Definition~\ref{def:Clifford_generator} that for each $x\in M$, we denote by $C_{x}$ the corresponding Clifford generator for $\A(M)$. 

\begin{Lem}\label{connection AM and Bn}
The element $C_{n,v}$ defined above satisfies the following:
\begin{enumerate}[(1)]
\item For any $x_1, \cdots, x_n \in M$ and any $i \in \{1, \ldots, n\}$, we have $\beta_{x_1, \ldots, x_n} (C_{n, e_i}) = C_{x_i}$. 
\item The assignment $v \mapsto C_{n,v}$ is linear, i.e., for any $v,w \in \IR^n$ and $a \in \IR$, we have $C_{n, v+w} = C_{n, v} + C_{n, w}$ and $C_{n, av} = aC_{n,v}$. 
\end{enumerate}
\end{Lem}

\begin{proof}
By definition, for any $(x,t) \in M \times \IR_+$, we have 
\[\begin{split}\beta_{x_1, \cdots, x_n} (C_{n,e_i}) (x,t)&= \left(\Cl\left(\rho_{x,t; x_1, \cdots, x_n}\right) \right) \left(C_{n,e_i}\left(A_{x,t; x_1, \cdots, x_n}\right)\right)\\
&=A^{1/2}_{x,t; x_1, \cdots, x_n}e_i=(-\log_x(x_i),t)=C_{x_i}(x,t).\end{split}\]
Linearity follows directly from the definition.
\end{proof}

\begin{Lem}\label{generators Bn}
The $C^*$-algebra $\B_0(n)$ is generated by the subset
\[\left\{ \prod_{i=1}^{n} f_i\left(C_{n, e_i}\right) \mid f_i \in C_0(\IR) \right\} \]
inside $\B(n)$. 
\end{Lem}

\begin{proof}
First of all, we shall show that $\prod_{i=1}^{n}f(C_{n,e_i})$ belongs to $\B_0(n)$ for any $f_1,\cdots,f_n\in C_0(\IR)$. It suffices to show that for any $\varepsilon>0$, there exists a compact subset $K\subseteq M_n(\IR)_+$ such that $\|\prod_{i=1}^{n}f(C_{n,e_i})(T)\|\leq\varepsilon$ for any $T\in M_n(\IR)_+ \setminus K$. To this end, we find $R>0$ such that
$$|f_i(t)|<\frac{\varepsilon}{1+\prod_{i=1}^{n}\|f_i\|}$$
for each $i\in\{1,\cdots,n\}$ and $|t|>R$. Define
$$K_R=\left\{T\in M_n(\IR)_+\mid \|T^{1/2}e_i\|_2\leq R\text{ for each }i\right\}.$$
One has that $|T^{1/2}_{ij}|\leq R$ for any $T\in K_R$, where $T^{1/2}_{ij}$ is the $(i,j)$-th entry of $T^{1/2}$. This shows that $K_R$ is a compact set in $M_n(\IR)_+$. Then by definition, for any $T\in K_R^c$, there exists $i\in\{1,\cdots,n\}$ such that $\|C_{n,e_i}(T)\|_{\Cl(\IR^n)}>R$. As a result, one has that
$$\|f_i(C_{n,e_i}(T))\|\leq \frac{\varepsilon}{1+\prod_{i=1}^{n}\|f_i\|}.$$
As a conclusion, one has that $\|\prod_{i=1}^n f_i(C_{n,e_i})(T)\|\leq\varepsilon$ if $T\in M_n(\IR)_+ \setminus K_R$. This shows that $\prod_{i=1}^n f_i(C_{n,e_i})$ belongs to $\B_0(n)$.

Let $\B'_0(n)$ be the $C^*$-subalgebra of $\B_0(n)$ generated by $\prod_{i=1}^{n} f_i(C_{n, e_i})$ for $f_1,\cdots,f_n\in C_0(\IR)$. Then it suffices to show $\B_0(n)=\B'_0(n)$. To this end, we shall first focus on $C_0(M_n(\IR)_+)$, embedded in the center of $\B_0(n)$ as scalar functions. We claim that $C_0(M_n(\IR)_+)\subseteq\B'_0(n)$. In view of the Stone-Weierstrass Theorem, it suffices to show that there exists a family of elements in $C_0(M_n(\IR)_+)\cap \B'_0(n)$ that separates points in $M_n(\IR)_+$.

To this end, for each $R>0$, we define an even function $f_{0}^{(R)}$ and an odd function $f_{1}^{(R)}$ in $C_0(\IR)$ by
\begin{align*}
f_{0}^{(R)} (x) &= 
\begin{cases}
1 , & |x| \leq R;\\
2 - \frac{|x|}{R}, & R < |x| \leq 2R;\\
0, & |x| > 2R.
\end{cases}\; \quad \text{and} \quad f_{1}^{(R)} (x) = x f_{0}^{(R)} (x) \; .\end{align*}
Hence, for any $i\in\{1,\cdots,n\}$, the element $f_0^{(R)}(C_{n,e_i})$ is in $C_b(M_n(\IR)_+)$, indeed, we have that
\[f_{0}^{(R)} \left(C_{n, e_i}\right) (T) = \begin{cases}
1 , & \left\| T^{1/2} e_i \right\|_2 \leq R; \\
2 - \frac{\left\| T^{1/2} e_i \right\|_2}{R}, & R < \left\| T^{1/2} e_i \right\|_2 \leq 2R; \\
0, & \left\| T^{1/2} e_i \right\|_2 > 2R.\end{cases}\]
for any $T\in M_n(\IR)_+$. It follows that the element
$$g^{(R)}_n=\prod_{i=1}^{n} f_0^{(R)}(C_{n, e_i})$$
takes values on $K_R$ and $0$ outside $K_{2R}$, and is thus in $C_0(M_n(\IR)_+)$. For any $i\in\{1,\cdots,n\}$, observing that $f_1^{(R)}(C_{n,e_i})=f_0^{(R)}(C_{n,e_i})\cdot C_{n,e_i}$, we also define the element
\[C_{n, e_i}^{(R)} = C_{n, e_i} g^{(R)}_n = \prod_{j = 1}^{n} f_{\delta_{ij}}^{(R)} \left(C_{n, e_j}\right) \in\B_0'(n) \; ,\]
where 
\[\delta_{ij} = \begin{cases}1, & i = j; \\
0, & i \not= j.\end{cases}\; \]
It then follows from the algebraic relation $vw + wv = 2 \langle v, w \rangle \cdot 1$ in the Clifford algebra $\Cl(\IR^n)$, where $v, w \in \IR^n$, that for any $i,j \in \{1, \cdots, n\}$ and $T \in M_n(\IR)_{+}$, we have
\begin{align*}
&\left( C_{n, e_{i}}^{(R)} C_{n, e_{j}}^{(R)} + C_{n, e_{j}}^{(R)} C_{n, e_{i}}^{(R)} \right) (T) \\
=& \left(C_{n, e_{i}} (T)C_{n, e_{j}}(T) + C_{n, e_{j}}(T) C_{n, e_{i}} (T) \right) \left(g^{(R)}_n (T)\right)^2 \\
=& 2 \left\langle T^{1/2} e_{i} , T^{1/2} e_{j} \right\rangle \cdot \left(g^{(R)}_n (T)\right)^2 \\
=& 2 \left( e_{i}^T T e_{j} \right) \cdot \left( g^{(R)}_n (T) \right)^2 \; ,
\end{align*}
where $e_i$ is viewed as a column vector. Hence the elements
\[C_{n, e_{i}}^{(R)} C_{n, e_{j}}^{(R)} + C_{n, e_{j}}^{(R)} C_{n, e_{i}}^{(R)} \quad \text{for}\quad i,j \in \{1,\cdots,n\}\]
fall in $C_0(M_n(\IR)_{+}) \cap \B_0'(n)$. 

We claim that this family separates points in $M_n(\IR)_{+}$. Indeed, fixing any distinct $T_1, T_2 \in M_n(\IR)_{+}$, we find $i_0, j_0 \in \{1, \ldots, n\}$ such that $e_{i_0}^{T} T_1 e_{j_0} \not= e_{i_0}^{T}T_2 e_{j_0}$, and set 
\[R = \max\{ \left\| T_k e_i \right\|_2 \mid k =1,2, \ i = 1, \ldots, n \}.\]
Hence we have $T_1, T_2 \in K_R$, which implies that 
\[\left( C_{n, e_{i_0}}^{(R)} C_{n, e_{j_0}}^{(R)} + C_{n, e_{j_0}}^{(R)} C_{n, e_{i_0}}^{(R)} \right) (T_k) = 2 \left( e_{i_0}^{T} T_k e_{j_0} \right) \cdot \left( g_n (T_k) \right)^2 = 2 \left( e_{i_0}^{T}T_k e_{j_0} \right) \cdot 1\]
for $k=1,2$. Therefore the element $C_{n, e_{i_0}}^{(R)} C_{n, e_{j_0}}^{(R)} + C_{n, e_{j_0}}^{(R)} C_{n, e_{i_0}}^{(R)}$ takes different values at $T_1$ and $T_2$, as desired. 
	
To finish the proof that $\B_0'(n) = \B_0(n)$, we can now use the fact that $\B_0'(n)$ is a $M_n(\IR)_{+}$-$C^*$-subalgebra of $\B_0(n)$ and reduce the problem to showing that for any $T \in M_n(\IR)_{+}$, the fiber $\B_0'(n)|_T$ is equal to $\B_0(n)|_T$, which is $\Cl\left(\operatorname{ran}\left({T}^{1/2}\right)\right)$, but this is clear since as long as $R$ is sufficiently large so that $T \in K_R$, we have
\[C_{n, e_{i}}^{(R)} (T) = C_{n, e_{i}} (T) = T^{1/2} e_{i}\]
for $i = 1, \cdots, n$, and thus these elements in $\B_0'(n)|_T$ generate the entire fiber $\B_0(n)|_T$.
\end{proof}

\begin{Def}
Let $M$ be a Hilbert-Hadamard space and let $x_1, \cdots, x_n$ be points in $M$. We define 
\[P_{M;x_1, \cdots, x_n} = \{ A_{x,t; x_1, \cdots, x_n} \mid x \in M, t \in [0,\infty) \} \subseteq M_n(\IR)_{+}.\]
\end{Def}

\begin{Rem}
Observe that as subsets of $M_n(\IR)$, we have
\[P_{M;x_1, \cdots, x_n} = P_{M;x_1, \cdots, x_n} + \IR_{+} \cdot (1)_{n \times n} = P_{M;x_1, \cdots, x_n}^0 + \IR_{+} \cdot (1)_{n \times n}\]	
where 
\[P_{M;x_1, \cdots, x_n}^0 = \left\{ \left( \left\langle \log_x \left( x_i \right), \log_x \left( x_j \right) \right\rangle \right)_{1 \leq i,j \leq n}  \mid x \in M, t \in [0,\infty) \right\}\]
and $(1)_{n \times n}$ is the matrix with every entry equal to $1$. This matrix is a multiple of the rank-$1$ projection onto the diagonal line in $\IR^n$. 
\end{Rem}

\begin{Thm}
Let $M$ be a Hilbert-Hadamard space and let $x_1, \cdots, x_n$ be points in $M$. Then the $*$-homomorphism $\beta_{x_1, \cdots, x_n}:\B(n) \to \Pi_b(M)$ maps $\B_0(n)$ onto $\A(M,\{x_1, \cdots, x_n\})$, with its kernel equal to 
\[C_0 \left(M_n(\IR)_{+} \setminus \overline{P_{M;x_1, \cdots, x_n}} \right) \cdot \B_0(n).\]
\end{Thm}

\begin{proof}
Notice that $C_0 \left(M_n(\IR)_{+} \setminus \overline{P_{M;x_1, \cdots, x_n}} \right) \cdot \B_0(n)$ is an ideal of $\B_0(n)$ with quotient algebra $C(\overline{P_{M;x_1, \cdots, x_n}})\cdot\B_0(n)$. The quotient map is clearly given by the restriction map
$$g\mapsto g|_{\overline{P_{M;x_1, \cdots, x_n}}}.$$
Then the $*$-homomorphism $\beta_{x_1,\cdots,x_n}$ descends to a $C^*$-homomorphism from the quotient algebra $C\left(\overline{P_{M;x_1, \cdots,x_n}}\right)\cdot\B_0(n)$ to $\A(M,\{x_1,\cdots,x_n\})$, it suffices to show this map is an isomorphism.

By Lemma \ref{connection AM and Bn}, one can easily check that
$$\beta_{x_1,\cdots,x_n}\left(f(C_{n,e_i})|_{\overline{P_{M;x_1, \cdots, x_n}}}\right)(x,t)=f(C_{x_i}(x,t))=\beta^M_{x_i}(f)\in\A(M).$$
For a fixed $i\in I$, assume that $f_i\in C_c(\IR)$ satisfying that $\supp(f_i)\subseteq[-R,R]$. Set $R_0=\max\limits_{j\ne i}\{d(x_i,x_j)\}$. We choose $f_j\in C_c(\IR)$ to be an even function such that $f_j(t)=1$ for any $t\in [-R_0-R,R_0+R]$ for each $j\ne i$. One then has that
$$\beta_{x_1,\cdots,x_n}\left(\prod_{i=1}^nf(C_{n,e_i})|_{\overline{P_{M;x_1, \cdots, x_n}}}\right)(x,t)=\beta_{x_i}^M(f_i).$$
Since the set $\left\{\prod^n_{i=1}f_i(C_{n,e_i})|_{\overline{P_{M;x_1, \cdots, x_n}}}\ \mid\ f_i\in C_c(\IR)\right\}$ generates $C\left(\overline{P_{M;x_1, \cdots,x_n}}\right)\cdot\B_0(n)$ by Lemma \ref{generators Bn} and $\left\{\beta_{x_i}(f)\mid f\in C_c(\IR)\right\}$ generates $\A(M,\{x_1,\cdots,x_n\})$, then the Bott map $\beta_{x_1,\cdots,x_n}$ is surjective.

For any $g\in C\left(\overline{P_{M;x_1, \cdots,x_n}}\right)\cdot\B_0(n)$ with $\beta_{x_1,\cdots,x_n}(g)=0$, one has that
$$\beta_{x_1,\cdots,x_n}(g)(x,t)=g(A_{x,t; x_1, \cdots, x_n})=0$$
for each $A_{x,t;x_1,\cdots,x_n}\in P_{M;x_1, \cdots,x_n}$. Since $P_{M;x_1, \cdots,x_n}$ is dense in $\overline{P_{M;x_1, \cdots,x_n}}$, we conclude that $g=0$. This shows that $C\left(\overline{P_{M;x_1, \cdots,x_n}}\right)\cdot\B_0(n)$ is isomorphic to $\A(M,\{x_1,\cdots,x_n\})$.
\end{proof}

As a direct conclusion, $\A(M,\{x_1,\cdots,x_n\})$ is a Type I $C^*$-algebra. By \cite[Lemma 7.2]{GWY2018}, one has that $\A(M)$ is a direct limit of $\A(M,F)$, where $F$ ranges over finite subsets of (an arbitrary dense subset of) $M$. When $M$ is separable, we can arrange it so that $F$ ranges over an increasing sequence of finite subsets that exhausts a countable dense subset $S$ of $M$, whence
$$\A(M)=\overline{\bigcup_{n}\A(M,F_n)} \quad \text{ with } F_1 \subseteq F_2 \subseteq \ldots \text{ finite and } S = \bigcup_n F_n  \; .$$
We then have the following corollary. 

\begin{Cor}\label{AofM in Bootstrap class}
Let $M$ be a separable Hilbert-Hadamard space. Then $\A(M)$ is a direct limit of a sequence of type I $C^*$-algebra, and is thus in the bootstrap class.
\end{Cor}

\section{$KK$-theoretic constructions for twisted localization algebras}\label{sec: KK-theoretic constructions for twisted localization algebras}

In this section, we shall introduce two constructions for twisted localization algebras which are inspired by Kasparov's equivariant $KK$-theory.

\subsection{A $KK$-product construction}\label{KK-product for localization alg}

Let $A, B$ be $\Ga$-$C^*$-algebras, $\Delta$ a compact separable $\Ga$-space, and $1_{C(\Delta)}\in KK^{\Ga}_0(C(\Delta),C(\Delta))$ the unit element. As a remark, we do not assume the $\Ga$-action on $\Delta$ to be isometric here. Then there exists a canonical homomorphism
$$\cdot\ox1_{C(\Delta)}:KK^{\Ga}(A,B)\to KK^{\Ga}(A\ox C(\Delta),B\ox C(\Delta))$$
which is defined by the $KK$-product:
$$[E,\phi,F]\mapsto [E\ox C(\Delta),\phi\ox 1,F\ox 1_{C(\Delta)}].$$
Let $C(\Delta)\rtimes_r\Ga$ be the reduce crossed product of $C(\Delta)$ and $\Ga$. Since there is a canonical inclusion $i: C(\Delta)\to C(\Delta)\rtimes_r\Ga$ defined by $f\mapsto f\cdot \lambda_e$, where $e\in\Ga$ is the unit element. The $\Ga$-action on $C(\Delta)\rtimes_r\Ga$ is given by the inner automorphisms $ad_{\lambda_\gamma}$. Combining the construction above, this induces a homomorphism on $KK$-theory
$$KK^{\Ga}(A,B)\to KK^{\Ga}(A\ox C(\Delta),B\ox(C(\Delta)\rtimes_r\Ga))$$
by
$$[E,\phi,F]\mapsto [E\ox(C(\Delta)\rtimes_r\Ga),\phi\ox i,F\ox 1_{C(\Delta)\rtimes_r\Ga}].$$

In this subsection, we shall introduce an analog of this homomorphism for twisted localization algebras. To be more specific, let $X$ be a proper $\Ga$-space, $(N, M,f)$ a coefficient system for $\Ga\car X$, we can then define the twisted Roe algebra and twisted localization algebra. We would like to construct a $C^*$-algebra $C^*_L(X\times\Delta,\A(M)\ox(C(\Delta)\rtimes_r\Ga))^{\Ga}$ which is an analog of twisted localization algebra for $X\times\Delta$ with coefficient in $\A(M)\ox (C(\Delta)\rtimes_r\Ga)$ and a map
$$K_*(C^*_L(X,\A(M))^{\Ga})\to K_*(C^*_L(X\times\Delta,\A(M)\ox(C(\Delta)\rtimes_r\Ga))^{\Ga}).$$

Since we are considering the localization algebra of $X\times\Delta$, we equip $X\times\Delta$ with an invariant metric as in \cite[Theorem 4.3.4]{metrizationGaspace}. As we discussed in Remark \ref{metric on product space} (also in Appendix \ref{Appendix A}), $K$-theory of the localization algebra does not depend on the choice of the metrization.
Inspired by Kasparov's $KK$-theory, it is natural to define the map to be
$$g\mapsto g\ox 1,$$
where $g\in C^*_L(X,\A(M))^{\Ga}$ and $g\ox 1_{C(\Delta)\rtimes_r\Ga}$ is viewed a function from $[0,\infty)$ to the $C^*$-algebra of all adjointable module homomorphism $\L\left(E^X_{\A(M)}\ox(C(\Delta)\rtimes_r\Ga)\right)$, where
$$E^X_{\A(M)}=\ell^2(N)\ox\ell^2(\Ga)\ox\H\ox\A(M),$$
and $E^X_{\A(M)}\ox(C(\Delta)\rtimes_r\Ga)$ is viewed as a Hilbert $\A(M)\ox(C(\Delta)\rtimes_r\Ga)$-module.
However, one should notice that this Hilbert module is not \emph{absorbing} (also called \emph{ample}, cf. \cite[Definition 2.1]{DWW2018}) since the representation $\sigma: C(\Delta)\to\L(C(\Delta)\rtimes_r\Ga)$ is not absorbing, where $\sigma$ is defined by
$$\sigma(f_1)(f_2\cdot\lambda_\gamma)=(f_1f_2)\cdot\lambda_\gamma\quad\text{for any }f_1,f_2\in C(\Delta),\gamma\in\Ga,$$
so that may bring some troubles when we consider the functorial properties. To solve this problem, one basic idea is to absorb $\sigma$ by a "larger" representation.

Since $C(\Delta)\rtimes_r\Ga$ is unital, $\L(C(\Delta)\rtimes_r\Ga)$ is isomorphic to $C(\Delta)\rtimes_r\Ga$. One can see that $\sigma$ is also given by $f\mapsto f\cdot\lambda_e\in C(\Delta)\rtimes_r\Ga$ for any $f\in C(\Delta)$. Choose a $\Ga$-invariant countable dense subset $W\subseteq \Delta$. For each $n\in\IN$, choose $\{W_n\}$ to be a sequence of subsets of $W$ such that $W_n$ is a $\frac 1n$-net, $W_n\subseteq W_{n+1}$ for each $n\in\IN$. Since $\Delta$ is compact, we shall assume that each $W_n$ is \emph{finite}. Let $\{\phi_w^n\}_{w\in W_n}$ be an $\ell^2$-partition of unity associated with the open cover $\{B(w,\frac 1n)\}_{w\in W_n}$. Moreover, we also assume that $W=\bigcup_{n\in\IN}\Ga\cdot W_n$. We should mention that the family $\{B(\gamma\cdot w,\frac 1n)\}_{w\in W_n}$ also forms an open cover for $\Delta$ and $\{\gamma\cdot\phi_w^n\}$ forms an associated $\ell^2$-partition of unity for each $\gamma\in\Ga$. Notice that $C(\Delta)\rtimes_r\Ga$ has a faithful representation on $\ell^2(\Ga)\ox\ell^2(W)$ given by
$$(\lambda_{\gamma})(\delta_{\eta}\ox\delta_w)=\delta_{\gamma\eta}\ox\xi\quad\text{and}\quad f(\delta_\eta\ox\delta_w)=f(\eta w)(\delta_\eta\ox \delta_w),$$
for any $\gamma,\eta\in\Ga$, $w\in W$ and $f\in C(\Delta)$.
We still write it $\sigma$ the composition of $\sigma:C(\Delta)\to C(\Delta)\rtimes_r\Ga$ and the faithful representation $C(\Delta)\rtimes_r\Ga\to\B(\ell^2(\Ga)\ox\ell^2(W))$, i.e.,
$$\sigma(f)(\delta_\gamma\ox\xi)=\delta_\gamma\ox(\gamma^{-1}\cdot f)\xi$$
where $(\gamma^{-1}\cdot f)(x)=f(\gamma x)$ and it acts on $\xi$ by pointwise multiplication for $\xi\in\ell^2(W)$.

\begin{Def}\label{absorbing representation}
Let $A$ be a $C^*$-algebra, let $H_1,H_2$ be Hilbert spaces. Let $\pi_i: A\to\B(H_i)$ be unital representations for $i=1,2$. We say that \emph{$\pi_1$ absorbs $\pi_2$} if there exists a sequence of isometries $\{V_n:H_2\to H_1\}$ such that
$$\|V_n\pi_2(a)-\pi_1(a)V_n\|\to 0\quad \text{as}\quad n\to\infty$$
for any $a\in A$.
\end{Def}

Let $H_{\Delta}=\ell^2(W)\ox(\ell^2(\Ga)\ox\ell^2(W))$. Define $\pi:C(\Delta)\to \B(H_{\Delta})$ by
$$\pi(f)(\delta_{w_1}\ox \delta_\gamma\ox\delta_{w_2})=f(w_1)\cdot(\delta_{w_1}\ox \delta_\gamma\ox\delta_{w_2}),$$
for any $f\in C(\Delta)$, $w_1,w_2\in W$, and $\gamma\in \Ga$. We now show that $\pi:C(\Delta)\to \B(H_{\Delta})$ can absorb $\sigma:C(\Delta)\to\B(\ell^2(\Ga)\ox\ell^2(W))$.

For each $n\in\IN$, define $V_n:\ell^2(\Ga)\ox\ell^2(W)\to H_{\Delta}$ by
$$\delta_\gamma\ox\xi\mapsto \sum_{w\in W_n}\delta_{\gamma w}\ox\delta_\gamma\ox\phi_w^n\cdot\xi.$$
for each $f\in C(\Delta)$, $\gamma\in\Ga$. One can check that $V_n$ is an isometry by a simple calculation:
\[\begin{split}\langle V_n(\delta_\gamma\ox\xi),V_n(\delta_\gamma\ox\xi)\rangle&=\left\langle\sum_{w\in W_n}\delta_{\gamma w}\ox\delta_\gamma\ox\phi_w^n\cdot\xi,\sum_{w\in W_n}\delta_{\gamma w}\ox\delta_\gamma\ox\phi_w^n\cdot\xi\right\rangle
\\&=\sum_{w\in W_n}\langle(\phi^n_w)^2\cdot\xi,\xi\rangle=\langle\delta_\gamma\ox\xi,\delta_\gamma\ox\xi\rangle.\end{split}\]
It is also easy to check that
$$V_n^*(\delta_{w}\ox\delta_\gamma\ox\xi)=\left\{\begin{aligned}&\delta_\gamma\ox\phi_{\gamma^{-1} w}^n\cdot\xi,&&\gamma^{-1} w\in W_n\\&0,&&\text{otherwise}.\end{aligned}\right.$$
We claim that $\pi$ absorbs $\sigma$ through the sequence of isometries $\{V_n\}_{n\in\IN}$. Indeed, for each $f\in C(\Delta)$ and $\delta_\gamma\ox\xi\in\ell^2(\Ga)\ox\ell^2(W)$, we have that
$$(V_n\circ \sigma(f))(\delta_\gamma\ox\xi)=V_n(\delta_\gamma\ox (\gamma^{-1}f)\cdot\xi)=\sum_{w\in W_n}\delta_{\gamma w}\ox\delta_\gamma\ox\phi_w^n(\gamma^{-1} f)\cdot\xi;$$
$$(\pi(f)\circ V_n)(\delta_\gamma\ox\xi)=\pi(f)\left(\sum_{w\in W_n}\delta_{\gamma w}\ox\delta_\gamma\ox\phi_w^n \cdot\xi\right)=\sum_{w\in W_n}\delta_{\gamma w}\ox\delta_\gamma\ox f(\gamma w)\phi_w^n\cdot\xi.$$
Thus
$$\|(V_n\circ\sigma(f_1)-\pi(f_1)\circ V_n)(\delta_\gamma\ox\xi)\|=\left\|\sum_{w\in W_n}\delta_{\gamma w}\ox\delta_\gamma\ox\phi_w^n(\gamma^{-1} f-f(\gamma w))\cdot\xi\right\|.$$
Since $\supp(\phi_{w}^n)\subseteq B(w,\frac 1n)$ and $\gamma^{-1}f$ is uniformly continuous, for any $\varepsilon>0$, there exists $N>0$ such that
$$|(\gamma^{-1}f)(x)-f(\gamma w)|=|f(\gamma x)-f(\gamma w)|<\varepsilon$$
for any $x\in B(w,\frac 1n)$ and $n>N$. Thus, we have that
\[\begin{split}\|(V_n\circ\sigma(f_1)-\pi(f_1)\circ V_n)(\delta_\gamma\ox\xi)\|&\leq\varepsilon\cdot\left\|\sum_{w\in W_n}\delta_{\gamma w}\ox\delta_\gamma\ox\phi_w^n\cdot\xi\right\|=\varepsilon\cdot\|\delta_{\gamma}\ox\xi\|\end{split}\]
for any $n>N$.
This shows that
$$\|V_n\circ\sigma(f_1)-\pi(f_1)\circ V_n\|\to 0\quad\text{as}\quad n\to\infty.$$

Let $p: X\times\Delta\to X$ be the canonical projection onto $X$, which is clearly $\Ga$-equivariant. Let $W\subseteq \Delta$ be the countable subset as above, $(N,M,f)$ the coefficient system for $\Ga\car X$. Then $(N\times W,M,f\circ p)$ forms a coefficient system for $\Ga\car X\times\Delta$. Thus we can define the \emph{algebraic twisted Roe algebra} $\IC[X\times\Delta,\A(M)]^{\Ga}$. We denote by
$$\IC[X\times\Delta,\A(M)]^{\Ga}\ox_{alg}(C(\Delta)\rtimes_r\Ga)$$
the algebraic tensor of $\IC[X\times\Delta,\A(M)]^{\Ga}$ and $C(\Delta)\rtimes_r\Ga$. Notice that it has a faithful representation on the Hilbert $\A(M)$-module $E=E^{X\times\Delta}_{\A(M)}\ox (\ell^2(\Ga)\ox\ell^2(W))$, where the representation is given by the tensor product of the canonical representation of $\IC[X\times\Delta,\A(M)]^{\Ga}$ on
$$E^{X\times\Delta}_{\A(M)}=\ell^2(N)\ox\ell^2(W)\ox\ell^2(\Ga)\ox\H\ox\A(M)$$
and the canonical representation of $C(\Delta)\rtimes_r\Ga$ on $\ell^2(\Ga)\ox\ell^2(W)$. One can also see that $E= E^X_{\A(M)}\ox H_\Delta$.
We denote by 
$$C^*(X\times\Delta,\A(M))^{\Ga}\ox (C(\Delta)\ox\Ga)$$
the completion of $\IC[X\times\Delta,\A(M)]^{\Ga}\ox_{alg}(C(\Delta)\rtimes_r\Ga)$ under the norm given by this representation

For any $T=\sum_{i=1}^N T_i\ox a_i\in\IC[X\times\Delta,\A(M)]^{\Ga}\ox (C(\Delta)\rtimes_r\Ga)$, the propagation of $T$, denoted by $\prop(T)$, is defined to be $\max\{\prop(T_i)\mid i=1,\cdots,N\}$. Moreover, we can still view $T$ as a $(N\times W)$-by-$(N\times W)$ matrix. For any $x,y\in N\times W$, the "$(x,y)$-th matrix coefficient of $T$", denoted by $T(x,y)$, takes value in 
$$\A(M)\ox\K(\ell^2(\Ga)\ox\H)\ox C(\Delta)\rtimes_r\Ga.$$
We define the support of $T(x,y)$ to be $\bigcup_{i=1}^N\supp(T_i(x,y))$ as a subset of $M\times \IR_+$.

\begin{Def}
Define $C^*_L(X\times\Delta,\A(M)\ox(C(\Delta)\rtimes_r\Ga))^{\Ga}$ to be the $C^*$-algebra generated by all bounded, uniformly continuous functions
$$g:\IR_+\to\IC[X\times\Delta,\A(M)]^{\Ga}\ox(C(\Delta)\rtimes_r\Ga)$$
such that\begin{itemize}
\item[(1)] (propagation condition) $\prop(g(s))$ tends to $0$ as $t$ tends to infinity;
\item[(2)] (twist condition) there exists $R> 0$ such that $\supp((g(s))(x,y))\subseteq B(f\circ p(x),R)\subseteq M\times\IR_+$ for all $s\in\IR_+$ and $x,y\in N\times W$.
\end{itemize}\end{Def}

\begin{Rem}\label{imitating}
Imitating the definition above, one can also define $C^*_L(X,\A(M)\ox(C(\Delta)\rtimes_r\Ga))^{\Ga}$ for any proper metric space $X$ with a proper $\Ga$-action and a coefficient system $(N,M,f)$. It is generated by all functions
$$g:\IR_+\to\IC[X,\A(M)]^{\Ga}\ox (C(\Delta)\rtimes_r\Ga)$$
satisfying that $\prop(g(t))$ tends to $0$ as $t$ tends to $\infty$ and the twist condition as above.
\end{Rem}

Let $\{V_n:\ell^2(\Ga)\ox\ell^2(W)\to H_{\Delta}\}$ be the sequence of isometries as before. For each $t\in[n,n+1]$, define $V_t:\Big(\ell^2(\Ga)\ox\ell^2(W)\Big)\oplus\Big(\ell^2(\Ga)\ox\ell^2(W)\Big)\to H_{\Delta}\oplus H_{\Delta}$ by
$$V_t=R_t^*(V_n\oplus V_{n+1})R_t$$
where
$$R_t=\begin{pmatrix}\cos(\pi\theta/2)&\sin(\pi\theta/2)\\-\sin(\pi\theta/2)&\cos(\pi\theta/2)\end{pmatrix}$$
and $\theta=t-n$. Notice that this family $\{V_t\}_{t\geq 0}$ is continuous by $t$. To simplify the notation, we shall also denote it by $V_t$ the isometry $1\ox V_t: E^X_{\A(M)}\ox\ell^2(\Ga)\ox\ell^2(W)\to E^{X}_{\A(M)}\ox H_\Delta=E$.

Define $\Phi:C^*_L(X,\A(M))^{\Ga}\to C^*_L(X\times\Delta,\A(M)\ox C(\Delta))^{\Ga}\ox M_2(\IC)$ by
$$\Big(\Phi(g)\Big)(t)\mapsto V_t\Big((g(t)\ox 1_{C(\Delta)\rtimes_r\Ga})\oplus 0\Big)V_t^*.$$
Here are some explanations. For each $t\geq 0$, $g(t)\ox 1_{C(\Delta)\rtimes_r\Ga}$ forms an adjointable endomorphism on $E^X_{\A(M)}\ox\ell^2(\Ga)\ox\ell^2(W)$. For each $n\in\IN$, one can see that $V_nV^*_n$ forms a bounded operator on $H_\Delta$. Fix $w_0\in W$. For any $\gamma\in\Ga$ and $\xi\in\ell^2(W)$, if $w_0\in W_n$, one has that
\begin{equation}\label{VNVN*}V_nV_n^*(\delta_{\gamma w_0}\ox\delta_\gamma\ox\xi)=\sum_{w\in W_n}\delta_{\gamma w}\ox\delta_\gamma\ox\phi^n_{w_0}\phi_w^n\cdot\xi.\end{equation}
Thus $V_nV^*_n$ can be seen as an operator on $\ell^2(\Ga\cdot W_n)\ox\ell^2(\Ga)\ox\ell^2(W)\subseteq H_\Delta$. From \eqref{VNVN*}, if we fix $w_1,w_2\in W_n$, one can see that the "$(\gamma w_1,\gamma w_2)$-th matrix coefficient" of $V_nV_n^*$ contains $\phi_{w_1}^n\phi^n_{w_2}\in C(\Delta)\rtimes_r\Ga$ as an addend for any $\gamma\in\Ga$. If there are $w_3,w_4\in W_n$ and $\gamma'\in\Ga$ such that $\gamma'w_3=\gamma w_1$ and $\gamma'w_4=\gamma w_2$, one can see that $\phi_{w_1}^n\phi^n_{w_2}$ is also an addend for the "$(\gamma w_3,\gamma w_4)$-th matrix coefficient" of $V_nV_n^*$.

Write $W_n=\{w_1,\cdots,w_N\}$. Then
$$(\Phi(g))(n)=\sum_{i,j=1}^N T^n_{ij}\ox\phi_{w_i}^n\phi_{w_j}^n$$
where $T^n_{ij}\in\IC[X\times\Delta,\A(M)]^{\Ga}$ is defined by
$$T^n_{ij}\Big((x_1,z_1),(x_2,z_2)\Big)=\left\{\begin{aligned}&T(x_1,x_2),&&\exists\gamma\in\Ga\text{ such that }\gamma z_1=w_i,\gamma z_2=w_j;\\&0,&&\text{otherwise.}\end{aligned}\right.$$
It is clear that $T^n_{ij}$ is a well-defined element of $\IC[X\times\Delta,\A(M)]^{\Ga}$.

Write $\prop_X(T^n_{ij})$ and $\prop_\Delta(T^n_{ij})$ for the propagation of $T^n_{ij}$ considered with respect to the $X$-module structure and the $\Delta$-module structure respectively (here the $\Delta$-propagation is associated with the original metric on $\Delta$). By definition, one has that $\prop_X(T^n_{ij})\leq\prop(g(n))$. Since $\diam(\supp(\phi_w^n))\leq\frac 1n$, one has that $\prop_\Delta(T^n_{ij})\leq\frac 2n$. Since the metric on $X\times\Delta$ induces the product topology, we conclude that the propagation of $(\Phi(g))(n)$ tends to $0$ as $n$ tends to infinity. This shows that $\Phi(g)$ is well-defined. It induces a homomorphism on $K$-theory:
$$\Phi_*:K_*(C^*_L(X,\A(M))^{\Ga})\to K_*(C^*_L(X\times\Delta,\A(M)\ox(C(\Delta)\rtimes_r\Ga))^{\Ga}).$$

\subsection{A K\"unneth theorem for twisted localization algebras}
\label{subsec:Kunneth}

In this subsection, we shall prove a K\"unneth formula for $K_*(C^*_L(X,\A(M)\ox(C(\Delta)\rtimes_r\Ga))^{\Ga})$ by using the result from Section \ref{Sec: Approximating AofM}, where the twisted localization algebra is defined as above, see Remark \ref{imitating}. To be specific, let $(N,M,f)$ be a coefficient system for $X$. For any $t\in [0,1]$, there is an associated coefficient system $(N,M^{[0,1]},f_t)$ as we discussed in Section \ref{sec: Twisted algebras and the Bott homomorphisms}. In particular, when $t=0$, the map $f_0$ is a constant map onto the based point of $M^{[0,1]}$ and the group action $\alpha_0$ on $M^{[0,1]}$ is the trivial action.

Imitating Remark \ref{imitating} and Definition \ref{twist Roe big}, one can similarly define this enlarged version of the twisted localization algebra $C^*_L(X,\A_{[0,1]}(M)\ox(C(\Delta\rtimes_r\Ga)))^{\Ga}$ and the evaluation map
$$ev_t:C^*_L\left(X,\A_{[0,1]}(M)\ox(C(\Delta\rtimes_r\Ga))\right)^{\Ga}\to C^*_L\left(X,\A(M^{[0,1]})\ox(C(\Delta\rtimes_r\Ga))\right)^{\Ga,(t)},$$
where the right-hand side is the twisted localization algebra associated with $(N,M^{[0,1]},f_t)$. Then the proof of Lemma \ref{deformation} applies verbatim to the following lemma:

\begin{Lem}\label{independent on t}
Let $X$ be a free and proper metric space that admits a uniformly locally finite $\F$-cover $\U$ such that each element in $\U$ is contractible, $(N.M,f)$ a coefficient system. Then for each $t\in[0,1]$, the $K$-theortic map
$$(ev_t)_*:K_*\left(C^*_L\left(X,\A_{[0,1]}(M)\ox(C(\Delta\rtimes_r\Ga))\right)^{\Ga}\right)\to K_*\left(C^*_L\left(X,\A(M^{[0,1]})\ox(C(\Delta\rtimes_r\Ga))\right)^{\Ga,(t)}\right)$$
is an isomorphism.
\end{Lem}

Our main theorem of this section is stated as follows.

\begin{Thm}\label{Kunneth formula}
Let $X$ be a proper metric space with a free and proper $\Ga$-action. Assume that $X$ admits a uniformly locally finite $\F$-cover $\U$ such that each element in $\U$ is contractible. Then the $KK$-product
$$\bigoplus\limits_{i+j=*} K_{i}\Big( C_L^*\left(X,\mathcal{A}(M^{[0,1]}) \right)^{\Gamma,(t)} \Big) \otimes K_{j} \Big( C(\Delta) \rtimes_r \Gamma  \Big)\to K_*\Big(C^*_L\left(X,\A(M^{[0,1]})\ox(C(\Delta)\rtimes_r\Ga)\right)^{\Ga,(t)}\Big).$$
is a rational isomorphism for each $t\in [0,1]$.
\end{Thm}

One should compare this theorem with Theorem \ref{kunneth}. We should mention that all the functorial properties we discussed in Appendix \ref{Appendix A} still hold for the functor
$$X\mapsto K_*\Big(C^*_L\left(X,\A(M^{[0,1]})\ox(C(\Delta)\rtimes_r\Ga)\right)^{(0)}\Big)$$
for any locally compact, Hausdorff, second countable space $X$ with trivial $\Ga$-action.

\begin{proof}
For any  $g\in C^*_L(X,\A(M^{[0,1]}))^{\Ga,(t)}$ and $a\in C(\Delta)\rtimes_r\Ga$, their tensor product $g\ox a$ clearly defines an element in $C^*_L(X,\A(M^{[0,1]})\ox(C(\Delta)\rtimes_r\Ga))^{\Ga,(t)}$. This gives a canonical inclusion
$$C^*_L\left(X,\A(M^{[0,1]})\right)^{\Ga,(t)}\ox (C(\Delta)\rtimes_r\Ga)\mapsto C^*_L\left(X,\A(M^{[0,1]})\ox(C(\Delta)\rtimes_r\Ga)\right)^{\Ga,(t)}.$$
Consider the homomorphism given by the composition
\[\begin{split}\bigoplus\limits_{i+j=*} K_{i}\Big( C_L^*\left(X,\mathcal{A}(M^{[0,1]}) \right)^{\Gamma,(t)} \Big) \otimes K_{j}\Big( C(\Delta) \rtimes_r \Gamma  \Big)&\to K_*\Big(C^*_L\left(X,\A(M^{[0,1]})\right)^{\Ga,(t)}\ox(C(\Delta)\rtimes_r\Ga)\Big)\\
&\to K_*\Big(C^*_L\left(X,\A(M^{[0,1]})\ox(C(\Delta)\rtimes_r\Ga)\right)^{\Ga,(t)}\Big).\end{split}\]
where the first map is given by Kasparov $KK$-product (cf. \cite{Kas88} or \cite[Section 2.10]{HIT2020} for some details of the external product in $K$-theroy) and the second map is induced by the canonical inclusion. We shall then prove this map is an isomorphism after tensoring both sides by $\IQ$.

Notice that $X$ admits a uniformly locally finite $\F$-cover $\U$ such that each element in $\U$ is contractible, by Lemma \ref{independent on t}, it suffices to show the $KK$-product above is a rational isomorphism for each $t=0$. In this case, the $\Ga$-action on $\A(M^{[0,1]})$ is trivial. Thus one can show that
$$K_*\left(C^*_L\left(X,\A(M^{[0,1]})\ox(C(\Delta)\rtimes_r\Ga)\right)^{\Ga,(0)}\right)\cong K_*\left(C^*_L\left(X/\Ga,\A(M^{[0,1]})\ox(C(\Delta)\rtimes_r\Ga)\right)^{(0)}\right),$$
$$K_*\left(C^*_L\left(X,\A(M^{[0,1]})\right)^{\Ga,(0)}\right)\cong K_*\left(C^*_L\left(X/\Ga,\A(M^{[0,1]})\right)^{(0)}\right),$$
by using an argument similar to \cite[Theorem 6.5.15]{HIT2020}, where the coefficient system is given by $(N/\Ga,M^{[0,1]},f_Q)$ and $f_Q$ is the constant map from $N/\Ga$ to the base point $m$ of $M$. Thus it suffices to show
$$\bigoplus\limits_{i+j=*} K_{i}\Big( C_L^*\left(Q,\mathcal{A}(M^{[0,1]}) \right)^{(0)} \Big) \otimes K_{j} \Big( C(\Delta) \rtimes_r \Gamma  \Big)\to K_*\Big(C^*_L\left(Q,\A(M^{[0,1]})\ox(C(\Delta)\rtimes_r\Ga)\right)^{(0)}\Big),$$
is a rational isomorphism for any proper metric space $Q=X/\Ga$ with a trivial $\Ga $-action.

First, we consider the case when $Q$ is compact. When $Q$ is a single point, by using the evaluation map, it suffices to show
$$\bigoplus\limits_{i+j=*} K_{i}\left(\mathcal{A}(M^{[0,1]})\right) \otimes K_{j} \Big( C(\Delta) \rtimes_r \Gamma  \Big)\to K_*\Big(\A(M^{[0,1]})\ox(C(\Delta)\rtimes_r\Ga)\Big),$$
is an isomorphism, where the map is given by the external $KK$-product. This certainly holds since $\A(M^{[0,1]})$ is in the Bootstrap class by Corollary \ref{AofM in Bootstrap class}. Since both the functors
$$X\mapsto \bigoplus\limits_{i+j=*} K_{i}\Big( C_L^*\left(X,\mathcal{A}(M^{[0,1]}) \right)^{(0)} \Big) \otimes K_{j} \Big( C(\Delta) \rtimes_r \Gamma  \Big)\ox\IQ,$$
$$X\mapsto K_*\Big(C^*_L\left(X,\A(M^{[0,1]})\ox(C(\Delta)\rtimes_r\Ga)\right)^{(0)}\Big)\ox\IQ$$
satisfy the Mayer-Vietoris theorem, by using a cutting and pasting argument and the Five Lemma, the theorem for the case when $X$ is a finite CW-complex. Since each compact Hausdorff space can be viewed as an inverse limit of finite CW-complexes, the theorem holds for the case when $X$ is compact.

For the general case, let $X^+$ be the one-point compactification of $X$. Then by using an analogue of Theorem \ref{excision}, one has a six-term exact sequence for both functors. Thus the theorem holds by using the Five Lemma.
\end{proof}

\section{Proof of the main theorem}\label{sec: Proof of the main theorem}

In this section, we shall use the constructions in the previous sections to prove Theorem \ref{main result with torsion}. 
Compared with the torsion-free case (\Cref{main result torsion-free}), 
the result in the general setting requires 
a key extra ingredient, namely to show the canonical map $\pi:\MRdnX\to P_d(X)$ induces a rational injection 
\[
    \pi_* \colon \lim_{d,n\to\infty} K_*(C^*_L(\MRdnX,\A(M))^{\Ga}) \to \lim_{d\to\infty} K_*(C^*_L(P_d(X),\A(M))^{\Ga})
\]
on the $K$-theory of the twisted localization algebras. 
In the special case that $\Gamma \curvearrowright X$ is cocompact, 
should we drop the coefficient algebra $\A(M)$ in the above, then by \Cref{rem:Milnor-Ripe-classifying}, it would amount to proving the canonical map $\pi_* \colon RK^{\Ga}_*(\E\Ga) \to RK^{\Ga}_*(\EG)$ is a rational injection, a fact already known in \cite[Section~7]{BCH1994} via group homological computations. 
When we include the coefficient algebra $\A(M)$ but still require cocompactness for $\Gamma \curvearrowright X$, the desired rational injectivity can be proved using $KK$-theory with real coefficients (see \cite[Section~6]{AAS2020}). 
Using the machinery we have developed for twisted localization algebras, we shall, in \Cref{pi is injective} below, establish rational injectivity without the cocompactness assumption, which will allow us to complete our proof of Theorem \ref{main result with torsion}.

Consider the following commuting diagram:
\begin{equation}\begin{tikzcd}\label{main diam}
K_{*+1}(C^*_L(\MRdnX)^{\Ga}) \arrow[r, "\pi_*", "(2)"'] \arrow[rr,  bend left=10] \arrow[d, "(\widetilde{\beta_L})_*", "(5)"'] & K_{*+1}(C^*_L(P_d(X))^{\Ga}) \arrow[r, "ev_*","(1)"'] \arrow[d, "(\beta_L)_*", "(4)"'] & K_{*+1}(C^*(P_d(X))^{\Ga}) \arrow[d, "\beta_*", "(3)"']   \\
K_*(C^*_L(\MRdnX,\A_{[0,1]}(M))^{\Ga}) \arrow[r, "\pi_*","(10)"'] \arrow[d, "(\widetilde{ev_1})_*", "(8)"']  & K_*(C^*_L(P_d(X),\A_{[0,1]}(M))^{\Ga})\arrow[r, "ev_*","(9)"'] \arrow[d, "(ev_1)_*","(7)"']& K_*(C^*(P_d(X),\A_{[0,1]}(M))^{\Ga}) \arrow[d, "(ev_1)_*","(6)"']\\
K_*(C^*_L(\MRdnX,\A(M^{[0,1]}))^{\Ga,(1)}) \arrow[r, "\pi_*","(12)"']& K_*(C^*_L(P_d(X),\A(M^{[0,1]}))^{\Ga,(1)})\arrow[r, "ev_*","(11)"'] & K_*(C^*(P_d(X),\A(M^{[0,1]}))^{\Ga,(1)})
\end{tikzcd}\end{equation}
We have introduced the right two columns of this diagram in Section \ref{sec: Twisted algebras and the Bott homomorphisms}. The map (5) and (8) in the diagram are defined similarly with map (4) and (7), respectively. The maps $(2),(10),(12)$ are defined by using the functoriality properties of the localization algebras. Since $\pi:\MRdnX\to P_d(X)$ is an equivariant continuous coarse equivalence, one can find a family of equivariant continuous covering isometries $(V_t)_{t\geq 0}$ for $\pi$ as in \cite[Proposition 4.5.12]{HIT2020} such that it induces homomorphisms on $K$-theory:
$$\pi_*=(ad_{V_t})_*:K_*(C^*_L(\MRdnX)^{\Ga})\to K_*(C^*_L(P_d(X))^{\Ga})$$
$$\pi_*=(ad_{V_t})_*:K_*(C^*_L(\MRdnX,\A_{[0,1]}(M))^{\Ga})\to K_*(C^*_L(P_d(X),\A_{[0,1]}(M))^{\Ga})$$
$$\pi_*=(ad_{V_t})_*:K_*(C^*_L(\MRdnX,\A(M^{[0,1]}))^{\Ga,(1)})\to K_*(C^*_L(P_d(X),\A(M^{[0,1]}))^{\Ga,(1)})$$
and these homomorphisms do not depend on the choice of $(V_t)_{t\geq0}$.

\begin{Lem}\label{pi is injective}
The homomorphism
$$\pi_*:\lim_{d,n\to\infty}K_*(C^*_L(\MRdnX,\A(M^{[0,1]}))^{\Ga,(1)})\to \lim_{d\to\infty}K_*(C^*_L(P_d(X),\A(M^{[0,1]}))^{\Ga,(1)})$$
is rationally injective.
\end{Lem}

\begin{proof}
Let $\Omega_{X}$ be the $\Ga$-space with property TAF as in \Cref{TAF space}, which consists of all linear orders on $X$. For any $d>0$ and $n\in\IN$, equip $X\times\Omega_X$ with the twist metric as we discussed in Section \ref{sec: A Concrete model for space with Property TAF}. Thus the canonical quotients
$$P_d(X)\times\Omega_X\to P_d(X)\quad\text{ and }\quad \MRdnX\times\Omega_X\to\MRdnX$$
are both equivariant coarse equivalence. Moreover, the $\Ga$-action on $P_d(X)\times\Omega_X$ and $\MRdnX\times\Omega_X$ are both free and proper. Let $\U$ be a uniformly locally finite $\F$-cover of $P_d(X)$. For each $U\in\U$, there exists a finite subgroup $\Ga_U\in\F$ such that $U$ is $\Ga_U$-invariant. Since $\Omega_X$ is compact and $\Gamma_U$ is finite, one can choose a finite open cover $\V_U$ on $\Omega_X$ such that for any $V\in\V_U$, $V\cap \gamma V=\emptyset$ for each $\gamma\in\Ga_U$. Then $\widetilde{\U}=\{U\times V\}_{V\in\V_U,U\in\U}$ is clearly a uniformly locally finite $\F$-cover for $P_d(X)\times\Omega_X$. Thus we have that $P_d(X)\times\Delta_{\Ga}\in\F\P(X,\Ga)$. Similarly, one also has that $\MRdnX\times\Delta_{\Ga}\in\F\P(X,\Ga)$. 

Denote $n=\max\limits_{x\in X}\#B(x,d)$ for each $d>0$. For any $\R\in \Omega_X$ and $F\subset X$ with $\diam(F)\leq d$, one can write
$$F=\{x_{\R,1},\cdots,x_{\R,k}\}$$
such that $x_{\R,i}<_{\R}x_{\R,j}$ for any $i< j$, where $k=\#F\leq n$. Then there is also a canonical ordering map $P_d(X)\times\Omega_X\to\MRdnX$ defined by 
\begin{equation}\label{eq:ordering_map}
    \left(\sum_{x\in F}t_xx,\R\right)\mapsto (t_{1}x_{\R,1},\cdots,t_{k}x_{\R,k},0,\cdots),
\end{equation}
where $t_i=t_{x_{\R,i}}$ for each $i$. This ordering map is continuous. Indeed, for any $\varepsilon<\min_{x\in F}{t_x}$, the set $W_{\varepsilon}$ of all points $\langle s,z\rangle\in\MRdnX$ satisfying that $|s_i-t_i|<\varepsilon$ and $z_i=x_i$ for any $1\leq i\leq k$ forms a open neighbourhood of $(t_{1}x_{\R,1},\cdots,t_{k}x_{\R,k},0,\cdots)$. Then
$$U=\left\{\sum_{x\in D}s_xx\ \Big|\ \diam(D)\leq d,F\subseteq D,|s_x-t_x|\leq \frac{\varepsilon}{2}\right\}$$
forms a open neighbourhood of $\sum_{x\in F}t_xx$. Let $V$ be the set of all $\R'\in\Omega_X$ satisfying that $\R'|_F=\R_F$ and $x<_{\R'}z$ for any $z\in D\backslash F$ and $x\in F$, where $D\subset X$ with $\diam(D)\leq d$ and $F\subset D$. Thus the image of $U\times V$ under the ordering map falls in $W_{\varepsilon}$, which means the ordering map is continuous. It is easy to check that the ordering map is also an equivariant coarse equivalence.

Similarly, one can also define a collapsing map
$$\MRdnX\times\Omega_{X}\to\MRdnX.$$
For any $\langle s,z\rangle$ and $\R\in\Omega_X$, write $\supp(\langle s,z\rangle)=\{x_{\R,1},\cdots,x_{\R,k}\}$ as above. The collapsing map is defined by
$$(\langle s,z\rangle,\R)\mapsto (t_{1}x_{\R,1},\cdots,t_{k}x_{\R,k},0,\cdots),$$
where $t_i=\sum_{z=x_{\R_i}} s_z$. The collapsing map is also an equivariant continuous coarse equivalence. We then have the following commuting diagram
\begin{equation}\label{collapsing-ordering}\begin{tikzcd}
P_d(X)\times\Omega_X \arrow[rd, "\text{ordering}"] \arrow[rrd, bend left=15,"\Ga-\text{homotopy}"] &                            &   \\
\MRdnX\times\Omega_X \arrow[r, "\text{collapsing}"] \arrow[u, "\pi"]       & \MRdnX \arrow[r, "{d<d',n<n'}"] & \widetilde{P_{d',n'}}(X),
\end{tikzcd}\end{equation}
where there exists $d'>d,n'>n$ such that the composition of the collapsing map and the inclusion is $\Ga$-homotopy to the canonical inclusion $\MRdnX\to\widetilde{P_{d',n'}}(X)$ by Theorem \ref{classify free and proper}.

By using the functoriality of $K$-theory of twisted localization algebra, we have the following commuting diagram:
\[\begin{tikzcd}
{K_*\Big(C^*_L(\widetilde{P_{d,n}}(X),\mathcal{A})^{\Gamma}\Big)} \arrow[r, "\pi_*"] \arrow[d, "\Phi_*"]\arrow[ddddd,"\IQ-\text{injection}" description, bend right=90] & {K_*\Big(C^*_L(P_{d}(X),\mathcal{A})^{\Gamma}\Big)} \arrow[d, "\Phi_*"]\\
{K_*\Big(C^*_L(\widetilde{P_{d,n}}(X) \times\Omega_X,\mathcal{A} \otimes (C(\Omega_X)\rtimes_r \Gamma))^{\Gamma}\Big)} \arrow[r, "\pi_*"] \arrow[d, "\text{collapsing}"]&{K_*\Big(C^*_L(P_{d}(X)\times \Omega_X,\mathcal{A}\otimes (C(\Omega_X)\rtimes_r\Gamma))^{\Gamma}\Big)} \arrow[ld, "\text{ordering}"', bend left=12]\arrow[ldd, bend left=21,"\Ga-\text{homotopy}"] \arrow[ldddd, bend left=38] \\
{K_*\Big(C^*_L(\widetilde{P_{d,n}}(X),\mathcal{A}\otimes(C(\Omega_X)\rtimes_r \Gamma))^{\Gamma}\Big)} \arrow[d, "{d<d',n<n'}"]&\\
{K_*\Big(C^*_L(\widetilde{P_{d',n'}}(X),\mathcal{A}\otimes(C(\Omega_X)\rtimes_r \Gamma))^{\Gamma}\Big)} \arrow[d, "\cong"', "\text{K\"unneth formula}"]&\\
{\bigoplus\limits_{i+j=*} K_{i}\Big( C_L^*\big( \widetilde{P_{d',n'}}(X), \; \mathcal{A} \big)^\Gamma \Big) \otimes_{\mathbb{Z}} K_{j} \Big( C(\Omega_X) \rtimes_r \Gamma  \Big)} \arrow[d, "\text{trace $\tau_{\Omega_X}$ on $(C(\Omega_X)\rtimes_r\Gamma)$}"]&\\
{K_{*}\Big( C_L^*\big( \widetilde{P_{d',n'}}(X), \; \A \big)^\Gamma \Big) \otimes_{\mathbb{Z}} \mathbb{R}}&
\end{tikzcd}\]
Here are some explanations for the diagram. The notation $\A$ in the diagram is short for $\A(M^{[0,1]})$ and the map $\Phi_*$ is defined as in Section \ref{KK-product for localization alg}. The first four rows are induced by the diagram \eqref{collapsing-ordering}, as we discussed above, this diagram indeed commutes. The fifth row is given by the K\"unneth formula for twisted localization algebra as in Theorem \ref{Kunneth formula}. The trace $\tau_{\Omega_X}$ is defined as in the last part of Section \ref{sec: A Concrete model for space with Property TAF} and acts on $K$-theory as in \cite[Definition 2.3.17]{HIT2020}.

Tracing along the left column, the map is given by $[g]\mapsto [i_*(g)]\ox 1$. Indeed, the first map is essentially given by $[g]\mapsto[g\ox 1_{C(\Omega_X)\rtimes_r\Ga}]$ and then absorbed by an ample representation. Since the composition of the collapsing map and the inclusion is $\Ga$-homotopic to the canonical inclusion map, thus the composition of the first four maps is given by $[g]\mapsto [i_*(g)\ox 1_{C(\Omega_X)\rtimes_r\Ga}]$, where $i_*$ is induced by the canonical inclusion map
$$i:\widetilde{P_{d,n}}(X)\to \widetilde{P_{d',n'}}(X).$$
By definition, $[i_*(g)\ox 1_{C(\Omega_X)\rtimes_r\Ga}]$ corresponds to $[i_*(g)]\ox[1_{C(\Omega_X)\rtimes_r\Ga}]$ under the K\"unneth formula, see Theorem \ref{Kunneth formula}. Since $\tau_{\Omega_X}\left([1_{C(\Omega_X)\rtimes_r\Ga}]\right)=1$, thus we conclude that the left column is given by $[g]\mapsto [i_*(g)]\ox 1$. Passing $d,n$ to infinity, we have that the left column is actually given by the identity map after tensoring $\IQ$. This means that the first row is a rational injection.
\end{proof}

By using a similar argument, one also has the following corollary.

\begin{Cor}\label{inj between classifying}
For any proper $\Ga$-space $X$, the homomorphism
$$\pi_*:\lim_{d,n\to\infty}K_*(C^*_L(\MRdnX)^{\Ga})\to \lim_{d\to\infty}K_*(C^*_L(P_d(X))^{\Ga})$$
is rational injective. If we moreover assume that $X$ is $\Ga$-compact, then the result above shows that the canonical homomorphism
$$RK^{\Ga}_*(\E\Ga)\ox\IQ\to RK^{\Ga}_*(\EG)\ox\IQ$$
is injective.\qed
\end{Cor}

The $\Ga$-compact case of Corollary \ref{inj between classifying} is proved in different ways from \cite{BC1988} and \cite{AAS2020}.

\begin{proof}[Proof of Theorem \ref{main result with torsion}]
Trace along the leftmost column and then the bottom row in the diagram \eqref{main diam}. In Section \ref{sec: A deformation trick}, we have shown that the composition of map (5) and map (8) is injective as $d,n$ tend to infinity, since $\MRdnX$ admits a uniformly locally finite $\F$-cover $\U$ such that each element in $\U$ is contractible for each $d\geq 0$ and $n\in\IN$. By Theorem \ref{twist assembly}, we have that the map (11)
is an isomorphism as $d$ tends to infinity, since $\Ga\car X$ has equivariant bounded geometry and $f:X\to M$ is an equivariant coarse embedding. By Lemma \ref{pi is injective}, we have that map (12) is rational injective as $d,n$ tend to infinity, then we have the composition of maps (5), (8), (12), and (11) is rational injective. This implies that the composition of the maps in the top row is rational injective, as desired.
\end{proof}

\begin{appendices}

\section{A K\"unneth formula for $K$-theory of twisted localization algebras}\label{Appendix A}

In this appendix, we provide a complete proof of Theorem \ref{Kunneth}, which can be viewed as a K\"unneth formula for $K$-theory of twisted localization algebras. The reader is referred to \cite{Duke1987} for analogous discussions in the framework of Kasparov's $KK$-theory.

To this end, we shall first introduce a new description of twisted localization algebras and use it to prove the $K$-theory of twisted localization algebra only depends on the topological structure (instead of the metric structure) of the base space. The new version of twisted localization algebras will provide better functorial properties that we need. The reader is referred to \cite{DWW2018} for some relevant discussion.

\subsection*{Localization algebras for topological spaces}

Throughout this section, let $X$ be a secondly countable, proper metric space with \emph{trivial group action}. Let $(N, M^{[0,1]},f_0)$ be the coefficient system as in Section 4. The twisted localization algebra is denoted by $C^*_L(X,\A(M^{[0,1]}))$ in this section (\emph{we omit $\Ga$ and $(0)$ on the top corner since the group action and the embedding are trivial}). Recall that $E_{\A(M^{[0,1]})}=\ell^2(N)\ox\H_0\ox\A(M^{[0,1]})$. An element $T\in\L(E_{\A(M^{[0,1]})})$ can be viewed as a $N$-by-$N$ matrix, where the $(x,y)$-entry of $T$, denoted by $T(x,y)$, is defined to be $\delta_xT\delta_y$. For trivial action case, we can omit $\ell^2(\Ga)$ in $E_{\A(M^{[0,1]})}$.

\begin{Def}\label{commutator}
Define $\sC_L[X,\A(M^{[0,1]})]$ to be the set of all bounded and uniformly continuous functions
$$g:[0,\infty)\to\L(E_{\A(M^{[0,1]})})$$
satisfying\begin{itemize}
\item[(1)] for each $t\in[0,\infty)$ and $x,y\in N$, $(g(t))(x,y)$ is an element of $\K(\H_0)\ox\A(M^{[0,1]})$ and $g(t)$ satisfies condition (3) in Definition \ref{twisted Roe small} for each $t\in\IR_+$;
\item[(2)] for any $a\in C_0(X)$, the commutator $[a,g(t)]$ tends to $0$ as $t$ tends to infinity, where $a\in C_0(X)$ is viewed as a self-adjoint operator on $\L(E_{\A(M^{[0,1]})})$;
\item[(3)] there exists $R>0$ such that
$$\supp((g(t))(x,y))\subseteq B(m,R)\subseteq M^{[0,1]}\times\IR_+$$
for all $x,y\in N$, where $m\in M\subseteq M^{[0,1]}$ is the base point of $M$, viewed as a point of $M^{[0,1]}$.
\end{itemize}
Define $\sC^*_L(X,\A(M^{[0,1]}))$ to be the completion of $\sC_L[X,\A(M^{[0,1]})]$ under the norm
$$\|g\|=\sup_{t\in[0,\infty)}\|g(t)\|.$$
\end{Def}

Notice that we require the commutator $[g(t), a]$ to tend to $0$ instead of $\prop(f(t))$ tends to $0$. For $a\in C_0(X)$ and $r>0$, the $r$-oscillation of $a$ is defined to be
$$\omega_r(a)=\sup_{d(x,y)\leq r}|a(x)-a(y)|.$$
Since $f$ is uniformly continuous, $\omega_r(a)\to 0$ as $r\to 0$. Notice that
\begin{equation}\label{commutator}\|[T,a]\|\leq 8\omega_{Prop(T)}(a)\|T\|,\end{equation}
see \cite[Lemma 6.12]{HIT2020} for a proof. Then we have that $C^*_L(X,\A(M^{[0,1]}))\subseteq \sC^*_L(X,\A(M^{[0,1]}))$. Notice that our new version of localization algebras is larger and does not depend on the metric structure of $X$, but only depends on the topological structure. When we compare these two algebras on $K$-theory level, we are only interested in the asymptotic behavior of the functions, while the precise values of these functions at any finite $t\in\IR_+$ are somewhat redundant information. Therefore it is useful to develop a formalism that does not depend on such redundant data.

\begin{Def}[Annihilator ideals]
The annihilator ideals in $\sC^*_L(X,\A(M^{[0,1]}))$ and $C^*_L(X,\A(M^{[0,1]}))$ are, respectively, defined to be
$$\mathrm{Ann}(\sC^*_L)=\left\{g\in\sC^*_L(X,\A(M^{[0,1]}))\mid \forall a\in C_0(X),\, a\cdot g\in C_0(\IR_+,\L(E_{\A(M^{[0,1]})}))\right\};$$
$$\mathrm{Ann}(C^*_L)=C^*_L(X,\A(M^{[0,1]}))\cap \mathrm{Ann}(\sC^*_L).$$
\end{Def}

To simplify notations, we will denote the quotient algebras by
$$\sC^*_{L,Q}=\frac{\sC^*_L(X,\A(M^{[0,1]}))}{\mathrm{Ann}(\sC^*_L)}\quad\text{ and }\quad C^*_{L,Q}=\frac{C^*_L(X,\A(M^{[0,1]}))}{\mathrm{Ann}(C^*_L)}.$$
Then the canonical inclusion induces a $C^*$-homomorphism on the quotients algebras
$$i_Q:C^*_{L,Q}\to \sC^*_{L,Q}.$$

\begin{Lem}\label{quo-iso}
The canonical inclusion induces an isomorphism between the quotient algebras, i.e.,
$$i_Q:C^*_{L,Q}\xrightarrow{\cong} \sC^*_{L,Q}.$$
\end{Lem}

\begin{proof}
It suffices to show $i_Q$ is surjective, i.e., for any $g\in\sC_L[X,\A(M^{[0,1]})]$, there exists $\wh g\in\IC_L[X,\A(M^{[0,1]})]$ such that $g-\wh g\in \mathrm{Ann}(\sC^*_L)$.

Fix a base point $x_0\in X$. For each $n\in\IN^*$, set $\U_n$ to be a locally finite open cover of $X$ such that $\diam(U)\leq 1/n$ for any $U\in\U_n$ and $\{\phi_U^n\}_{U\in\U_n}$ to be an $\ell^2$-partition of unity associated to $\U_n$. Define
$$\U_{n,x_0}=\{U\in\U_n\mid U\cap B(x_0,n)\neq\emptyset\}.$$
Since $\U_n$ is locally finite, there are only finitely many elements in $\U_{n,x_0}$. Thus there exists $T_n>0$ such that
$$\|[g(t),\phi^n_U]\|\leq\frac{1}{n\cdot\#\U_{n,x_0}}$$
for any $U\in\U_{n,x_0}$ and $t\geq T_n$. Without loss of generality, assume that $T_{n+1}>T_n+1$ for each $n\in\IN$.

Define $\Phi_n(T)=\sum_{U\in\U_n}\phi^n_UT\phi_U^n$ for each $T\in\L(E_{\A(M^{[0,1]})})$ with the convergence in the strong operator topology. It is clear that $\Phi_n$ is an unital positive linear map for each $n\in\IN$. Thus $\Phi_n$ is bounded and $\|\Phi_n\|\leq 1$ for each $n$. For each $t\in[T_n,T_{n+1}]$, we define
\[\wh g(t)=\frac{T_{n+1}-t}{T_{n+1}-T_n}\Phi_n(g(t))+\frac{t-T_n}{T_{n+1}-T_n}\Phi_{n+1}(g(t)).\]
It is clear that $\wh g$ is bounded, $\wh g(t)\in\IC[X,\A(M^{[0,1]})]$ for each $t\in\IR_+$ and $\prop(\wh g(t))\to0$ as $t\to\infty$. By using a similar argument as \cite[Lemma 2.2]{QR2010}, one can also show that $\wh g$ is uniformly continuous. Then we conclude that $\wh g$ is an element in $\IC_L[X,\A(M^{[0,1]})]$.

It remains to show that $g-\wh g\in \mathrm{Ann}(\sC^*_L)$. For any $a\in C_c(X)$, there exists $N>0$ such that $\supp(a)\subseteq B(x_0,N)$. Then for each $n>N$,
\[\begin{split}
\|a(g(T_n)-\wh g(T_n))\|&=\Big\|a\cdot g(T_n)-a\sum_{U\in\U_{n,x_0}}\phi^n_U g(T_n)\phi^n_U\Big\|\\
&\leq \Big\|a\cdot g(T_n)-a\sum_{U\in\U_{n,x_0}}(\phi^n_U)^2\cdot g(T_n)\Big\|+\Big\|a\sum_{U\in\U_{n,x_0}}\phi^n_U\cdot[g(T_n),\phi^n_U]\Big\|\\
&\leq 0+\|a\|\cdot\frac 1n.
\end{split}\]
If $t\in[T_n,T_{n+1}]$ for $n\geq N$, then we have that
\[\begin{split}
\|a(g(t)-\wh g(t))\|&= \Big\|a\cdot g(t)-a\Phi_n(g(t))-a\frac{t-T_n}{T_{n+1}-T_n}\Phi_{n+1}(g(t))\Big\|\\
&\leq \frac{T_{n+1}-t}{T_{n+1}-T_n}\cdot\Big\|a\cdot g(t)-a\cdot\Phi_n(g(t))\Big\|+\frac{t-T_n}{T_{n+1}-T_n}\cdot\Big\|a\cdot g(t)-a\cdot\Phi_{n+1}(g(t))\Big\|\\
&\leq \|a\|\cdot\frac 1n\to 0\quad\text{ as }n\to\infty.
\end{split}\]
This shows that $g-\wh g\in \mathrm{Ann}(\sC^*_L)$.
\end{proof}

\begin{Lem}\label{ann-zero}
The quotient maps induce isomorphisms on $K$-theory:
$$K_*(C^*_L(X,\A(M^{[0,1]})))\xrightarrow{\cong} K_*(C^*_{L,Q})\quad \text{and}\quad K_*(\sC^*_L(X,\A(M^{[0,1]})))\xrightarrow{\cong}K_*(\sC^*_{L,Q}).$$
\end{Lem}

\begin{proof}
It suffices to show that $K_*(\mathrm{Ann}(C^*_L))$ and $K_*(\mathrm{Ann}(\sC^*_L))$ are trivial groups. It follows using an Eilenberg swindle argument similar to \cite[Lemma 6.4.11]{HIT2020}, which we leave to the reader.
\end{proof}

\begin{Thm}\label{comm-prop}
The $K$-theortic homomorphism induced by the canonical inclusion
$$i_*:K_*(C^*_L(X,\A(M^{[0,1]})))\to K_*(\sC^*_L(X,\A(M^{[0,1]})))$$
is an isomorphism.
\end{Thm}

\begin{proof}
Notice that we have the following commuting diagram:
$$\begin{tikzcd}
0 \arrow[r] & \mathrm{Ann}(C^*_L) \arrow[r] \arrow[d] & C^*_L(X,\A(M^{[0,1]})) \arrow[r] \arrow[d] & C^*_{L,Q} \arrow[d, "\cong"] \arrow[r] & 0 \\
0 \arrow[r] & \mathrm{Ann}(\sC^*_L) \arrow[r]           & \sC^*_L(X,\A(M^{[0,1]})) \arrow[r]           & \sC^*_{L,Q} \arrow[r]           & 0
\end{tikzcd}$$
It induces a commuting diagram on $K$-theory. Then the theorem follows directly from Lemma \ref{quo-iso} and Lemma \ref{ann-zero}.
\end{proof}

From the definition, one can see that $\sC^*_L(X,\A(M^{[0,1]}))$ only depend on $C_0(X)$, i.e., the topological structure of $X$. For metric space $(X,d)$, let $d'$ be another metric on $X$ such that $(X,d)$ is homeomorphic to $(X,d')$. Then for both metric spaces, the $K$-theory of their twisted localization algebras associated to $(N, M^{[0,1]},f_0)$ are isomorphic by Theorem \ref{comm-prop}. Similarly, we also have the following result:

\begin{Thm}\label{non-coe comm-prop}
Let $X$ be a secondly countable locally compact Hausdorff space, $H_X$ an ample $X$-module. Define $\sC^*_L(X)$ to be the $C^*$-algebra generated by all bounded and uniformly continuous functions $g:\IR_+\to B(H_X)$ such that $g(t)$ is locally compact for each $t\in\IR_+$ and
$$[g(t),a]\to 0\quad\text{as}\quad t\to\infty$$
for any $a\in C_0(X)$. Then $K_*(\sC^*_L(X))$ is naturally isomorphic to $KK_*(C_0(X),\IC)$, the $K$-homology group of $X$.
\end{Thm}

\begin{proof}
Equip $X$ with a metric that induces the topology of $X$. Then one can define $C^*_L(X)$ under this metric. By \cite[Theorem 3.5]{QR2010}, one has that $KK_*(C_0(X),\IC)$ is isomorphic to $K_*(C^*_L(X))$. By using a similar argument as Theorem \ref{comm-prop}, one has that $K_*(C^*_L(X))$ is isomorphic to $K_*(\sC^*_L(X))$ and this isomorphism does not depend on the choice of metric.
\end{proof}

The $K$-theory of the localization algebra provides an analytic-geometric model for the $K$-homology of $X$ to study the Baum-Connes conjecture. However, it has a disadvantage in that the functorial properties of localization algebras rely on the metric structures of spaces and Lipschitz maps between spaces (or continuous coarse maps) which is not as flexible as Kasparov's $K$-homology. Thanks to Theorem \ref{comm-prop}, for any second countable locally compact Hausdorff spaces $X, Y$, and a proper continuous map $f: X\to Y$, we can equip $X$ and $Y$ with bounded metrics which induce their topologies, and thus $f$ is clearly a coarse map in this view. Moreover, in this point of view, the strong Lipschitz condition can be dropped for the homotopy equivalence of $K$-theory of localization algebras. For the non-twisted case, it follows directly from Theorem \ref{non-coe comm-prop} and the homotopy invariance of $KK$-theory. For the twisted case, one can follow an outline of the proof of \cite[Proposition 6.4.14]{HIT2020} to show $K_*(\sC^*_L(X\times\IR_+,\A(M^{[0,1]})))=0$ for any $X$, then the homotopy invariance follows from a Mayer-Vietoris argument as in \cite[Theorem 6.4.16]{HIT2020}.

Before we can prove the K\"unneth formula, we still need to prove an excision theorem for twisted localization algebras for a special case. Let $X$ be a secondly countable \emph{compact} Hausdorff space. Let $K\subseteq X$ be a closed subset, denote $U=X\backslash K$. Equip $X$ with a metric that induces the topology of $X$. Since $X$ is compact, the metric on $X$ must be bounded. Let $(N,M^{[0,1]},f_0)$ be a coefficient system for $X$. Then
$$(N\cap K,M^{[0,1]},f_0|_K)\quad \text{and}\quad (N\cap U,M^{[0,1]},f_0|_U)$$
form the coefficient system for $K$ and $U$ respectively. We can then define the twisted localization algebras for $K$ and $U$. The norms of these two algebras are given by the representation on the Hilbert modules
$$E_K=\ell^2(N\cap K)\ox\H_0\ox\A(M^{[0,1]})\quad\text{ and }\quad E_{U}=\ell^2(N\cap U)\ox\H_0\ox\A(M^{[0,1]}).$$
To simplify the notation, we denote $E_X=\ell^2(N)\ox\H_0\ox\A(M^{[0,1]})$ temporarily. Then $E_X=E_K\oplus E_{U}$. Thus norms of the twisted algebras of $K$ and $U$ are given from the canonical representations on $E_X$ via the inclusions $E_{K}\to E_X$ and $E_U\to E_X$.

\begin{Def}
Let $X$ and $K$ be as above. Define $C^*_{L}(X,K,\A(M^{[0,1]}))$ to be the closed subalgebra of $C^*_{L}(X,\A(M^{[0,1]}))$ generated by all elements $g$ such that for any $\varepsilon>0$, $\exists T>0$ such that $(g(t))(x,y)=0$ for all $t>T$ and $x\notin B(K,\varepsilon)$, where $B(K,\varepsilon)$ is the $\varepsilon$-neighborhood of $K$.
\end{Def}

Using an argument similar to \cite[Lemma 3.10]{Yu1997}, one can show that $K_*(C^*_{L}(X,K,\A(M^{[0,1]})))$ is isomorphic to $K_*(C^*_{L}(K,\A(M^{[0,1]})))$. It is also clear that $C^*_{L}(X,K,\A(M^{[0,1]}))$ is a closed two-sided ideal of $C^*_L(X,\A(M^{[0,1]}))$.

\begin{Lem}\label{pi_Q}
Let $\pi:C^*_L(X,\A(M^{[0,1]}))\to C^*_L(U,\A(M^{[0,1]}))$ be the linear map defined by
$$\pi(g)=\chi_Ug\chi_U,$$
where $X$ is the characteristic function on $X$. Then $\pi$ induces an $C^*$-isomorphism on the quotient algebras:
$$\pi_Q:\frac{C^*_L(X,\A(M^{[0,1]}))}{C^*_L(X,K,\A(M^{[0,1]}))}\to\frac{C^*_L(U,\A(M^{[0,1]}))}{\mathrm{Ann}(C^*_L(U,\A(M^{[0,1]})))}\quad [g]\mapsto [\pi(g)],$$
where the ideal on the right side is the annihilator ideal of $X$, i.e., any element of the ideal will become a $C_0$ function after multiplying with a function in $C_0(U)$.
\end{Lem}

\begin{proof}
To see $\pi_Q$ is well-defined, we shall prove that for any $g\in C^*_{L}(X,K,\A(M^{[0,1]}))$, $\pi(g)$ must be an element of the annihilator ideal. Indeed, for any $a\in C_c(U)$, $\supp(a)$ must be a compact subsets of $U$. Since $\supp(a)\cap K=\emptyset$, then $d(\supp(a),K)=\varepsilon>0$. By definition, there exists $T>0$ such that $g(t)(x,y)=0$ for all $x\in\supp(a)$ and $t>T$. This shows that $\pi(g)$ is in the annihilator ideal and $\pi_Q$ is a well-defined set-theoretic map.

By using the inequality \eqref{commutator}, one can see clearly that $[\chi_U,g]$ is in the annihilator ideal of $U$ for any $g\in C^*_L(X,\A(M^{[0,1]}))$. This shows that $\pi_Q$ is a well-defined $C^*$-homomorphism. Since $\pi$ is surjective, then $\pi_Q$ is clearly a surjection. It suffices to show that for any $g\in C^*_L(X,\A(M^{[0,1]}))$, if $\pi(g)\in \mathrm{Ann}(C^*_L(U,\A(M^{[0,1]})))$, then we must have $g\in C^*_{L,K}(X,\A(M^{[0,1]}))$.

To see this, for any $\varepsilon>0$, set $C=\overline{X\backslash B(K,\varepsilon)}$. Since $X$ is compact, $C$ must be a compact subspace of $U$. Let $a\in C_0(U)$ such that $a(x)=1$ for all $x\in C$. Since $\pi(g)\in \mathrm{Ann}(C^*_L(U,\A(M^{[0,1]})))$, without loss of generality, there exists $T>0$ such that $a\pi(g(t))=0$ for all $t>T$. This means that $a\cdot g(t)\cdot\chi_U=0$, combining $\prop(g(t))\to 0$, we conclude that $(g(t))(x,y)=0$ for all $x\in C$. We therefore have that $g\in C^*_{L,K}(X,\A(M^{[0,1]}))$. This completes the proof.
\end{proof}

We should mention that our proof only works for the case when $X$ is compact. For the case when $X$ is locally compact but non-compact, the reader is referred to Section 6.4 in \cite{HIT2020} or \cite[Proposition 4.6]{GLWZ2022} for some introductions. One can also find a general proof by using the functorial property in \cite[Proposition B.2.3]{HIT2020}

\begin{Thm}\label{excision}
There is a six-term exact sequence on $K$-theory:
$$\begin{tikzcd}
K_0(C^*_L(K,\A(M^{[0,1]}))) \arrow[r] & K_0(C^*_L(X,\A(M^{[0,1]}))) \arrow[r] & K_0(C^*_L(U,\A(M^{[0,1]}))) \arrow[d] \\
K_1(C^*_L(U,\A(M^{[0,1]}))) \arrow[u] & K_1(C^*_L(X,\A(M^{[0,1]}))) \arrow[l] & K_1(C^*_L(K,\A(M^{[0,1]}))) \arrow[l]
\end{tikzcd}$$
associated to the $C^*$-algebraic exact sequence
$$0\to C_0(U)\to C(X)\to C(K)\to 0.$$
\end{Thm}

\begin{proof}
Consider the following $C^*$-algebraic exact sequence:
$$0\to C^*_L(X,K,\A(M^{[0,1]}))\to C^*_L(X,\A(M^{[0,1]}))\to\frac{C^*_L(X,\A(M^{[0,1]}))}{C^*_L(X,K,\A(M^{[0,1]}))}\to 0.$$
It induces a six-term exact sequence on $K$-theory. Then the theorem holds certainly from Lemma \ref{ann-zero} and Lemma \ref{pi_Q}.
\end{proof}

\subsection*{The K\"unneth formula}

\begin{Thm}\label{kunneth}
Let $X$ be a locally compact, Hausdorff, second countable space with a \emph{trivial group action}. Let $(N,M^{[0,1]},f_0)$ be a coefficient system as in Section 4, where $f_0:N\to M^{[0,1]}$ is a constant map to the base point $m$ of $M$ viewed as a subspace of $M^{[0,1]}$. Then
$$K_i(C^*_L(X,\A(M^{[0,1]}))^{(0)})\ox\IQ\cong \bigoplus_{j\in\IZ_2\IZ}K_j(X)\ox K_{i-j}(\A(M^{[0,1]}))\ox\IQ,$$
where $C^*_L(X,\A(M^{[0,1]}))^{(0)}$ is the twisted localization algebra defined as in Definition \ref{twisted localization small} with trivial group action.
\end{Thm}

\begin{proof}
From above, we have shown that the localization algebra does not depend on the choice of metric which induces the topology of $X$. Thus we can view $X$ as a metric space induced with a bounded metric, thus the localization algebras are well-defined.

There exists a canonical inclusion
$$\IC_L[X]\ox\A(M^{[0,1]})\to\IC_L[X,\A(M^{[0,1]})]^{(0)},$$
since any element $g\ox a\in\IC_L[X]\ox\A(M^{[0,1]})$ can be naturally seen as an element of $\IC_L[X,\A(M^{[0,1]})]^{(0)}$. To see this, for any $\varepsilon>0$, there exists $R>0$, and $\wh a\in\A(M^{[0,1]})$ such that $\|a-\wh a\|\leq \varepsilon$ and $\supp(\wh a)\subseteq B(m,R)\subseteq M^{[0,1]}\times\IR_+$. Then $g\ox\wh a$ defines a well-defined element in $\IC_L[X,\A(M^{[0,1]})]^{(0)}$, where
$$((\wh g\ox\wh a)(t))(x,y)=(\wh g(t))(x,y)\ox\wh a\in\K(\H)\ox\A(M^{[0,1]}),$$
and $\H=\ell^2(\Ga)\ox\H_0$ as in Remark \ref{Roe finite}. This map induces a $C^*$-homomorphism
$$C^*_L(X)\ox\A(M^{[0,1]})\to C^*_L(X,\A(M^{[0,1]}))^{(0)},$$
Then we have the following homomorphism on $K$-theory:
$$\bigoplus_{j\in\IZ_2\IZ}K_j(X)\ox K_{i-j}(\A(M^{[0,1]}))\to K_i(C^*_L(X)\ox\A(M^{[0,1]}))\to K_i(C^*_L(X,\A(M^{[0,1]}))^{(0)}).$$
The first map is given by the Kasparov $KK$-product (see \cite{Kas88} or \cite[Section 2.10]{HIT2020} for an introduction to external products in $K$-theroy), and the second map is induced by the canonical inclusion. We shall then prove this map is an isomorphism after tensoring both sides by $\IQ$.

First, we consider the case when $X$ is a compact space. This case is similar to \cite[Lemma 2.4]{GWY2018}. When $X$ is a point, the theorem holds clearly since $K_0(X)\cong\IZ$, $K_1(X)\cong 0$ and
$$K_*(C^*_L(X,\A(M^{[0,1]}))^{(0)})\cong K_*(\A(M^{[0,1]})).$$
By Theorem \ref{comm-prop} and Theorem \ref{non-coe comm-prop}, the localization algebras are homotopy invariant on the $K$-theory level. Thus the theorem holds for the case when $X$ is contractible. Since both functors
$$X\mapsto K_*(X)\ox K_*(\A(M^{[0,1]}))\ox\IQ\quad\text{and}\quad X\mapsto K_*(C^*_L(X,\A(M^{[0,1]}))^{(0)})\ox\IQ$$
satisfy the Mayer-Vietoris theorem. We should mention that we take tensor products with $\IQ$ to make sure the Mayer-Vietoris sequence is still exact after taking tensor $K_*(\A(M^{[0,1]}))$. By using a cutting and pasting argument and the Five Lemma, the theorem holds for the case when $X$ is a compact CW-complex. Since any compact Hausdorff space can be written as an inverse limit of finite CW complexes (see \cite[Theorem X.10.1]{fundaofalgtop}), the theorem holds for all compact spaces.

For the case when $X$ is locally compact, denote $X^+=X\cup\infty$ the one-point compactification. Equip $X^+$ with a metric. Then $(N\cup\{\infty\},M^{[0,1]},\widetilde{f_0})$ forms a coefficient system for $X^+$, where $\widetilde{f_0}$ is the constant map onto the base point $m\in M$. Then $\{\infty\}$ is a closed subspace of $X^+$ and $X=X^+\backslash\{\infty\}$. By Theorem \ref{excision}, we have a six-term exact sequence of $K$-theory of twisted localization algebras associated to the $C^*$-exact sequence:
$$0\to C_0(X)\to C(X^+)\to \IC\to 0.$$
Since the theorem holds for both $X^+$ and $\{\infty\}$, the theorem also holds for $X$ by using the Five Lemma.
\end{proof}

\end{appendices}

\bibliographystyle{alpha}
\bibliography{ref}
\end{document}